\documentclass{amsart}

\usepackage{amssymb}
\usepackage[all]{xy}
\usepackage{hyperref}

\usepackage{enumitem}   

\setlist[enumerate]{itemsep=.2em,topsep=.2em,leftmargin=1.25em,itemindent=2.0em}

\newtheorem{thm}{Theorem}

\newtheorem{lem}[thm]{Lemma}
\newtheorem{cor}[thm]{Corollary}

\newtheorem{prop}[thm]{Proposition}

\newtheorem{conj}[thm]{Conjecture}
   
\theoremstyle{definition}
\newtheorem{defn}[thm]{Definition}

\newtheorem{say}[thm]{}
\newtheorem{exmp}[thm]{Example}

\newtheorem{rem}[thm]{Remark}          

\newtheorem*{ack}{Acknowledgments}      
\newtheorem{notation}[thm]{Notation}   
  
\newtheorem{defn-thm}[thm]{Definition--Theorem}  
\newtheorem{defn-lem}[thm]{Definition--Lemma}  

\theoremstyle{remark}

\let \cedilla =\c
\renewcommand{\c}[0]{{\mathbb C}}  

\renewcommand{\o}[0]{{\mathcal O}} 
\newcommand{\z}[0]{{\mathbb Z}}
\newcommand{\n}[0]{{\mathbb N}}
\renewcommand{\r}[0]{{\mathbb R}} 

\renewcommand{\a}[0]{{\mathbb A}}

\newcommand{\p}[0]{{\mathbb P}}
\newcommand{\f}[0]{{\mathbb F}}
\newcommand{\q}[0]{{\mathbb Q}}
\newcommand{\map}[0]{\dasharrow}
\newcommand{\qtq}[1]{\quad\mbox{#1}\quad}
\newcommand{\spec}[0]{\operatorname{Spec}}
\newcommand{\pic}[0]{\operatorname{Pic}}
\newcommand{\pics}[0]{\operatorname{\mathbf{Pic}}}

\newcommand{\pico}[0]{\operatorname{\mathbf{Pic}}^{\circ}}

\newcommand{\alb}[0]{\operatorname{\mathbf{Alb}}}

\newcommand{\ns}[0]{\operatorname{NS}}

\newcommand{\rank}[0]{\operatorname{rank}}
\newcommand{\qrank}[0]{\operatorname{rank}_{\q}}

\newcommand{\supp}[0]{\operatorname{Supp}}    
\newcommand{\red}[0]{\operatorname{red}}    
\newcommand{\codim}[0]{\operatorname{codim}}    
\newcommand{\im}[0]{\operatorname{im}}

\newcommand{\coker}[0]{\operatorname{coker}}

\newcommand{\bs}[0]{\operatorname{Bs}}

\newcommand{\sing}[0]{\operatorname{Sing}}

\newcommand{\res}[0]{\operatorname{\mathcal R}} 
 
\newcommand{\jacs}[0]{\operatorname{\mathbf{Jac}}} 
\newcommand{\jac}[0]{\operatorname{Jac}}

\newcommand{\chr}[0]{\operatorname{char}}

\newcommand{\cl}[0]{\operatorname{Cl}}

\newcommand{\clo}[0]{\operatorname{Cl^{\circ}}}
\newcommand{\clos}[0]{\operatorname{\mathbf{Cl}^{\circ}}}
\newcommand{\cls}[0]{\operatorname{\mathbf{Cl}}}

\newcommand{\len}[0]{\operatorname{length}}

\newcommand{\onto}[0]{\twoheadrightarrow}

\newcommand{\simq}[0]{\sim_{\q}}

\newcommand{\tsum}[0]{\textstyle{\sum}} 
\newcommand{\tr}[0]{\operatorname{tr}}

\newcommand{\wdiv}[0]{\operatorname{WDiv}} 
 
\newcommand{\pdiv}[0]{\operatorname{PDiv}}

\newcommand{\mg}[0]{{\mathbb G}_m}

\def\into{\DOTSB\lhook\joinrel\to}

\def\loccoh#1.#2.#3.#4.{H^{#1}_{#2}(#3,#4)}

\DeclareMathAlphabet{\mathchanc}{OT1}{pzc}%
                                {m}{it}

\newcommand{\gm}[0]{{\mathbb G}_m}
\newcommand{\ga}[0]{{\mathbb G}_a}

\usepackage[all]{xy}\xyoption{dvips}

\newcommand{\tprod}[0]{\textstyle{\prod}} 

\newcommand{\norm}[0]{\operatorname{norm}}

\newcommand{\albm}[0]{\operatorname{alb}}

\newcommand{\BH}[0]{\operatorname{BH}}

\newcommand{\sims}[0]{\sim_{\rm s}}
\newcommand{\simsa}[0]{\sim_{\rm sa}}

\begin{document}
\bibliographystyle{amsalpha}


 \title{What determines a variety?}
\author{J\'anos Koll\'ar}
\begin{abstract}
We show that, over a field of characteristic 0,  a normal, projective variety of dimension at least 4
 is uniquely determined by its underlying topological space. 
The proof builds on previous work of Lieblich and Olsson. 
\end{abstract}

 \maketitle

\tableofcontents

\section{Introduction}

By definition, a  scheme $X$ is a topological space---which we denote by $|X|$---and a sheaf of rings  $\o_X$  on the open subsets of $|X|$. 

As a continuation of \cite{k-lo, k-lo-2},
we study the following  natural but seldom considered questions. 
\begin{itemize}
\item How to read off properties of $X$ from $|X|$?
\item Does $|X|$ alone determine $X$?
\end{itemize}

In some cases, the answer is clearly negative.

\begin{itemize}
\item  If $C$ is an irreducible curve over an infinite field $k$,
then $|C|$ has 1 generic point and $|k|$ closed points, so $|C|$ depends only on the cardinality of $k$. 
\item For any variety $X$, the seminormalization  map $X^{\rm sn}\to X$
is a homeomorphism, so we should clearly restrict attention to
seminormal varieties. 
The same applies to the weak normalization in positive characeristic.
\item Let $K/L$ be a finite field extension and $X$ an irreducible,  normal $K$-variety. We can also view it  as an $L$-variety  $X_L$, and then $|X|=|X_L|$. Note that  $X_L$ has at least $\deg (K/L)$ geometric irreducible components, we can thus avoid this problem by working only with  geometrically irreducible varieties.
\item   Purely inseparable morphisms in positive characteristic
give homeomorphisms between quite dissimilar varieties. 
\item A more unexpected example is constructed  in \cite{MR624904}:  $|\p^2_{\bar\f_p}|\sim |\p^2_{\bar\f_q}|$
for any 2 primes $p, q$. In general it seems that  surfaces over finite fields and their algebraic extensions give many similar examples (\ref{quasifinite.field.thm.cor}).
\item   Let $S, T$ be real algebraic surfaces birational to $\p^2_{\r}$
and $\Phi:S(\r)\sim T(\r)$ a homeomorphism  in the classical topology.
By \cite{k-mang},  $\Phi$  can be approximated by maps that are homeomorphisms  
both in the classical and the Zariski topology.
\end{itemize}

Despite these examples, it is  possible that the above are exceptional  instances and the answer is positive in almost all other cases.

In Sections~\ref{set.ci.prop.sec}--\ref{proj.sp.sec} 
we study curves $C\subset X$ such that every  finite subset
$Z\subset C$ is obtained (or almost obtained) as $C\cap D$ for some divisor $D\subset X$. 
It turns out that this depends mostly on the base field $k$.
Studying a similar question for the union of a curve and a divisor
  allows us to distinguish
projective spaces from all other varieties.

\begin{thm}\label{planes.in.Pn.thm.2.intro}
Let $L$ be a  field  of characteristic 0 and  $K$    an arbitrary  field.
Let $Y_L$ be a normal, projective, geometrically irreducible  variety of dimension $n\geq 2$ over  $L$  such that $|Y_L|$  is homeomorphic to
 $|\p^n_K|$.    Then   $K\cong L$ and 
 $Y_L\cong \p^n_L$.
\end{thm}

The linkage of divisors is investigated in 
Sections~\ref{sec.zero.set.sect}--\ref{top.pencils.sec}.  This leads to an  answer to the basic problem in characteristic 0, save in low dimensions.
(We comment on the dimension restrictions in (\ref{basic.plan.say}).)

\begin{thm}\label{home.iso.thm}
 Let $K, L$ be   fields of characteristic 0  and $X_K, Y_L$   normal, projective, geometrically irreducible  varieties 
over $K$  (resp.\ $L$). Let
$\Phi:|X_K|\sim |Y_L|$ be a homeomorphism.  Assume that 
\begin{enumerate}
\item either $\dim X\geq 4$,
\item or $\dim X\geq 3$ and  $K, L$ are   finitely generated  field extensions of $\q$.
\end{enumerate}
Then $\Phi$ is the composite of  a field isomorphism  $\phi:K\cong L$ and an algebraic isomorphism of $L$-varieties $X_L^{\phi}\cong Y_L$. 
\end{thm}

For completeness, let us state the most general form of (\ref{home.iso.thm})  that could be true.

\begin{conj}\label{home.iso.conj}
 Let $K, L$ be   fields    and $X_K, Y_L$   seminormal,  geometrically irreducible  varieties 
over $K$  (resp.\ $L$). Let
$\Phi:|X_K|\sim |Y_L|$ be a homeomorphism. 
Assume that $\chr L=0$  and $\dim X_L\geq 2$.

Then $\Phi$ is the composite of  a field isomorphism  $\phi:K\cong L$ and an algebraic isomorphism of $L$-varieties $X_L^{\phi}\cong Y_L$. 
\end{conj}

It is possible that  modifications of the current methods will handle the seminormal and  proper case. Our approach has many difficulties with open varieties; though some cases are treated in \cite{k-lo, k-lo-2}.  
\medskip

Ancillary results are collected in Sections~\ref{compl.int.sec}--\ref{spec.fields.sec}. These are mostly known but hard to find references  for the forms that we need. The notion of weakly Hilbertian fields introduced in Section~\ref{spec.fields.sec}
is new and may be of independent interest.

\begin{say}[Terminology and notation] For ease of reference, here is a list of new or  non-standard terminology and notation that we use.

{\it Albanese variety,}  classical  $\alb(X)$, \`a la Grothendieck  $\alb^{\rm gr}(X)$  (\ref{pic.alb.say}).

{\it Absolutely scip,}   variant of scip (\ref{ci.defn.3.abs}). 

{\it Ample-ci,} complete intersection of ample divisors   (\ref{ci.sci.say}). 

{\it Ample-sci,} set-theoretic complete intersection of ample divisors   (\ref{ci.sci.say}).

{\it Anti--Mordell-Weil field,}  (\ref{anti-mordell.defn}). 

{\it Bertini-Hilbert dimension of a field,} (\ref{ber-hil.dim.defn}).

$\BH(k)$,  the Bertini-Hilbert dimension; it is either 1 or 2  (\ref{ber-hil.dim.defn}). 

{\it Detects  linear similarity,}    (\ref{rest.linsim.say}).

{\it Field-locally thin subset,}   (\ref{thin.flt.defn}).

{\it Generically scip,} variant of scip   (\ref{ci.defn}). 

{\it Linked, $L$-linked,}  (\ref{t-link.gen.defn}).  

{\it Linking is free,}  (\ref{t-link.ZW.defn}). 

{\it Linking on  $W_2$  determines linking on  $W_1$,}
 (\ref{red.int.top.prop.defn}). 

{\it Linking is minimally restrictive,}   (\ref{red.int.top.prop.defn}). 

{\it Locally finite field,}  algebraic extension of a finite field  (\ref{loc.fin.defn}). 

{\it Mordell-Weil field,}  (\ref{m.w.f.defn}).

{\it Non-Cartier center,} (\ref{non.car.cent.defn}).  

{\it Scip,}  set-theoretic complete intersection property.  
\nopagebreak

\qquad For irreducible subsets  (\ref{ci.defn}) and for  reducible subsets  (\ref{ci.defn.red.3}). 

{\it Scip with finite defect,}  variant of scip  (\ref{ci.defn.3}). 

 
{\it Thin subset,}   (\ref{thin.flt.defn}). 

{\it Stably dense,}  (\ref{QD.defn}). 

{\it Topological pencil, t-pencil,}  (\ref{t.pencil.defn}). 

{\it Weakly Hilbertian field,}  (\ref{weak.hilb.say}).

 $\Gamma^B(X,L)$ , sections of (powers of) $L$ of that vanish exactly along $B$  (\ref{full.supp.sects.say}). 

$\res_V^B(X,L)$,   restriction of $\Gamma^B(X,L)$ to $V\subset X$  (\ref{ZW.H.link.R.defn}). 

$\res^A_k\gm$,   Weil restriction of the multiplicative group (\ref{gm.not.defn}). 

$\rho(X)$,  Picard number using Cartier divisors.

$\rho^{\rm cl}(X)$, Picard number (or class rank) using non-Cartier divisors (\ref{cl.norm.say}).

$\Sigma(X)$, $\{x\in X: \dim\o_{x,X}=0, \mbox{ or } =1 \mbox{ and   not regular, or not } S_2.\}$  (\ref{section.zeros.say}.3).

$\sim$,  linear equivalence of divisors $D_1, D_2$.

$\simq$,  $\q$-linear equivalence:  $mD_1\sim mD_2$ for some $m\neq  0$.

$\sims$,  linear similarity: $m_1D_1\sim m_2D_2$ for some $m_1, m_2\neq  0$ (\ref{QD.defn}). 

$\simsa$,  linear similarity +  $D_1, D_2$ ample and irreducible (\ref{simsa.defn}).
\end{say}

\medskip
{\bf Introduction to Sections~\ref{set.ci.prop.sec}--\ref{proj.sp.sec}}
\medskip
  
Let $X$ be a normal, projective  $K$-variety over some field $K$ and
$C\subset X$ an irreducible curve. In Sections~\ref{set.ci.prop.sec}--\ref{Mordell-Weil-fields.sec}
we study which finite subsets of $C$ are obtained as the intersection of $C$ with some divisor; a condition that depends only on the topology of the pair $|C|\subset |X|$.

Somewhat surprisingly, at the most basic level  the answer is
governed by the base field $K$. More precisely, by the  qualitative behavior of the group of $K$-points  of Abelian varieties over $K$. 
 There are 3 classes of fields $K$ for which the `size' of
$A(K)$  is about the same for every Abelian variety over $K$.
\begin{itemize}
\item (Finite fields)  For these  $A(K)$ is  finite. More generally,
if $K$ is {\it locally finite}  then $A(K)$ is a torsion group.
\item (Number fields)  For these  $A(K)$ is  finitely generated by the Mordell-Weil theorem.
More generally, the same holds for fields that are finitely generated over a prime field \cite{MR0102520}.
\item (Geometric fields)  For these  $A(K)$ has  infinite rank. This holds for example if $K$ is algebraically closed, except for $\bar \f_p$.
 These  are called  anti--Mordell-Weil fields \cite{2019arXiv190204011I}; see 
(\ref{anti-mordell.defn}) for the main examples. 
\end{itemize}

Roughly speaking, our results show how to read off closely related  properties of $K$ from the topology of $|X|$ in the first 2 cases and to recognize  
rational curves on $X$ in the third case.  

The conclusions are most complete  for  locally finite fields.

\begin{thm} Let $X$ be an irreducible, quasi-projective variety of dimension $\geq 2$ over a perfect  field $K$. The following are equivalent.
\begin{enumerate}
\item For every irreducible curve $C\subset X$ and every finite, closed  subset $P\subset C$, there is a divisor $D\subset X$ such that $C\cap D=P$  (as sets).
\item $K$ is locally finite.
\item  $A(K)$ is a torsion group for every commutative algebraic group  $A$ over $K$.
\end{enumerate}
\end{thm}

We have a more complicated characterization of the Mordell-Weil case.

\begin{thm} Let $X$ be an irreducible, quasi-projective variety of dimension $\geq 2$ over a  field $K$. The following are equivalent.
\begin{enumerate}
\item For every irreducible curve $C\subset X$ there is a finite, closed subset $\Sigma\subset C$ such that for  every finite, closed subset $P\subset C$, there is a divisor $D\subset X$ such that $P\subset C\cap D\subset P\cup \Sigma$  (as sets).
\item  $A(K)$ has finite $\q$-rank   for every Abelian variety $A$ over $K$.
\end{enumerate}
\end{thm}

In the geometric cases our considerations yield a Zariski-topological characterization of  rational curves.
The complete statement in (\ref{ci.char.thm.char0.rel}) needs several definitions, so here we state it somewhat informally.

\begin{prop}  Let $X$ be an irreducible, quasi-projective variety of dimension $n\geq 2$ over an anti--Mordell-Weil   field  and
$C\subset X$ an irreducible, geometrically reduced curve. One can decide using only the topology of the pair
$|C|\subset |X|$ whether $C$ is rational  or not. 
\end{prop}

The last result is especially useful if $X$ contains many rational curves, for example for $X=\p^n$. However, in  Sections~\ref{red.scip.sec}--\ref{proj.sp.sec} 
we get better results by observing that,
from the topological point of view, 
$$
(\mbox{hyperplane})\cup (\mbox{line})\subset \p^n
$$
is a very unusual configuration. This leads to the proof of (\ref{planes.in.Pn.thm.2.intro}).

\medskip
{\bf Introduction to Sections~\ref{sec.zero.set.sect}--\ref{top.pencils.sec}}
\medskip

A program towards showing that $|X|$ determines $X$ was started in \cite{k-lo}.  As an intermediate step,  \cite{k-lo} proposed to look for additional data that, together with $|X|$, determine $X$. A  candidate is {\it linear equivalence} of divisors, denoted by $D_1\sim D_2$.  
In (\ref{simsa.defn}) we introduce a variant, called
{\it linear  similarity of ample divisors} and denoted by $\simsa$.
Then we show first that  $|X|$ determines $\simsa$, then that
$\bigl( |X|, \simsa\bigr)$ determines $\sim$ and finally \cite{k-lo} proves that 
$\bigl( |X|, \sim\bigr) $ determines $X$. Symbolically:
$$
|X|\stackrel{(\mbox{\ref{main.thm.pf.step.1}})}{\longrightarrow}
\bigl( |X|, \simsa\bigr)
\stackrel{(\mbox{\ref{main.thm.pf.step.2}})}{\longrightarrow}
\bigl( |X|, \sim\bigr)
\stackrel{\mbox{\cite{k-lo}}}{\longrightarrow}
X.
$$


\begin{defn}[Linear  similarity of ample divisors]\label{simsa.defn}
 Let  $X$ be a normal variety and $\pdiv(X)$ the set of prime divisors on $X$. 
We define a relation on  $\pdiv(X)\times \pdiv(X)$ by
declaring that  $D_1\simsa D_2$ iff
\begin{enumerate}
\item $D_1, D_2$ are $\q$-Cartier, ample and
\item $m_1D_1\sim m_2D_2$ for some $m_1, m_2>0$.
\end{enumerate}
Note that  $\qrank\cl(X)=1$ iff  $D_1\simsa D_2$ for any 2 prime divisors on $X$. In these cases the relation $\simsa$ carries no extra information.
\end{defn}

The first step of the proof finds  $\simsa$.

\begin{prop}\label{main.thm.pf.step.1}   Let $X$ be a normal, projective variety of dimension $\geq 3$ over a field $k$. Assume that $k$ is not locally finite. Then $|X|$ determines $\simsa$.
\end{prop}

The proof  is actually  a quite short argument in Section~\ref{lin.sim.sec}, which is surprising since
$\sim$ and $\simsa$ seem very closely related at first sight.
We show how to recognize
\begin{enumerate}
\item irreducible, ample  $\q$-Cartier divisors (\ref{ample.is.top.lem}),
\item linear similarity of  irreducible, ample  $\q$-Cartier divisors
 (\ref{ample.Qiso.top.lem}), and
\item irreducible $\q$-Cartier divisors  (\ref{ample.Qiso.top.lem.var}).
\end{enumerate}

Then a longer argument shows that once we know $|X|$ and $\simsa$,
then we also know $\sim$.

\begin{thm}\label{main.thm.pf.step.2}
 Let $X$ be a normal, projective, geometrically irreducible variety over a field $k$ of characteristic 0. Assume that
\begin{enumerate}
\item either $\dim X\geq 4$,
\item or $\dim X\geq 3$ and $k$ is a   finitely generated  field extension of $\q$.
\end{enumerate}
Then  $|X|$ and $\simsa$ together determine $\sim$.
(See (\ref{lineq.top.prop.thm}) for a more general version.)
\end{thm}

This is  longer; 
  we step-by-step recognize the following objects/properties.
\begin{enumerate}
\item  $k$-points  (\ref{ZW.H.link.prop.1p.cor.2}). 
\item Isomorphism of residue fields of closed points (\ref{ZW.H.link.prop.1p.cor.3}). 
\item Isomorphism of reduced, 0-dimensional subschemes (\ref{iso.0-cyc.sinsa.say}). 
\item Transversality of  0-dimensional intersections of subvarieties (\ref{num.eq.top.say.det}). 
\item Two irreducible curves having the same degree (\ref{secs.isom.zeros.prop.c1}).
\item Two irreducible divisors having the same degree (\ref{secs.isom.zeros.prop.c2}).
\item $\q$-linear equivalence of ample divisors  (\ref{num.eq.top.say}.5). 
\end{enumerate} 
We can now use  (\ref{t.penc.alg.deg.char})  to recognize  algebraic pencils of divisors. 

The remaining problem is to decide which members of a topological pencil
are members of the corresponding algebraic pencil; this is the {\it true membership} problem  \cite[5.3.5]{k-lo-2}
As in \cite[Chap.5]{k-lo-2}, instead of  a complete solution, we only deal with  `well behaved' linear systems.  
We find sufficient  conditions for
\begin{enumerate}\setcounter{enumi}{7}
\item linearity of a pencil  (\ref{lin.test.lem}),
\item membership in a pencil  (\ref{true.memb.char.lem}) and
\item linear equivalence of reduced divisors (\ref{lin.eq.thm}). 
\end{enumerate} 
Then, as in \cite[5.3.8]{k-lo-2}, we see that  linear equivalences between reduced divisors generate the full linear equivalence relation. 
The next result then completes the proof of  
(\ref{home.iso.thm}). 

\begin{thm}\cite[Main Thm.]{k-lo} \label{main.thm.pf.step.0}
 Let $X$ be a normal, projective, geometrically irreducible variety of dimension $\geq 2$ over a field. Then  $|X|$ and $\sim$  together determine $X$. \qed
\end{thm}

The main technical tool for all this is the study of 
linkage of divisors.

\begin{say}[Linkage of divisors]\label{basic.plan.say}
 Let $X$ be a  normal, projective, irreducible  $k$-variety and
$L$  an ample line bundle on $X$.  Let
$Z_1, Z_2\subset X$ be closed, irreducible subvarieties such that $\dim (Z_1\cap Z_2)=0$.
We consider the following 
\medskip

{\it Linkage problem \ref{basic.plan.say}.1.} We say that  $s_i\in H^0(X, L^{m_i})$  with zero sets  $H_i:=(s_i=0)$ are {\it $L$-linked} on $Z_1\cup Z_2$ if there is an 
 $s\in H^0(X, L^r)$ with zero set  $H:=(s=0)$ such that
$$
\supp (H_1\cap Z_1)\cup  \supp (H_2\cap Z_2)=\supp \bigl((Z_1\cup Z_2)\cap H\bigr).
$$
Note that if the Picard number of $X$ is 1, then this is clearly a question involving only the underlying topology $|X|$. In fact, by (\ref{main.thm.pf.step.1}), this is almost always a question about $|X|$.

\medskip
{\it Sufficient condition  \ref{basic.plan.say}.2.}  If $s_1^{r_1}|_{Z_1\cap Z_2}=c\cdot  s_2^{r_2}|_{Z_1\cap Z_2}  $ for some
nonzero $c \in k^\times$ and $r_1, r_2\in \n$, then 
 the $H_i:=(s_i=0)$ are {\it $L$-linked} on $Z_1\cup Z_2$.
\medskip 

Next note that $s_i|_{Z_1\cap Z_2}$ can be an arbitrary element of
$H^0(Z_1\cap Z_2, L|_{Z_1\cap Z_2})\cong  H^0(Z_1\cap Z_2, \o_{Z_1\cap Z_2}) $ for $L$ sufficiently ample. Thus we obtain the following.
\medskip

{\it Proposition  \ref{basic.plan.say}.3.} If the sufficient condition  (\ref{basic.plan.say}.2) is necessary and any two $H_1, H_2$ are $L$-linked, then $H^0(Z_1\cap Z_2, \o_{Z_1\cap Z_2})\cong k $. That is, $Z_1\cap Z_2$ is a reduced $k$-point of $X$. 

This gives us a topological way to identify $k$-points and also check whether an intersection is transverse or not.

\medskip

{\it When is condition  (\ref{basic.plan.say}.2)  necessary?   \ref{basic.plan.say}.4.}  We know that 
$\supp (s_i|_{Z_i}=0)=\supp (s|_{Z_i}=0)$. Thus if
\begin{enumerate}
\item[(a)]  the $Z_i$ are normal, geometrically irreducible,  and
\item[(b)] the  $\supp (s_i|_{Z_i}=0)$ are irreducible, then  
\end{enumerate}
$$
s_i^{m_i}|_{Z_i}=c_i\cdot s^{n_i}|_{Z_i}
\eqno{(\ref{basic.plan.say}.4.c)}
$$  for some
nonzero  $c_i\in k^\times$ and $m_i, n_i\in \n$. Thus we conclude that there is  a constant $c\in k^\times$ such that
$$
s_1^{m_1n_2}|_{Z_1\cap Z_2}  =c\cdot s_2^{m_2n_1}|_{Z_1\cap Z_2}, 
\eqno{(\ref{basic.plan.say}.4.d)}
$$ 
as needed.  
\medskip

It remains to deal with the conditions (\ref{basic.plan.say}.4.a--b).
It turns out that  normality can be avoided,   this is worked out in Section~\ref{sec.zero.set.sect}.  Instead of geometric irreducibility,  the  key condition is geometric connectedness, which is guaranteed if the $Z_i$ are set-theoretic complete intersections of ample divisors. Thus the troublesome condition is (\ref{basic.plan.say}.4.b).
\medskip

{\it Irreducibility of   $H\cap Z$ \ref{basic.plan.say}.5.}
Let $Z\subset X$ be an irreducible subvariety and $H\subset X$ a general ample divisor. When is $H\cap Z$ irreducible?

Bertini's theorem says that this holds if $\dim Z\geq 2$. Since we also need
$\dim (Z_1\cap Z_2)=0$, we must have  $\dim X\geq 4$. For the best results we also need
$Z_1\cap Z_2$ to be a single point, which usually can be arranged only if
$\dim X\geq 5$. This is the case when the methods work best.

There are 2 ways to lower the dimension. First, it turns out that we get almost everything if one of the $Z_i$ has dimension $\geq 2$, the other can be a curve. Thus we get all results if $\dim X\geq 4$, as in (\ref{home.iso.thm}.1).

If $k$ is finitely generated over $\q$, then a theorem of Hilbert guarantees that $H\cap Z$ is irreducible for most ample divisors $H$ even if $\dim Z=1$.
This property defines Hilbertian fields (\ref{hilb.field.say}), and,   for such fields,  we can work with varieties of dimension $\geq 3$; leading to  (\ref{home.iso.thm}.2).

In order to give unified treatments, we  introduce the Bertini-Hilbert dimension of fields in (\ref{ber-hil.dim.defn}).

\end{say}

\medskip
{\bf Introduction to Sections~\ref{compl.int.sec}--\ref{spec.fields.sec}}
\medskip

These sections provide background material, most of which is known but may be hard to find in the form that we need. 
Section~\ref{compl.int.sec} summarizes properties of complete intersections and various Bertini-type theorems.

The theory of Picard group, Picard variety and  Albanese variety is recalled in Section~\ref{pic.cl.alb.sec}. The literature is much less complete about the  class group and its scheme version, which does not even seem to have a name. 

Basic results on commutative algebraic groups and the multiplicative group of
Artin algebras are studied in Section~\ref{comm.ag.sec}.

Section~\ref{spec.fields.sec} recalls the definitions and basic properties of 
locally finite, Mordell-Weil, anti--Mordell-Weil and Hilbertian fields.
We introduce the new notion of weakly Hilbertian fields  (\ref{weak.hilb.say}).

\begin{ack} My interest in this subject was started by a lecture of M.~Olsson about \cite{k-lo} at MSRI. After that, a series of  discussions with M.~Lieblich and  M.~Olsson led to further results \cite{k-lo-2}. These had a major influence on this work.

I received help from  M.~Larsen  with algebraic tori, B.~Poonen  with Bertini and Noether-Lefschetz theorems, J.~Silverman with N\'eron's theorem, 
K.~Smith     with discrete valuation rings and C.~Voisin with Noether-Lefschetz theorems; I am grateful to all of them.

Partial  financial support    was provided  by  the NSF under grant number
DMS-1901855.
\end{ack}

\section{Conjectures}

We discuss several  conjectures that came up while studying the above questions, but which are of independent interest.

\subsection*{Noether-Lefschetz theorem over countable fields}{\ }
\medskip

The following is probably  old, \cite{terasoma} attributes a version of it to T.~Shioda.

\begin{conj} \label{Noether-Lefschetz.conj}
Let $k$ be a field that is not locally finite, $X$ a normal projective $k$-variety of dimension $\geq 3$ and $H$ an ample Cartier divisor. Then, for $m\gg 1$ and for
`most' $k$-divisors  $D\in |mH|(k)$, the restriction maps
$$
\pic(X)\to \pic(D)\qtq{and} \cl(X)\to \cl(D) \qtq{are isomorphisms.}
$$
\end{conj}

The traditional satement of the Noether-Lefschetz theorem says
that the conclusion holds outside a countable union of
proper, closed subvarieties of $|mH|$; see \cite{sga2} for the Picard group and  \cite{MR2219849, MR2567426, ji-tocome} for the class group.
This gives a positive answer to the conjecture whenever  $k$ is uncountable.

For countable algebraically closed fields  a positive answer is given in 
\cite{and96, mau-poo, 2018arXiv181006481A, 2018arXiv181006550C}. 
See also \cite{terasoma} for similar results   over $\q$ for complete intersections in $\p^n$.

\subsection*{Independence of intersection points}{\ }
\medskip

To start with an example, 
let $E\subset \p^2$ be an elliptic curve,  and $L\subset \p^2$ a very general line
intersecting $E$ at 3 points $p_1, p_2, p_3$. Let $[H]\in \pic(E)$ denote the hyperplane class. One can  see that $m_1p_1+m_2p_2+m_3p_3\sim nH$ holds only for $m_1=m_2=m_3=n$.  This is easy over $\c$, gets quite a bit harder over
 $\q$ if we want the line to be also defined over $\q$.
The following conjecture says that a similar claim holds for all smooth curves.

\begin{conj} \label{indep.zeros.conj.0} Let $k$ be a field that is not locally finite.
Let $C$ be a smooth, projective curve  of genus $\geq 1$ over $k$ and 
 $L$  a very ample line bundle on $C$, 
 For a section $s\in H^0(C, L)$
write $\{p_i(s):i\in I\}$ (resp.\  $\{\bar p_i(s):i\in \bar I\}$) for the closed points  (resp.\  $\bar k$-points)   of $(s=0)$.
Then, for `most'  sections, we have injections
$$
\begin{array}{lcl}
\oplus_{i\in I} \z[p_i(s)]  &\into& \pic(C)
\qtq{(weak form),}\\[1ex]
\oplus_{i\in \bar I} \z[\bar p_i(s)] &\into& \pic(C_{\bar k})
\qtq{(strong form).}
\end{array}
$$
\end{conj}

It is not clear what  `most' should mean.  It is possible that 
this holds outside a field-locally thin set (\ref{thin.flt.defn}), but some heuristics suggest otherwise.
 
For the proof of (\ref{ZW.H.link.defn.conj}) we would need the following  stronger variant.
If true, it would  allow us to prove  (\ref{home.iso.thm}) for 3-folds as well.

\begin{conj} \label{indep.zeros.conj} Using the notation of 
(\ref{indep.zeros.conj.0}), let
$A\subsetneq \pico(C)$ be an Abelian subvariety and
$\Gamma\subset \pic(C)$ a finitely generated subgroup that contains $[L]$.
Then, for
`most'  sections, we have injections
$$
\begin{array}{lcl}
\oplus_{i\in I} \z[p_i(s)]\bigl\slash\tsum_{i\in I}   [p_i(s)]  &\into& \pic(C)\bigr\slash \langle  A(k), \Gamma\rangle
\qtq{(weak form),}\\[1ex]
\oplus_{i\in \bar I} \z[\bar p_i(s)]\bigl\slash\tsum_{i\in \bar I}   [\bar p_i(s)]  &\into& \pic(C_{\bar k})\bigr\slash \langle  A(\bar k), \Gamma\rangle
\qtq{(strong form).}
\end{array}
$$
\end{conj}

\medskip

\subsection*{Sections with few zeros}{\ }
\medskip

The next 2 conjectures posit that, for `most' ample line bundles, every section has many zeros. 

\begin{conj}\label{few.zeros.conj.2}
 Let $K$ be an algebraically closed field other than $\bar\f_p$. Let $C$ be a smooth, projective curve over $K$. Then, for `most' ample line bundles $L$, every section of $L^m$ has at least $g(C)$ zeros for every $m\geq 1$.
\end{conj}

Line bundles of degree $d$ that have  a section with fewer than  $g$ zeros
form a closed subset of dimension $ g-1$ of $\pics^d(C)$  obtained as the image of the maps
$$
\phi_{\mathbf m}: C^{g-1}\to  \pics^d(C)\qtq{given by}
(c_1,\dots, c_{g-1})\mapsto \o_C\bigl(\tsum_i m_i[c_i]\bigr),
$$
where ${\mathbf m}:=(m_1,\dots, m_{g-1})$ such that  $\tsum m_i=d$. 
Thus (\ref{few.zeros.conj.2})  is true if $K$ is uncountable. The most interesting open case is probably $\bar\q$. By (\ref{2.zeroes.lem.1}), there is a curve $C$ and a line bundle  $L$ over $\bar\q$,  such that  every section of $L^m$ has at least 2 zeros for every $m\geq 1$.

We prove the nodal rational curve cases of (\ref{few.zeros.conj.2}) in 
(\ref{barQ.not.wh.cor.cor}).

Thinking of the curve $C$ as a subvariety of its Jacobian leads to the following stronger form.

\begin{conj}\label{few.zeros.conj.1}
 Let $K$ be an algebraically closed field other than $\bar\f_p$. Let $A$ be an Abelian variety over $K$ and $Z_i\subset A$ subvarieties such that $\sum_i\dim Z_i<\dim A$.
Then, for `most' $p\in A(K)$, the equation
$$
n[p]=\tsum_i m_i[z_i]  \quad  n, m_i\in \z, z_i\in Z_i(K),
$$
has only the trivial solution  $n=m_i=0$. 
\end{conj}

Next we give an example with only one $Z_i$ where this holds.

\begin{exmp} Let $k$ be any field. Assume that $A=B\times E$ where $B$ is a  simple Abelian variety, $E$ an elliptic curve, and  we have only one $Z=Z_1\subset A$ of dimension $\leq \dim A -2$. Assume also that $Z$ does not contain any translate of $E$.

Let $\pi:A\to B$ be the coordinate projection.
If $p\in E(k), z\in Z(k)$ and $n[p]=m[z]$, then $m[\pi(z)]=0$, that is, $\pi(z)$ is a torsion point in $\pi(Z)$. By  \cite{MR1609518} there are only finitely many such,
so there are only finitely many  $\{z_j\in Z: j\in J\}$ for which there is an 
 $m_j> 0$  such that $m_j[z_j]\in E$. 

Thus if $p\in E(k)$ then   $n[p]=m[z]$ has a nontrivial solution  
iff $p$ is in the saturation of  $m_j[z_j]$ for some $j\in J$. 

If $\qrank E(k)\geq 2$, then finitely many subgroups of $\q$-rank 1 do not cover $E(k)$. 

To get such Jacobian examples, fix an elliptic curve $E$ over $\bar \q$ and
let $C$ be a sufficiently general member of a very ample linear system on $E\times \p^1$.  Then, by \cite[1.6]{koch-fl}, $\jacs(C)$ is isogeneous to  the product of $E$ and of a  simple Abelian variety $B$. 
\end{exmp}

\section{Set-theoretic complete intersection property}\label{set.ci.prop.sec}

\begin{defn}\label{ci.defn} Let $X$ be an irreducible scheme and $Z\subset X$ a closed, irreducible  subset.
We say that $Z$ has the {\it set-theoretic complete intersection property}---or that $Z$ is  {\it scip}---if the following holds.
\begin{enumerate}
\item  Let $D_Z\subset Z$ be a closed  subset  of pure codimension 1.   Then there is an effective divisor $D_X\subset X$ 
such that $\supp(D_X\cap Z)=D_Z$.
\end{enumerate}
In some cases only `nice' subvarieties $D_Z\subset Z$  are set-theoretic complete intersections. It is usually hard to formulate this in general, but the next variant allows us to ignore finitely many `bad' points of $Z$. 

We say that $Z$  is  {\it generically scip} if there is a finite (not necessarily closed) subset $\Sigma_Z\subset Z$ such that the following holds.
\begin{enumerate}\setcounter{enumi}{1}
\item  Let $D_Z\subset Z$ be a closed  subset  of pure codimension 1 that is disjoint from $\Sigma_Z$.   Then, 
for every finite (not necessarily closed) subset $\Sigma_X\subset X\setminus D_Z$,
there is an effective divisor $D_X\subset X$ that is disjoint from $\Sigma_X$ and such that $\supp(D_X\cap Z)=D_Z$.
\end{enumerate}
As a simple example, the quadric cone  $Q\subset \p^4$ is not scip  (over $\c$) but it is  generically scip with $\Sigma_Z=\{\mbox{vertex}\}$ and $\Sigma_X$ arbitrary.

The introduction of $\Sigma_Z$ means that we do not have to worry about some very singular points on $Z$. This is especially clear on curves, where we may assume that $\Sigma_Z$ contains all singular points.

The introduction of $\Sigma_X$ at first seems to make finding $D_X$ harder.
However, if $X$ is normal and $\Sigma_X$ contains all non-Cartier centers of $X$   (\ref{non.car.cent.defn}), then $D_X$ is a Cartier divisor. 
Thus we can usually work with the Picard group of $X$ (for which there are solid references), rather than  the class group  (for which modern references seem to be lacking). 

It is clear that these notions  depend only on the topological pair $|Z|\subset |X|$.

At the beginning we study the case when $Z$ is an irreducible  curve,
but later we need to understand many cases when $Z$ is reducible, and not even pure dimensional. 

We check in (\ref{genscip.pullback.lem}) that being generically scip  is invariant under purely inseparable morphisms and purely inseparable base field extensions. Thus, in order to save  considerable trouble with non-reduced group schemes, we usually work over perfect base fields.
\end{defn}

\medskip

\begin{lem} \label{ci.to.rc.lem}
Let $X$ be a normal, quasi-projective variety over a perfect field $k$ and $C\subset X$ an integral,   generically scip  curve.
Then 
$$\coker\bigl[\pic(X)\to \pic(C)\bigr]\qtq{is a torsion group.}
$$
\end{lem}

Proof.  We may assume that $\Sigma_Z\supset \sing C$ and
$\Sigma_Z\cup \Sigma_X$ contains all non-Cartier centers of $X$   (\ref{non.car.cent.defn}).
  Let $p\in C\setminus \Sigma_Z$ be a  point. By assumption 
there is an effective, Cartier  divisor $D_p$ 
such that $\supp(D_p\cap C)=\{p\}$. We do not know the intersection multiplicity at $p$,
so we can only say that  $\o_X(D_p)|_C\cong \o_C(m[p])$
for some $m>0$. (Here we use that $p$ is a regular point.) That is, $\o_C(p)$ is  torsion in
$\coker\bigl[\pic(X)\to \pic(C)\bigr]$. Since the $\o_C(p)$ generate
$\pic(C)$, we are done.\qed
\medskip

The rest of Sections~\ref{set.ci.prop.sec}--\ref{Mordell-Weil-fields.sec}  is essentially devoted to trying the understand the converse of (\ref{ci.to.rc.lem}).
Let us see first that the direct converse does not hold.

\begin{exmp} \label{rc.to.ci.1.exmp}
Let $C\subset \p^2$ be a smooth cubic over a number field. By the Mordell-Weil theorem 
$\pic(C)$ is finitely generated. Choose points  $\{p_i\in C: i\in I\}$ such that
$\o_C(1)$ and the $\o_C(p_i)$ form a basis of $\pic(C)\otimes\q$. 

Let $X$ be obtained by blowing up the points $p_i\in C\subset \p^2$.
Let $E_i\subset X$ be the exceptional curves
and let $C_X\subset X$  denote the  birational transform of $C$.
Note that $\pic(X)$ is spanned by the $E_i$  and the pull-back of $\o_{\p^2}(1)$.  Thus  $\pic(X)\to \pic(C_X)$ is an injection with torsion cokernel.  
\medskip

{\it Claim \ref{rc.to.ci.1.exmp}.1.}   $C_X\subset X$ is not generically scip.
\medskip

Proof. Choose $n_i>0$ and let $p\in C_X$ be a closed point such that
$[p]\sim \sum n_i[p_i]$.  Assume that $\{p\}=\supp(C_X\cap D)$ for some effective divisor $D\subset X$. Then  
$$
\o_X(D)|_{C_X}\cong \o_{C_X}(m[p])\cong \o_{C_X}\bigl(\tsum mn_i[p_i]\bigr)
$$
for some $m>0$. Since $\pic(X)\to \pic(C_X)$ is an injection,
this implies that
$D\sim \sum mn_i[E_i]$. But then $D= \sum mn_i[E_i]$ and
so $C_X\cap D=\{p_i:i\in I\}$. \qed
\end{exmp}

The following  is a partial converse of (\ref{ci.to.rc.lem}).

\begin{lem} \label{rc.to.ci.lem}
Let $X$ be a projective variety over a field $k$ and $C\subset X$ a reduced, irreducible  curve.
Assume that $$\coker\bigl[\pic^\circ(X)\to \pic^\circ(C)\bigr]\qtq{is a torsion group.}$$ Then $C$ is scip.
\end{lem}

Proof. Let $L$ be an ample line bundle such that $H^1(X, L^m\otimes T\otimes I_C)=0$ for every $m\geq 1$ and every $T\in \pic^\circ(X)$, where $I_C\subset \o_X$ denotes the ideal sheaf of $C$. 
Then $H^0(X, L^m\otimes T)\to H^0(C, L^m\otimes T|_C)$ is surjective for every $m\geq 1$. Set $d=\deg_CL$. For a  point  $p\in C$ let $P$ be a Cartier divisor on $C$ whose support is $p$ and set $r=\deg P$. (We can take $P=p$ if $p$ is a regular point.) 
Then   $L^r_C(-dP)\in \pic^\circ(C)$, thus
there is an  $m\geq 1$ and $T\in \pic^\circ(X)$ such that
$$
L^{mr}_C(-md P)\cong T^{-1}|_C.
$$
 This gives a section  $s_C\in H^0(C, L^m\otimes T|_C)$ whose divisor is $mdP$. It lifts to a section $s\in H^0(X, L^m\otimes T)$ and $D:=(s=0)$ works. \qed
\medskip

The next example shows that $C$ can be scip even if
$\coker\bigl[\pic^\circ(X)\to \pic^\circ(C)\bigr]$ is non-torsion.

\begin{exmp}\label{rc.to.ci.lem.2g.exmp}
Again let $C\subset \p^2$ be a smooth cubic over a number field.
 Choose a finite subset $L_i\in \pic^\circ(C)$, closed  under inverse,  that generates a full rank subgroup.

Choosing general sections in each  $L_i^{-1}(3)$, their zero sets $P_i\subset C$
are irreducible and distinct. Let
$S_i\to \p^2$ denote the blow-up of $P_i$ and  $C_i\subset S_i$  the birational transform of $C$. Then $\o_{S_i}(C_i)$ is a nef line bundle on $S_i$ 
and $ \o_{S_i}(C_i)|_{C_i}\cong L_i$. 

Finally consider the diagonal embedding   $C\subset \prod_i C_i\subset \prod_i S_i=:X$.
\medskip

{\it Claim \ref{rc.to.ci.lem.2g.exmp}.1.}   $C\subset X$ is scip.
\medskip

Proof. The key point is that
$\{ T|_C:  T\in \pic(X) \mbox{ and $T$ is nef}\}\subset \pic^\circ(C)$
contains a full  rank subgroup of $\pic^\circ(C)$. 
By Fujita's vanishing theorem  \cite[1.4.35]{laz-book},
there is  an ample line bundle $L$ such that $H^1(X, L^m\otimes T\otimes I_C)=0$ for every $m\geq 1$ and every nef $T$. The rest of the
argument in the proof of (\ref{rc.to.ci.lem}) works.
\end{exmp}

\begin{say}[Cokernel of $\pic(X)\to \pic(Y)$] \label{commensurate.lem}
(See (\ref{pic.cl.norm.say}) for definitions and notation involving the Picard group.)

If $X$ is a proper scheme then $\pic(X)$ is an extension of
$\ns(X)$ by $\pic^\circ(X)$. While
$\ns(X)$ is always a finitely generated abelian group, 
$\pic^\circ(X)$ can be trivial or very large, depending on the ground field and $X$. However, $\pics(X)$ is an algebraic group and  $\pics(X)(k)/\pic^\circ(X)$ is torsion. Thus, if $p:Y\to X$ is a morphism,   we aim to understand
$p^*:\pic(X)\to \pic(Y)$, in terms of 
 $$
p^*:\pico(X)\to \pico(Y) \qtq{and} p^*:\ns(X)\to \ns(Y).
\eqno{(\ref{commensurate.lem}.1)}
$$
We have better theoretical control of these  maps since the first is a map of Abelian varieties and the second a map of finitely generated abelian groups. 
\end{say}

\begin{prop}\label{basic.lower.bound.prop} Let $p:Y\to X$ be a morphism of proper $k$-schemes. Then
$$
\begin{array}{l}
\qrank\coker\bigl[\pic(X)\to \pic(Y)\bigr]\geq \\[1ex]
 \qquad \qquad\geq \qrank\coker\bigl[\pico(X)\to \pico(Y)\bigr](k)-\rank_{\q} \ns(X).
\end{array}
$$
\end{prop}

Proof. Let $\pic^*(X)\subset \pic(X)$ denote the preimage of
$\pic^\circ(Y) $.  Then $\pic^\circ(X)\subset \pic^*(X)$ and the
quotient is a subgroup of $\ns(X)$. Thus we see that
$$
\begin{array}{l}
\qrank\coker\bigl[\pic(X)\to \pic(Y)\bigr]\\
\geq  \qrank\coker\bigl[\pic^*(X)\to \pic^\circ(Y)\bigr]\\
\geq \qrank\coker\bigl[\pic^\circ(X)\to \pic^\circ(Y)\bigr]-\qrank \ns(X)\\
= \qrank\coker\bigl[\pico(X)\to \pico(Y)\bigr)(k)-\qrank \ns(X),
\end{array}
$$
where the last equality holds since the maps
$\pic^\circ(Z)\otimes \q\to \pico(Z)(k)\otimes \q$ are isomorphisms for proper $k$-schemes 
and  $A\mapsto A(k)\otimes \q$ is an exact functor of commutative group varieties (\ref{sch.onto.K.onto.lem}). \qed
\medskip

The first application is a characterization of locally finite fields.

\begin{thm}\label{quasifinite.field.thm}
 Let $X$ be an irreducible, quasi-projective variety of dimension $\geq 2$ over a perfect field $k$. 
The following are equivalent.
\begin{enumerate}
\item $k$ is locally finite.
\item Every irreducible curve  $C\subset X$ is scip.
\item Every irreducible curve  $C\subset X$ is generically scip.
\end{enumerate}
\end{thm}

Proof. 
Assume first that $k$ is locally finite and let $\bar X\supset X$
be a compactification.
Let $C\subset X$ be an irreducible curve with closure $\bar C$. If $k$ is locally finite then
$\pic^\circ(\bar C)(k)$ is torsion by (\ref{Q.rank.say}.1), hence 
 $\bar C$ is scip  by (\ref{rc.to.ci.lem}) and so is $C\subset X$.
(\ref{quasifinite.field.thm}.2) $\Rightarrow$ 
(\ref{quasifinite.field.thm}.3) is clear.

It remains to prove that (\ref{quasifinite.field.thm}.3) $\Rightarrow$ 
(\ref{quasifinite.field.thm}.1).  Note that if (\ref{quasifinite.field.thm}.3) holds for $X$ then it holds for every open subset of it, we may thus assume that $X$ is normal (or even smooth).
Let $\bar X\supset X$
be a normal compactification.

If $k$ is not locally finite then let $\bar C\subset \bar X^{\rm ns}$ be an irreducible  curve with at least 1 node contained in $X$. Note that $\pico(\bar X)$ is an Abelian variety
(\ref{pic.cl.norm.say}) and $\pico(\bar C)$ contains a $k$-torus (\ref{jac.say}). 
Thus $\coker\bigl[\pico(\bar X)\to \pico(\bar C)\bigr]$ contains a $k$-torus, 
hence its $\q$-rank is infinite (\ref{Q.rank.say}). 
Thus $\qrank \coker\bigl[\pic(\bar X)\to \pic(\bar C)\bigr]=\infty$ by
(\ref{basic.lower.bound.prop}) and so 
$\bar C$ is not generically scip by  (\ref{ci.to.rc.lem}).

Assume that the normalization of $\bar C$ has $m_\infty$ points lying over
$ \bar C\setminus  X$.  Then the kernel of the restriction map
$\pic(\bar C) \to \pic(C)$ has rank $\leq m_\infty$, thus we still have
$\qrank \coker\bigl[\pic(X)\to \pic(C)\bigr]=\infty$, hence
$C$ is not generically scip by  (\ref{ci.to.rc.lem}). \qed

\medskip

We can use this to strengthen the results of \cite{MR624904}.

\begin{prop} \label{quasifinite.field.thm.cor} Let $S$ be a normal, projective surface over a field $k$.
Then $|S|$ is homeomorphic to $\p^2_{\f_2}$ if and only if
\begin{enumerate}
\item $k$ is locally finite, and  
\item any 2 curves in $S$ have a non-empty intersection.
\end{enumerate}
\end{prop}

Proof.  By (\ref{quasifinite.field.thm}) every irreducible curve in
$\p^2_{\f_2}$ is scip, hence the same holds for $S$, so $k$ is locally finite by (\ref{quasifinite.field.thm}). 
Any 2 curves in $\p^2_{\f_2}$ have a non-empty intersection, hence the same holds for $S$. 

The converse, which is the combination of the next 2 claims,  is essentially proved  in \cite{MR624904}.

\medskip
{\it Claim \ref{quasifinite.field.thm.cor}.3.} 
 Up to homeomorphism there is a unique countable, noetherian, 2-dimensional topological space $X$ with the following properties.
\begin{enumerate}
\item[(a)] Every curve $C\subset X$ contains infinitely many points.
\item[(b)] Any 2 irreducible curves $C_1, C_2\subset X$ have a non-empty intersection.
\item[(c)] Given any 1-dimensional, closed $D\subset X$ and a 0-dimensional, closed $P\subset D$ such that $P\cap D_i\neq \emptyset$ for every irreducible component $D_i\subset D$, there is an  irreducible curve $C\subset X$
such that $D\cap C=P$.
\end{enumerate}

 Proof. Assume that we have 2 such $X, X'$. Choose well orderings of their points and curves  $Z_1, Z_2, \dots$ and  $Z'_1, Z'_2, \dots$.
Assume that we have a subset $I\subset \n$ and an injection
$\phi_I:I\to \n$ such that $\phi_I$ induces a homeomorphism
$$
\{Z_i: i\in I\}\longrightarrow \{Z'_{\phi_I(i)}: i\in I\}.
$$
Set $j:=\inf\{\n\setminus I\}$. We extend $\phi_I$ to
$J:=I\cup\{j\}$ as follows.

If $Z_j$ is a point on a $Z_i$ (resp.\ on none of them), we pick $Z'_{\phi_J(j)}$ to be any new point on  $Z'_{\phi_I(i)}$ (resp.\ on none of them).
If $Z_j$ is a curve, then first we pick images of all the new points in
$Z_j\cap Z_i$ for every $ i\in I$ and then use (c) to pick $Z'_{\phi_J(j)}$.

Continuing this way would only give an injection $\phi_{\infty}:\n\to\n$.
We fix this by alternating the construction between $\phi_I$ and
$\phi_I^{-1}$. \qed

\medskip
{\it Claim \ref{quasifinite.field.thm.cor}.4.} Let $K$ be locally finite field and $S$ a normal, projective surface over $K$. Then $|S|$ satisfies
(\ref{quasifinite.field.thm.cor}.3.a--c) iff  any 2 irreducible curves $C_1, C_2\subset S$ have a non-empty intersection.
\medskip

Proof. We need to prove (\ref{quasifinite.field.thm.cor}.c). First we blow up $P$ and normalize to get $S_1\to S$. Repeatedly blowing up points over $P$ we get
$S_r\to S$ such that the intersection matrix of $D_r$  (the birational transform of $D$) is negative definite. By \cite{Artin62},  $D_r\subset S_r$ can be contracted  to get  $\pi:S_r\to T$. 
By \cite{cha-poo},  there is  an irreducible hypersurface section $C_T\subset T$ that is disjoint from $\pi(D_r)$.  Let $C\subset S$ be its birational transform. \qed

\begin{lem} \label{genscip.pullback.lem}
Let $p:X'\to X$ be a  morphism between normal, projective varieties. Let $C'\subset X'$ be an irreducible  curve. Set $C:=p(C')$ and assume that
$k(C') /k(C)$ is purely inseparable. 

Then $C$ generically  scip    $\Rightarrow$ $C'$ generically  scip.
\end{lem}

Proof. Let $\Sigma_X\subset X$ be a finite subset 
that contains all non-Cartier centers and such that 
  $C'\setminus p^{-1}(\Sigma_X)\to C\setminus \Sigma_X$ is a bijection.

Pick any $q'\in C'\setminus p^{-1}(\Sigma_X)$ and set $q=p(q')$. 
There is a divisor $D(q)$ such that   $C\cap D(q)=\{q\}$. 
Then $D(q)$ is Cartier,  hence its pull-back gives a divisor
$D(q')$ such that  $C'\cap D(q')=\{q'\}$. \qed

\section{Mordell-Weil fields}\label{Mordell-Weil-fields.sec}

The Mordell-Weil theorem says that if $A$ is an Abelian variety over a number field $k$ then  $A(k)$ is a finitely generated group. 
Our results are not sensitive to torsion in $A(k)$, this is why we
need the concept of $\q$--Mordell-Weil fields  where
$\qrank A(k)$ is always finite; see (\ref{m.w.f.defn}).

 $\q$--Mordell-Weil fields have a nice characterization involving complete intersections on curves.

\begin{defn}\label{ci.defn.3} Let $X$ be a scheme and  $C\subset X$ a curve. 
\begin{enumerate}
\item   $C$  is  {\it scip with  defect} $\Sigma\subset C$ if, 
for every closed, finite subset $P\subset C$, there is an effective divisor $D=D(C,P)\subset X$ 
such that $P\subset \supp(D\cap C)\subset P\cup \Sigma$.
\item  $C$  is  {\it scip with finite defect} if it is
 scip with  defect $\Sigma$
for some finite subset   $\Sigma\subset C$.
\end{enumerate}



It is clear that these depend only on the topological pair $|C|\subset |X|$.

Note that being   scip with finite defect is invariant under birational  maps that are isomorphisms at the generic point of $C$, we just need to add to $\Sigma$ all  the indeterminacy points that lie on $C$.

Let  $X,Y$ be irreducible varieties and $\pi:Y\to X$  a dominant, finite  morphism. Let $C_X\subset X$ be a curve with reduced preimage
$C_Y\subset Y$. If $C_Y$ is scip with finite defect $\Sigma_Y$ then
$C_X$ is scip with finite defect $\Sigma_X:=\pi(\Sigma_Y)$. 

Thus most questions about these notions can be reduced to  normal, projective varieties.
 \end{defn}

\begin{lem}\label{sfd.fin.rank.prop.sm}  Let $X$ be a normal, projective variety over a perfect field $k$ and $C\subset X$ an irreducible  curve with normalization $\pi:\bar C\to C$.
Then $C$ is  scip  with finite  defect iff 
$$
\qrank\coker\bigl[ \pico(X)\to \jacs(\bar C)\bigr](k)<\infty.
\eqno{(\ref{sfd.fin.rank.prop.sm}.1)}
$$
\end{lem}

Proof.  Note first that the difference between
$$
\begin{array}{l}
\qrank\coker\bigl[ \pico(X)\to \jacs(\bar C)\bigr](k) \qtq{and}\\
\qrank\coker\bigl[ \pic(X)\to \pic(\bar C)\bigr]
\end{array}
$$
is at most $\qrank \ns(X)$.  

Assume that  $C$ is scip  with  finite  defect. We may enlarge $\Sigma_X$ to achieve that
$C\setminus \Sigma_X$ is smooth and every Weil divisor on $X$ not containing $\Sigma_X$ is Cartier along $C\setminus \Sigma_X$  (\ref{non.car.cent.thm.cor}). Let the defect set be $\Sigma\subset C\setminus \Sigma_X$ and 
 $\bar \Sigma\subset \bar C$ its preimage.
Let $\Gamma\subset \pic(\bar C)$ be the subgroup generated by
all $\bar q\in \bar\Sigma$.  

Pick any closed  point $p\in C\setminus \Sigma_X$. By assumption we have an effective Cartier divisor
$D_p$ such that $p\in \supp(D_p\cap C)\subset \{p\}\cup \Sigma$. 
This shows that, for some $m>0$,
$$
m[p]\in \bigl\langle \Gamma, \im\bigl[ \pic(X)\to \pic(\bar C)\bigr]\bigr\rangle.
$$
Since these $\{[p]: p\in C\setminus \Sigma\}$ generate $\pic(\bar C)$,
we get that 
$$
\qrank \coker\bigl[ \pic(X)\to \pic(\bar C)\bigr]\leq \qrank \Gamma.
$$

Conversely, assume that (\ref{sfd.fin.rank.prop.sm}.1) holds.
Curve singularities can be resolved, thus we ay assume that $C$ is smooth. 
Then there is  a finite subset
 $\{F_i:i\in I\}\subset \jac (\bar C)$, closed under inverses,
that generates $\coker\bigl[ \pic^\circ(X)\to \jac(\bar C)\bigr] $ modulo torsion.  Fix a point $p_0\in C\setminus \Sigma$ and an ample line bundle $L$ on $X$ such that  $\deg_C L=c\deg p_0$ for some $c>0$. 
Set $d_0=\deg p_0$ and choose $r_1$ such that
$r_1[p_0]+F_i\sim G_i$,
where the $G_i$  are effective and disjoint from $\Sigma_X$. 
 
Now pick any $p\in C\setminus \Sigma_X$ and
set $d=\deg p$. The line bundle
$\pi^*L^r\bigl(rc[p_0]+d[p_0]-d_0[p]\bigr)$ is in $\jac(\bar C) $. 
So, by assumption, there are
nonnegative $m_i$ and $T\in \pic^\circ(X)$ such that
$$
\pi^*L^r\bigl(rc[p_0]+d[p_0]-d_0[p]\bigr)\simq \o_{\bar C}\bigl(\tsum_i m_iF_i\bigr)\otimes T^{-1}.
$$
We can rewrite this as
$$
\pi^*(L^r\otimes T)\simq
\o_{\bar C}\bigl(d_0[p]+(rc-d-r_1\tsum m_i)[p_0]+\tsum_im_iG_i\bigr)
 $$
Every section of
$\pi^*(L^r\otimes T)|_C$ lifts back to a section of
$L^r\otimes T$ for large enough $r$. 
Thus $C$ is  scip  with   defect
$\Sigma=\supp(p_0+\tsum G_i)$. \qed

\begin{thm}\label{mor.wei.field.thm}
 Let $X$ be an irreducible, quasi-projective variety of dimension $\geq 2$ over a perfect field $k$. 
The following are equivalent.
\begin{enumerate}
\item $k$ is $\q$--Mordell-Weil.
\item Every irreducible curve  $C\subset X$ is scip with finite  defect.
\end{enumerate}
\end{thm}

Proof.  As we observed in (\ref{ci.defn.3}), it is enough to prove
(\ref{mor.wei.field.thm}.1) $\Rightarrow$ (\ref{mor.wei.field.thm}.2)
for normal, projective varieties.
If $k$ is $\q$--Mordell-Weil then (\ref{mor.wei.field.thm}.2) holds for these by (\ref{sfd.fin.rank.prop.sm}).

Conversely, if every irreducible curve  $C\subset X$ is scip with finite  defect then the same holds for the smooth locus  $X^{\rm ns}\subset X$.
Let $X'\supset X^{\rm ns}$ be a normal compactification. If 
 $ C'\subset X'$ has nonempty intersection with $X^{\rm ns}$
then $C'$ is also scip with finite  defect.

Now let $A$ be an Abelian variety over $k$. By (\ref{coker.any.quot})   there is an  irreducible,  projective curve $ C'\subset  X'$ and a
$\q$-injection
$$
A\into_{\q}  \coker\bigl[ \pico( X')\to \jacs(\bar C')\bigr].
$$
By (\ref{sfd.fin.rank.prop.sm})  $\coker\bigl[ \pico(X')\to \jacs(\bar C')\bigr](k)$ has  finite $\q$-rank, and  so does $A(k)$. 
Thus $k$ is $\q$--Mordell-Weil. 
 \qed

\medskip

By (\ref{sfd.fin.rank.prop.sm}), being  scip  with finite  defect depends on the interaction of $\pico(X)$ and   $\jacs(\bar C)$. The following definition is designed to get rid of the influence of $\pico(X)$.

\begin{defn}\label{ci.defn.3.abs} Let $X$ be a scheme,  $C\subset X$ an  irreducible curve.  We say that $C$ is  {\it absolutely scip  with finite defect} if
 the following holds.
\begin{enumerate}
\item  Let  $C'\neq C$  be any irreducible  curve. Then there  are finite subsets $\Sigma\subset C$ and $\Sigma'\subset  C'$ such that 
for every  finite subset $P\subset C$ there is an effective divisor $D\subset X$ 
such that $P\subset \supp\bigl(D\cap (C\cup C')\bigr)\subset P\cup \Sigma\cup \Sigma'$.

\end{enumerate}
 \end{defn}

Note that $P$ is a subset of $C$ only.
This has the following effect.

Let $\Gamma'\subset \jac(\bar C')$ be the subgroup generated by the 
preimages of $\Sigma'$. Let $\Gamma'_X\subset \pic(X)$
be the  preimage of $\Gamma'$ under $\pic(X)\to \jac(\bar C')$. 
Then the class of $D$ has to be in $\Gamma'_X$.

If $C'$ is general ample curve then the kernel of $\pic(X)\to \jac(\bar C')$
is torsion, thus  $\Gamma'_X$ is a finitely generated group. 

Now when we run the proof of (\ref{sfd.fin.rank.prop.sm})
for $C\subset X$, instead of the whole $\pic(X)$, we have only
$\Gamma'_X$ to choose $D$ from. 
The condition (\ref{sfd.fin.rank.prop.sm}.1) now becomes
$$
\qrank\coker\bigl[\Gamma'_X\to \pic(\bar C)\bigr]<\infty.
$$
Since $\Gamma'_X$  a finitely generated, this holds iff
$\qrank \pic(\bar C)<\infty$,
and we get the following.

\begin{prop}\label{sfgd.fin.rank.prop.rel}
  Let $X$ be a normal, projective variety of dimension $\geq 2$ over a perfect field $k$ and $C\subset X$ an irreducible   curve. 
Then $C$ is absolutely  scip  with finite  defect iff
$\pic(\bar C)$ has finite $\q$-rank. \qed
\end{prop}

This is especially useful over  fields where the opposite of the Mordell-Weil theorem  happens, these are the  {\it anti--Mordell-Weil} fields  (\ref{anti-mordell.defn}).
For varieties over such fields we can recognize rational curves 
using only their set-theoretic intersection properties.

Putting together (\ref{sfgd.fin.rank.prop.rel})  with (\ref{anti-mordell.defn}) gives the topological characterization of rational curves.

\begin{cor}\label{ci.char.thm.char0.rel} 
 Let $k$ be a perfect,   anti--Mordell-Weil field, 
 $X$ an irreducible, quasi-projective $k$-variety of dimension $\geq 2$ and  $C\subset X$   an irreducible   curve. Then $C$ is absolutely  scip  with finite  defect iff  every irreducible component of $C_{\bar k}$ is rational. \qed
\end{cor}

\begin{lem}[Curves with large Jacobians]\label{coker.any.quot}  Let $X$ be a geometrically normal, projective variety of dimension $\geq 2$  and  $A$  an Abelian variety over $k$.
Then there is an  irreducible,  projective curve $C\subset  X^{\rm ns}$ such that there is a  $\q$-injection  (that is, with finite kernel)
$$
A\into_{\q}  \coker\bigl[ \pico(X)\to \jacs(\bar C)\bigr].
$$
\end{lem}

Proof. 
Let $\bar C\subset A^\vee\times X$ be a curve that is a general, irreducible, complete intersection of
ample divisors; we use \cite{MR2385639} in case $k$ is finite. Let   $C\subset X$ be the image of the    second coordinate projection $\pi:\bar C\to C$. Then $C\subset X^{\rm ns}$ and $\bar C$ is the normalization of $C$. (In fact, $\bar C\cong C$ if $\dim X\geq 3$ and $k$ is infinite.)
By (\ref{ci.curve.onto.alb.cor})  the natural map
 $$
 A\times \pico(X) = \pico(A^\vee)\times \pico(X) \to \jacs(\bar C)
$$ is $\q$-injective, hence  
 we have a $\q$-injection
$A\into_{\q}  \coker\bigl[ \pico(X)\to \jacs(\bar C)\bigr]$.  \qed
\medskip

\section{Reducible scip subsets}\label{red.scip.sec}

\begin{defn}\label{ci.defn.red.3} Let $X$ be a scheme and $Z\subset X$ a closed   subset.
We say that  $Z$ is  {\it scip}  if 
the following holds.
\begin{enumerate}
\item  Let $D_Z\subset Z$ be a closed  subset  of pure codimension 1 that  has nonempty intersection with every irreducible component of $Z$.   Then there is an effective divisor $D_X\subset X$ 
such that $\supp(D_X\cap Z)=D_Z$.
\end{enumerate}
We say that  $Z$ is   {\it generically scip} if the following holds.
\begin{enumerate}\setcounter{enumi}{1}
\item  There is a (not necessarily closed) finite subset $\Sigma_Z\subset Z$ such that, if $D_Z$ in (\ref{ci.defn.red.3}.1) is disjoint from $\Sigma_Z$, then, 
for every (not necessarily closed) finite subset $\Sigma_X\subset X\setminus D_Z$,
we can find
$D_X\subset X$  as in (\ref{ci.defn.red.3}.1) that is also    disjoint from $\Sigma_X$.
\end{enumerate}
It is clear that these depend only on the topological pair $|Z|\subset |X|$.
Also, if (\ref{ci.defn.red.3}.2) holds for some $\Sigma_Z$ then it also holds for every larger $\Sigma'_Z\supset \Sigma_Z$.  We usually just take the union 
$\Sigma:=\Sigma_Z\cup \Sigma_X$ large enough.

If $Z$ is scip (resp.\ generically scip) then any union of its irreducible components is also scip (resp.\ generically scip).
\end{defn}

\begin{exmp}\label{ci.defn.red.3.exmp}  
In $\p^n$ with coordinates $x_0,\dots, x_n$, set
$L_1=(x_1=\cdots=x_i=0)$ and $L_2=(x_{i+1}=\cdots=x_n=0)$. We claim that
$L_1\cup L_2$ is generically scip.  Indeed, given divisors $D_{Z_i}\subset L_i$ not containing
$L_1\cap L_2=(1{:}0{:}\cdots{:}0)$, they can be given by
equations
$$
D_{Z_1}=\bigl(g_1(x_0, x_{i+1},\dots,x_n)=0\bigr)\qtq{and}
D_{Z_2}=\bigl(g_2(x_0, x_{1},\dots,x_i)=0\bigr).
$$
We may assume that   $g_i(1,0,\dots, 0)=1$. Then
$$
D_X:=\bigl(g_1^{\deg g_2}+g_2^{\deg g_1}=x_0^{\deg g_1\deg g_2}\bigr)
$$
satisfies  $\supp \bigl(D_X\cap (L_1\cup L_2)\bigr)=D_{Z_1}\cup D_{Z_2}$.

The following  strong converse  also illustrates the big difference between fields of characteristic 0, fields of positive characteristic and subfields of $\bar\f_p$.

\medskip
{\it Claim \ref{ci.defn.red.3.exmp}.1.}   
 Let $k$ be a field, $D\subset \p^n_k$ an irreducible divisor and $C\not\subset D$ an  irreducible, ample-sci  (\ref{ci.sci.say}) 
curve. Then
$C\cup D$ is scip (resp.\ generically scip) iff one of the following holds.
\begin{enumerate}
\item[(a)] $\chr k=0$, $C$ is a line and $D$ is a hyperplane.
\item[(b)]  $\chr k>0$, $\supp(C\cap D)$ is a single $k^{\rm ins}$-point and $C$, $D$ are both scip (resp.\ generically scip).
\item[(c)] $k$ is locally finite   and $D$ is scip (resp.\ generically scip).
\end{enumerate}

Proof.  Assume first that $\chr k=0$. By
(\ref{line.hyp.pl.X.lem})  $(C\cdot D)=1$, so $\deg C=1$ and $\deg D=1$. 
If  $\chr k>0$ but  $k$ is not locally finite, then 
 $\supp(C\cap D)$ is a single $k^{\rm ins}$-point by (\ref{Z.cap.W.is.geo.conn}.1). 
If $C\cup D$ is scip (resp.\ generically scip) then $C$ and $D$ are both scip (resp.\ generically scip). 
This shows that the conditions of (b) are necessary.
Their sufficiency also follows from (\ref{Z.cap.W.is.geo.conn}).
\qed

\medskip
Note that if $C$ is smooth and rational, then it is scip. If 
$n\geq 4$ and $D$ is smooth then it is scip.  
Let $D\subset \p^n$ be a smooth hypersurface and $n\geq 3$.
If $\deg D\leq n$ then there are  lots of smooth rational curves that meet $D$ in 1 point only.  If $\deg D\geq n+2$ then there should be few such curves, but there are  examples of arbitrary large degree. 
\end{exmp}

Next we prove  a general result about reducible scip subschemes.

\begin{notation}\label{k[X].not} For a $k$-scheme $Y$ we use  
$k[Y]:=H^0(Y, \o_Y)$  to denote the ring of regular functions.
If $Y$ is normal and proper then  $k[Y]=k$ iff $Y$ is geometrically integral.

If $Y$ is reduced  then $k(Y)$  denotes the ring of rational functions.
$Y$ is irreducible iff $k(Y)$ is a field.  If $Y_i$ are the irreducible components of $Y$ then  $k(Y)\cong \oplus_i k(Y_i)$. 
\end{notation}

\begin{prop} \label{Z.cap.W.is.geo.conn}
 Let $X$ be a normal, projective $k$-variety such that  $\rho(X)=1$.
 Let $Z,W\subset X$ be reduced, irreducible subvarieties such  that  $Z\cap W$ is 0-dimensional.  Assume that  $k$ is not locally finite.
Then   $Z\cup W$ is generically scip (\ref{ci.defn.red.3}) iff the following hold.
\begin{enumerate}
\item $Z$ and $W$ are generically scip,
\item $Z\cap W$ is irreducible,
\item either $k[\red(Z\cap W)]/k[W]$ or  $k[\red(Z\cap W)]/k[Z]$ is purely inseparable, and
\item if $\chr k=0$  then    $Z\cap W$ is reduced.
\end{enumerate}
\end{prop}

Proof.  Assume first that $Z\cup W$ is generically scip.
Choose any $\Sigma$ that contains $\Sigma(Z\cup W)$ (\ref{section.zeros.say}.3)
and the non-Cartier centers of $X$ (\ref{non.car.cent.defn}).
Let $L$ be an ample line bundle on $X$ such that
$H^0(X, L)\to H^0(Z\cap W, L_{Z\cap W})$ is surjective. 
Choose  sections $s_Z, s_W\in  H^0(X, L)$ that are nowhere zero on
$\Sigma$.    Write $ (s_Z|_Z=0)=\sum_i a_iA_i$ and 
 $ (s_W|_W=0)=\sum_j b_jB_j$. 
By assumption, for every $i,j$ there is
a divisor $D_{ij}\subset X$ such that 
$D_{ij}| _{Z\cup W}=c_{ij}A_i+d_{ij}B_j$ for some $c_{ij}, d_{ij}>0$. 
The $D_{ij}$ are Cartier since they are disjoint from $\Sigma$. 

For each $j$ a suitable positive linear combination of the $D_{ij}$ gives a divisor $D_j$ such that $D_j|_Z$ is a multiple of  $(s|_Z=0)$ and 
$D_j|_W$ is a multiple of $B_j$. Then we can take a
suitable positive linear combination $D$ of the $D_{j}$
such that $D|_Z$ is a multiple of  $(s_Z|_Z=0)$ and 
$D|_W$ is a multiple of $(s_W|_W=0)$.
 
Since $\rho(X)=1$, after passing to a suitable power we may assume that 
$$
D=(s=0) \qtq{for some} s\in H^0(X, L^m)\qtq{and} m>0. 
\eqno{(\ref{Z.cap.W.is.geo.conn}.5)}
$$
As we check in  (\ref{section.zeros.say}.4), this implies that 
$$
\begin{array}{lll}
s|_Z=u_Z s_Z^m|_Z & \mbox{for some} & u_Z\in k[Z]^{\times} \qtq{and}\\ 
s|_W=u_W s_W^m|_W & \mbox{for some} & u_W\in k[W]^{\times},
\end{array}
\eqno{(\ref{Z.cap.W.is.geo.conn}.6)}
$$
and hence
$$
(s_Z/s_W)^m|_{Z\cap W}=u_W|_{Z\cap W}\cdot u^{-1}_Z|_{Z\cap W}
\in \im \bigl[k[W]^\times \times k[Z]^{\times}\to 
k[Z\cap W]^{\times}\bigr].
\eqno{(\ref{Z.cap.W.is.geo.conn}.7)}
$$
We can arrange $s_Z/s_W$ to be an arbitrary element of 
$k[Z\cap W]^{\times} $, hence we conclude that
$$
k[Z\cap W]^{\times}\big\slash k[W]^\times \times k[Z]^{\times}
\qtq{is a torsion group.}
\eqno{(\ref{Z.cap.W.is.geo.conn}.8)}
$$
Now  (\ref{2.fields.mult.gen.lem}) shows that 
either $k[\red(Z\cap W)]/k[W]$ or  $k[\red(Z\cap W)]/k[Z]$ is purely inseparable, proving (\ref{Z.cap.W.is.geo.conn}.3). 
Since $Z, W$ are irreducible, $k[Z] $ and $k[W]$ are finite field extensions of $k$. Thus $k[\red(Z\cap W)]$ is a finite field extension of $k$, hence
$Z\cap W$ is irreducible. Finally (\ref{Z.cap.W.is.geo.conn}.4) follows from (\ref{ker.fin.Qrank.lem.2}).

Conversely, assume that (\ref{Z.cap.W.is.geo.conn}.1--4) hold and let $A_Z$ and $B_W$ be effective divisors on $Z$ and $W$ that are disjoint from $\Sigma$. By (\ref{Z.cap.W.is.geo.conn}.1)
they are both restrictions of Cartier divisors from $X$. 
Since $\rho(X)=1$, there is a power  $L^m$ and sections
$$
\sigma_Z\in H^0(Z, L^m|_Z)\qtq{and} \sigma_W\in H^0(W, L^m|_W).
$$
As in (\ref{Z.cap.W.is.geo.conn}.4--5) we see that suitable powers of 
$\sigma_Z^r$ and $\sigma_W^r$ can be glued to a section of
$\sigma_{Z\cup W}\in H^0\bigl(Z\cap W, L^{mr}|_{Z\cup W}\bigr)$ iff
 $$
\sigma_Z\sigma_W^{-1}\in k[Z\cap W]^{\times}\big\slash k[W]^\times \times k[Z]^{\times}
\qtq{is  torsion.}
\eqno{(\ref{Z.cap.W.is.geo.conn}.8)}
$$
This is guaranteed by (\ref{Z.cap.W.is.geo.conn}.3--4) using (\ref{2.fields.mult.gen.lem}).

Once $A_Z\cup B_W$ is defined as the zero set of a section
$\sigma_{Z\cup W}$, we can lift (a possibly higher power of) it to a section 
$\sigma_X\in H^0(X, L^N)$, and $D_X:=\supp (\sigma_X=0)$ shows that
$Z\cup W$ is generically scip.  
\qed

\begin{cor}\label{line.hyp.pl.X.lem}   Let $X$ be a smooth, projective $k$-variety such  that $\rho(X)=1$ and $\chr k=0$.
Let $C\subset X$ be a geometrically connected curve and $D\subset X$ a divisor.
If $C\cup D$ is generically scip then $(C\cdot D)=1$. \qed
\end{cor}

Looking at the above proof shows that there should be even fewer 
generically scip reducible subsets if $\rho(X)>1$,
but for now we have the following  slightly weaker result.

\begin{prop} \label{Z.cap.W.is.general}
 Let $X$ be a normal, projective $k$-variety where  $k$ is not locally  finite.
 Let $Z,W\subset X$ be reduced, irreducible subvarieties such  that  $Z\cap W$ is 0-dimensional. Assume that  $Z\cup W$ is generically scip (\ref{ci.defn.red.3}). Then
\begin{enumerate}
\item  $Z\cap W$ is irreducible,
\item either $k[\red(Z\cap W)]/k[W]$ or  $k[\red(Z\cap W)]/k[Z]$ is purely inseparable, and
\item if $\chr k=0$  then  
$$
\dim_k \ker\bigl[k[Z\cap W]\to k[\red(Z\cap W)]\bigr]\leq \frac{\rho(X)-1}{\deg [k:\q]}.
$$
\end{enumerate}
\end{prop}

Proof.  Choose $\Sigma$ to contain $\Sigma(Z\cup W)$ and the
non-Cartier centers of $X$ (\ref{non.car.cent.defn}).
As in (\ref{pic.Z.X.def}), let  
$\wdiv(X, \Sigma)$
denote the group of Weil divisors whose support is disjoint from
  $\Sigma$. These are all Cartier by our choice of $\Sigma$.
We get restriction maps
$r_Z:\wdiv(X, \Sigma)\to \wdiv(Z, \Sigma)$ and $r_W: \wdiv(X, \Sigma)\to \wdiv(W, \Sigma)$. 
These  descend to maps of the Picard groups
$\bar r_Z: \pic(X)\to \pic(Z)$ and  $ \bar r_W: \pic(X)\to \pic(W)$, which do not depend on $\Sigma$.
Set  $K_Z(X):=\ker \bar r_Z$, $K_W(X):= \ker  \bar r_W$ and
$K_{ZW}(X)=K_Z(X)\cap K_W(X)$.

As in  (\ref{pic.Z.X.def}),
the kernels of $\bar r_Z$ and $\bar r_W$ define  closed subgroups 
$ {\mathbf K}_Z(X)\subset \pics(X)$  and
$ {\mathbf K}_W(X)\subset \pics(X)$.  
Their intersection is denoted by $ {\mathbf K}_{ZW}(X)$.

Let $B$ be a  divisor in $\wdiv(X,\Sigma)$ whose class  $[B]$ lies in $K_{ZW}(X)$.  Then
$B|_Z=(s_Z)$, where $s_Z$ is unique up to $k[Z]^{\times} $, and $B|_W=(s_W)$,  where $s_W$ is unique up to $k[W]^{\times} $; here we use that $\Sigma\supset \Sigma(Z\cup W)$ and (\ref{section.zeros.say}).  Restricting both to $Z\cap W$  we get
$$
s_Z|_{Z\cap W}\cdot s_W^{-1}|_{Z\cap W}\in  k[Z\cap W]^{\times}\bigr\slash \bigl(k[W]^\times \times k[Z]^{\times}\bigr),
$$
which  defines  a homomorphism 
$$
K_{ZW}(X)\to  k[Z\cap W]^{\times}\bigr\slash \bigl(k[W]^\times \times k[Z]^{\times}\bigr).
$$
As in (\ref{gm.not.defn}),  we get   a homomorphism 
of algebraic groups
$$
\partial_{ZW}: {\mathbf K}_{ZW}(X) \to \bigl(\res^{Z\cap W}_k\gm\bigr)\bigr\slash \bigl(\res^{W}_k\gm\times \res^{Z}_k\gm\bigr).
$$
Note that ${\mathbf K}_{ZW}(X,\Sigma)\cap \pico(X)$ is an Abelian variety (\ref{pic.alb.say}.7), 
hence a positive dimensional subgroup of it has no morphisms to a linear algebraic group.  Set
$$
 \ns_{ZW}(X):=K_{ZW}(X)\bigr\slash\bigl(K_{ZW}(X)\cap \pico(X)(k)\bigr).
$$
Thus $\partial_{ZW}$  factors through
$$
\ns_{ZW}(X)\to \bigl(\res^{Z\cap W}_k\gm\bigr)\bigr\slash \bigl(\res^{W}_k\gm\times \res^{Z}_k\gm\bigr).
$$
Let $\Gamma_{ZW}(X)$  denote its image. 

By the above, $\Gamma_{ZW}(X)$ is the image of a subgroup of  $\ker\bigl[\ns(X)\to \ns(Z)\bigr]$ (modulo torsion).
All we need from this is that
$$
\qrank \Gamma_{ZW}(X)\leq \rho(X)-1.
$$

Now we start to follow the proof of (\ref{Z.cap.W.is.geo.conn}).
The departure from it happens at (\ref{Z.cap.W.is.geo.conn}.5), where  now 
$\sigma$ is not a  section of $L^m$, but of some $L^m(B)$ for some  $B\in K_{ZW}(X)$.
Thus we conclude that 
$$
(s/t)^m|_{Z\cap W}=u_W|_{Z\cap W}\cdot u^{-1}_Z|_{Z\cap W}\cdot \gamma
\in \im \bigl[k[W]^\times \times k[Z]^{\times}\to 
k[Z\cap W]^{\times}\bigr],
\eqno{(\ref{Z.cap.W.is.general}.3)}
$$
for some  $\gamma\in \Gamma_{ZW}(X)$. 
We can arrange $s/t$ to be an arbitrary element of 
$k[Z\cap W]^{\times} $, hence we conclude that
$$
k[Z\cap W]^{\times}\big\slash k[W]^\times \times k[Z]^{\times}\times \Gamma_{ZW}(X)
\qtq{is a torsion group.}
\eqno{(\ref{Z.cap.W.is.general}.4)}
$$

Now we use (\ref{2.fields.mult.gen.lem}) to get that 
$Z\cap W$ is irreducible and 
 (\ref{ker.fin.Qrank.lem.2}) implies (\ref{Z.cap.W.is.general}.2). Finally we get that 
$$
\qrank \ker\bigl[k[Z\cap W]\to k[\red(Z\cap W)]\bigr]\leq \qrank\Gamma_{ZW}(X)\leq 
\rho(X)-1.
$$
Thus (\ref{Z.cap.W.is.general}.3) follows from (\ref{Q.rank.say}.3). \qed

\medskip

In characteristic 0 we can reformulate the bound of (\ref{Z.cap.W.is.general}) as follows.

\begin{cor} \label{Z.cap.W.is.general.cor}
 Let $X$ be a normal, projective variety over a field $k$ of characteristic 0.
 Let $Z,W\subset X$ be reduced, irreducible subvarieties such  that  $Z\cap W$ is 0-dimensional. Assume that  $Z\cup W$ is generically scip (\ref{ci.defn.red.3}). Then
  $Z\cap W$ is irreducible and
$$
\dim_k k[Z\cap W]\leq \max\{\dim_k k[Z], \dim_k k[W]\}+ \frac{\rho(X)-1}{\deg [k:\q]}. \qed
$$
\end{cor}

\begin{exmp}\label{k.pts.dim.7.exmp} Combining the ideas of (\ref{ci.defn.red.3.exmp}) with (\ref{Z.cap.W.is.geo.conn}) we get a method to recognize $k$-points. The assumptions are restrictive, but this gives the first indication that one can get detailed scheme-theoretic information from the topology. However, scip and generically scip turn out to be too restrictive in general; searching for a more flexible variant  lead to the notion of  linkage in Section~\ref{link.divs.res.sect}.

\medskip
\noindent{\it Claim \ref{k.pts.dim.7.exmp}.1.}  Let $X$ be a smooth, projective $k$-variety of dimension $\geq 7$ such that  $\rho(X)=1$.
Then $p\in X$ is a $k$-point iff there are 3-dimensional, set-theoretic complete intersections    $Z,W\subset X$
such that 
\begin{enumerate}
\item $\supp(Z\cap W)=\{p\}$ and
\item $Z\cup W$ is generically scip.
\end{enumerate}
\medskip

Proof. Assume that $p\in X$ is a $k$-point and let
$Z,W\subset X$ be 3-dimensional, smooth,  complete intersections such that
$Z\cap W=\{p\}$ as schemes. Lefschetz theorem tells us that if $D_Z\subset Z$ is any divisor then (some multiple of) it is a complete intersection. 
Arguing as in (\ref{ci.defn.red.3.exmp}) we see that $Z\cup W$ is generically scip.

Conversely, $k[Z]=k[W]=k$ since $Z,W$ are set-theoretic complete intersections (\ref{sci.say}.1), thus  (\ref{Z.cap.W.is.geo.conn}) says that $p\in X(k)$.\qed
\medskip

Note that the bound $\dim X\geq 7$ can be improved to $\dim X\geq 5$ if
the Noether-Lefschetz theorem applies over $k$; see (\ref{Noether-Lefschetz.conj}) for such cases.
\end{exmp}

\section{Projective spaces}\label{proj.sp.sec}

We study the scip property for the union of a curve and of a divisor.
As we observed in Section~\ref{red.scip.sec},
this happens very rarely, and it leads to the following  stronger version of 
(\ref{planes.in.Pn.thm.2.intro})

\begin{thm}\label{planes.in.Pn.thm.2}
Let $L$ be a  field  of characteristic 0 and  $K$    an arbitrary  field.
Let $Y_L$ be a normal, projective, geometrically irreducible  $L$-variety of dimension $n\geq 2$  and 
 $\Phi:|\p^n_K|\sim |Y_L|$  a homeomorphism.  Then 
\begin{enumerate}
\item $Y_L\cong \p^n_L$,
\item $K\cong L$, and 
\item $\Phi$  
is the composite of a field isomorphism  $\phi:K\cong L$  and of an  automorphism of $\p^n_K$. 
\end{enumerate}
\end{thm}

We start with an easy to prove but interesting special case of 
(\ref{planes.in.Pn.thm.2}).

\begin{say}[Proof of (\ref{planes.in.Pn.thm.2}) when $Y_L\cong \p^n_L$]\label{planes.in.Pn.thm.pf.1} Let $H\subset \p^n_K$ be a hyperplane and
 $\ell\subset \p^n_K$  a line not contained in $H$. Then
$\ell\cup H$ is scip by (\ref{ci.defn.red.3}.2) hence
so is $\Phi(\ell)\cup \Phi(H)\subset \p^n_L$. So 
$\Phi(H)\subset \p^n_L$ is a hyperplane
by (\ref{ci.defn.red.3.exmp}.1.a). By taking intersections, we see that
$\Phi$ gives an isomorphism of the projective geometries
$K\p^n$ and $L\p^n$. By  the Veblen-Young theorem \cite{MR1506049}
this is induced by a field  isomorphism $\phi:K\cong L$.

Composing $\Phi$ with the natural isomorphism induced by $\phi^{-1}$,
we get a homeomorphism $\Psi:|\p^n_K|\to |\p^n_K|$ that is the identity on $K$-points. It remains to show that it is the identity on all points. 
Let $C\subset \p^n_K$ be a $K$-rational curve. It has infinitely many $K$-points and these are fixed by $\Psi$.   Thus $C\cap \Psi(C)$ is infinite, hence
$C=\Psi(C)$. However, we do not yet know that $\Psi|_C$ is the identity.

Let $p\in \p^n_K$ be a closed point.  Assume that there are $K$-rational curves $C_{\lambda}\subset \p^n_K$ such that $\{p\}=\cap C_{\lambda}$. 
Then
$\{\Psi(p)\}=\cap \Psi(C_{\lambda})=\cap C_{\lambda}=\{p\}$,
as needed. 

It remains to construct such curves $C_{\lambda}$. For this we can work in an affine chart  $p\in \a^n_K\subset \p^n_K$ with coordinates $x_i$.
Note that  $K(p)/K$ is a finite, separable extension, hence can be generated by a single element $z_p\in K(p)$. 
We can thus write
$x_i(p)=h_i(z_p)$ for some  $h_i\in K[t]$.

Let  $g(t)\in K[t]$ be the  minimal polynomial of $z_p$; we can then identify $z_p$ with a root of $g$ in $\bar K$.

For $i\in \{1,\dots, n\}$  and $a\in K$ let $C_{i,a}$ be the image of
$$
\tau_{i,a}:t\mapsto   \bigl(h_1(t), h_2(t),\dots, h_n(t)\bigr)+ ag(t)e_i,
$$
where $e_i$ is the $i$th standard basis vector.

The $C_{i,a}$ are $K$-rational curves, hence stabilized by $\Phi$. 

\medskip
{\it Claim \ref{planes.in.Pn.thm.pf.1}.1.}  $\cap_{i,a} C_{i,a}=\{p\}$.   
\medskip

Proof. 
First,  the $\tau_{i,a}$ all map $z_p\in \a^1(\bar K)$ to $p$, so $p\in \cap_{i,a} C_{i,a}$.  To see the converse, assume that
$p\neq q\in \cap_{i,a} C_{i,a}$. After permuting the coordinates we may assume that  $p_n\neq q_n$.  If $q=\tau_{1,a}(z')$ then $h_n(z')=q_n$.
The equation  $h_n(*)=q_n$ has  finitely many solutions  $z'_j$ and they are all different from $z_p$. Then, for all but finitely many  $a\in K$, $h_1(z'_i)+ag(z'_i)\neq q_i$  for every $i$. Thus $q\not\in C_{1,a}$. \qed
\medskip

In positive characteristic  the above proof and (\ref{Z.cap.W.is.geo.conn}) give the following.
\medskip

{\it Claim \ref{planes.in.Pn.thm.pf.1}.2.}
 Let $K, L$ be  perfect  fields that are not locally finite,  $n\geq 2$ and 
$\Phi:|\p^n_K|\sim |\p^n_L|$ a  homeo\-morphism. Then $\Phi$ induces a bijection
$\p^n(K)\leftrightarrow \p^n(L)$. \qed
\end{say}

For the proof of (\ref{planes.in.Pn.thm.2}), the key step is the following.

\begin{lem}\label{planes.in.Pn.thm.2.lem}
Using the notation and assumptions of (\ref{planes.in.Pn.thm.2}), 
let $|H|$ denote the linear system of all
hyperplanes in $\p^n_K$. Then 
 $\{\Phi(H): H\in |H|(K)\}$ is a bounded family of divisors on $Y_L$.
\end{lem}

\begin{say}[Proof of (\ref{planes.in.Pn.thm.2}) assuming (\ref{planes.in.Pn.thm.2.lem})]
First note that $K$ is not locally finite  by (\ref{quasifinite.field.thm}).

Set $Z:=\Phi^{-1}(\sing(Y_L))$.  Let $|D|\subset |H|$ 
be a pencil of hyperplanes 
whose base locus is not contained in $Z$,
and
$\{D_{\lambda}: \lambda\in \Lambda\}$ the corresponding t-pencil (\ref{t.pencil.defn}).
Thus $\{\Phi(D_{\lambda}): \lambda\in \Lambda\}$ is a t-pencil on $Y_L$.

There are infinitely many hyperplanes among the $\{D_{\lambda}: \lambda\in \Lambda\}$ and their $\Phi$-images form a bounded family of divisors
by  (\ref{planes.in.Pn.thm.2.lem}). 
Thus   $\{\Phi(D_{\lambda}): \lambda\in \Lambda\}$ is
algebraic (\ref{t.penc.alg.deg.char}), linear (\ref{lin.test.lem})  and
the images of the $K$-hyperplanes are true members (\ref{true.memb.char.lem}).
Thus, by (\ref{build.linsys.lem}),  the $\{\Phi(H): H\in |H|(K)\}$ span an $n$-dimensional  linear system 
 $|H|^Y$, which is  basepoint-free
since already the $\{\Phi(H): H\in |H|(K)\}$ have no point in common.
So $|H|^Y$ gives a morphism  $g:Y\to \p^n_L$. 
Since any hyperplane $H$ has nonempty intersection with every curve, the same holds for $\Phi(H)$, so $g:Y\to \p^n_L$ is finite and  $|H|^Y$ is ample. 

Therefore  members of $|H|^Y$ are geometrically connected, and so are the $\Phi$-images of the lines, since they are set-theoretic complete intersections of members of  $|H|^Y$  (\ref{sci.say}.1).
By (\ref{Z.cap.W.is.geo.conn}) 
this shows that 
$\supp(\Phi(\ell)\cap \Phi(H)\bigr)$ is an $L$-point whenever
 it is a smooth point of $Y_L$. 
We also obtain this point as 
$\Phi(H_1)\cap\cdots\cap  \Phi(H_n)$, or as a fiber of
$g:Y\to \p^n_L$. Since $\chr L=0$, general fibers of $g$ are reduced.
 Thus   $g:Y\onto \p^n_L$ is finite and of degree 1, hence an isomorphism.  The rest now follows from (\ref{planes.in.Pn.thm.pf.1}).
\qed
\end{say}

The following lemma, which essentially says that pencils determine higher dimensional linear systems,  is longer to state than to prove.

\begin{lem}\label{build.linsys.lem} Let $Y$ be a normal, projective variety over a field $L$. Let $K$ be an infinite field and $e_0, \dots, e_n\in K\p^n$
independent points. Assume that we have a map
$$
\Phi: K\p^n\to (\mbox{effective Weil divisors on $Y$})
$$
with the following property. 
\begin{enumerate}
\item For $r=1,\dots, n$ there are Zariski open subsets $\emptyset\neq U_{r-1}\subset \langle e_0, \dots, e_{r-1}\rangle$ such that, for every $p\in U_{r-1}$, the divisors
$\{\Phi(q): q\in \langle p, e_{r}\rangle\}$
are $L$-members of a linear pencil on $Y$. 
\end{enumerate}
Then there is a Zariski open subset $W\subset  K\p^n$ such that 
the divisors
$$
\{\Phi(q): q\in W\}
$$
are $L$-members of a linear system of dimension $\leq n$  on $Y$. \qed
\end{lem}

\begin{say}[Proof of  (\ref{planes.in.Pn.thm.2.lem}) when $\dim Y= 2$]
\label{planes.in.Pn.thm.2.lem.pf.2}
Let $\ell,\ell'$ be 2 lines in $\p^2_K$. 
We aim to show that  $\bigl(\Phi(\ell)\cdot \Phi(\ell')\bigr)$ is bounded.
If $\Phi(\ell)\cap \Phi(\ell')$ is a smooth point of $Y$, then
$$
\bigl(\Phi(\ell)\cdot \Phi(\ell')\bigr)=\dim_kk\bigl[\Phi(\ell)\cap \Phi(\ell')\bigr],
$$
and, by (\ref{Z.cap.W.is.general.cor}), the latter is bounded by
$$
\max\bigl\{ \dim_kk[\Phi(\ell)],  \dim_kk[\Phi(\ell)]\bigr\}+\frac{\rho(Y)-1}{\deg [L:\q]}.
$$
If $C\subset Y$ is any integral curve then $\dim_k k[C]$ is
the number of connected components of $C_{\bar k}$. By (\ref{PTBPS.thm}) 
then either $\dim_k k[C]<\rho(Y_{\bar k})$, or $C$ 
is contained in members of a basepoint-free pencil. In the latter case it is disjoint from other members of the pencil, which does not happen for
$\Phi(\ell)$. 

Thus we conclude that 
$\bigl(\Phi(\ell)\cdot \Phi(\ell')\bigr)$ is bounded, except possibly for lines that meet at a singular point of $Y$.

We check in (\ref{inf.cur.on.surf.lem}) that almost all 
$\Phi(\ell)$ have positive self-intersection number.
Since  $\Phi(\ell)$ has positive intersection number with every other curve, 
we see that almost all $\Phi(\ell)$  are ample.
A general such $\Phi(\ell)$ is contained in the smooth locus of $Y$, hence it has bounded intersection number with every other $\Phi(\ell')$.
Hence the $\Phi(\ell)$ form a bounded family. 
\qed
\end{say}

\begin{lem} \label{inf.cur.on.surf.lem}
Let $S$ be a normal, projective surface and  $\{C_i: i\in I\}$ an infinite collection of irreducible curves on $S$ such that the 
intersection numbers $\{(C_{i_1}\cdot C_{i_2}): i_1\neq i_2\in I\}$ are positive and bounded from above. Then $(C_i^2)>0$ for almost all $i\in I$.
\end{lem}

Proof.  Set $d_1:=\sup \{(C_{i_1}\cdot C_{i_2}): i_1\neq i_2\in I\}$. 
First we claim that the  self-intersections $(C_i^2):i\in I$ are bounded from below  by some $-d_2<0$. If not then for every $r$ there is a subset  $J_r\subset I$ of $r$ elements such that $(C_j^2)<-rd_1$ for every $j\in J_r$.  Then the
intersection matrix of the curves  $\{C_j:j\in J_r\}$ is negative definite. This is only possible for  $r<\rho^{\rm cl}(S)$  (\ref{cl.norm.say}).

Note that the local intersection numbers at a point $s\in S$ are bounded from below by some $\beta(s)>0$  (the inverse of the determinant of the intersection matrix of the exceptional curves on  a resolution). 
So there is a $\beta>0$ such that  
$(C_{i_1}\cdot C_{i_2})\geq \beta$ for every $ i_1\neq i_2\in I$.

Let $J\subset I$ be any finite subset of $>1+d_2/\beta$ elements.
Then  $D:=\tsum_{j\in J} C_j$ is nef and big. Since
$(D\cdot C_i)\leq d_1|J|$ for every $i\in I\setminus J$,  these $C_i$ are contained in a  family of curves parametrized by a $k$-scheme of finite type $W$ by (\ref{nef.big.bdd.lem}). The  family of curves parametrized by $W$ contains only finitely many
curves with negative self-intersection,  since these are isolated points of  $W$. Also, $W$ contains only 
 finitely many basepoint-free pencils, and  $\{C_i: i\in I\}$ can contain at most one member of each such pencil since $(C_{i_1}\cdot C_{i_2})>0$ for $ i_1\neq i_2\in I$. Thus $(C_i^2)>0$ for all but finitely many  $i\in I$.
 \qed

\begin{lem} \label{nef.big.bdd.lem}
Let $S$ be a normal, projective surface and $M$ a   big $\q$-divisor on $S$. Then for every $m$, the set of curves
$$
\{C\colon   C \mbox{ is irreducible and } (M\cdot C)\leq m\}
$$
is a bounded family.
\end{lem}

Proof. Note that on a normal surface the intersection numbers  make sense even for divisors that are not $\q$-Cartier.  After a resolution and pull-back, we may assume that $S$ is smooth.
By Kodaira's lemma  \cite[2.60]{km-book}  we can write  $M\simq A+E$ where $A$ is ample and $E$ is effective.  Thus either $C\subset E$ or $(A\cdot C)\leq m$. \qed
\medskip

\begin{say}[Proof of  (\ref{planes.in.Pn.thm.2.lem}) when $\dim Y\geq 3$]
\label{planes.in.Pn.thm.2.lem.pf.3}
In this case the $\Phi(H)$ are ample by (\ref{ample.is.top.lem}), hence
 geometrically connected, and so are the $\Phi$-images of the lines since they are set-theoretic complete intersections of members of  $|H|^Y$.  Thus
$$
\dim_k\bigl[\Phi(\ell)\cap \Phi(H)\bigr]\leq 1+\frac{\rho(Y)-1}{\deg [k:\q]}
\eqno{(\ref{planes.in.Pn.thm.2.lem.pf.3}.1)}
$$
by (\ref{Z.cap.W.is.general.cor}). 
If $\Phi(\ell)\cap \Phi(H)$ is a smooth point of $Y$, then
$$
\bigl(\Phi(\ell)\cdot \Phi(H)\bigr)=\dim_k\bigl[\Phi(\ell)\cap \Phi(H)\bigr].
$$
So the $\Phi(H)$ have bounded intersection number with
a curve that is an intersection of ample divisors, hence they form a bounded family.  \qed
\end{say}

The following result is proved in \cite{bopisi}, sharpening earlier versions of \cite{MR1786508, MR2177196}. 

\begin{thm} \label{PTBPS.thm} Let $X$ be a  normal, projective, irreducible  $k$-variety and $\{D_i: i\in I\}$ pairwise disjoint divisors.
Then
\begin{enumerate}
\item either  $|I|\leq \rho^{\rm cl}(X)-1$, 
\item or all the $D_i$ are contained in members of a basepoint-free pencil. \qed
\end{enumerate}
\end{thm}

\section{Sections and their zero sets}\label{sec.zero.set.sect}

We discuss foundational results about sections and their zero sets that are needed in our study of linkage.

\begin{say}\label{section.zeros.say}
Let $X$ be a normal, geometrically integral, proper $k$-variety,
$L$  a line bundle on $X$ and $s_1, s_2\in H^0(X, L)$ sections of $L$ with corresponding divisors  $(s_i=0)$. Then  $(s_1=0)=(s_2=0)$ iff
$s_1=s_2\cdot c$ for some  $c\in k^\times$. 
Our aim is to relax normality as much as possible and still keep the conclusion.

 Let $Y$ be a reduced scheme,  $L$  a line bundle on $Y$ and $s\in H^0(Y, L)$ a section. It has  {\it  scheme-theoretic zeros}  $(s=0)$ and 
 {\it  divisor-theoretic zeros;} the latter is the Weil divisor
$\sum_\eta \len_{k(\eta)}(L/s\o_Y) [\bar\eta]$, where the summation is over all codimension 1 points of $Y$. The scheme-theoretic zeros determine the 
divisor-theoretic zeros, but the converse does not always hold.

 We consider 2 genericity conditions.
\begin{enumerate}
\item Every generic point of $\supp(s=0)$ is a regular, codimension 1 point of $Y$.
\item $Y$ is  $S_2$ along $\supp(s=0)$.
\end{enumerate}
If (\ref{section.zeros.say}.1) holds then $(s=0)$ is well-defined as a  Weil divisor that is generically Cartier. 
If (\ref{section.zeros.say}.2) holds then the zero set  $(s=0)$ has no embedded points. Thus if (\ref{section.zeros.say}.1--2) both hold then the divisor-theoretic zeros determine the scheme-theoretic zeros. In these cases we do not distinguish the 2 concepts and simply talk about the {\it zero set,}  and denote it by $(s=0)$.  

Both of these conditions hold, except at some special points.

\medskip
{\it Definition \ref{section.zeros.say}.3.} For a  reduced scheme $Y$, 
let $ \Sigma(Y)\subset Y$ denote the set of points $y\in Y$ such that
$\o_{y,Y}$ is either of dimension 0,  or of dimension 1 but not regular, or  not $S_2$. 

If $Y$ is of finite type then $ \Sigma(Y)$ is finite.
Note that $s$ satisfies the conditions (\ref{section.zeros.say}.1--2) iff  $(s=0)$ is disjoint from  $ \Sigma(Y)$. 

The usual correspondence between divisors and sections holds outside $\Sigma(Y)$. 

\medskip
{\it Claim \ref{section.zeros.say}.4.}  Let $Y$ be a reduced scheme,
 $L$  a line bundle on $Y$ and $s_1, s_2\in H^0(Y, L)$ sections that do  not vanish at any point of  $\Sigma(Y)$. 
Then  $(s_1=0)=(s_2=0)$ iff  $s_1=s_2\cdot u$ for some $H^0(Y, \o_Y)^\times$. \qed

\medskip

{\it Example \ref{section.zeros.say}.5.} Let $Y\subset \p^4_k$ be the union of
$(x_1=x_2=0)$ and of $(x_3=x_4=0)$. Note that  $(1{:}0{:}0{:}0{:}0)$ is a non-$S_2$ point and $H^0(Y, \o_Y)=k$. 

 Consider  $s(a,c)=ax_1+cx_3\in H^0(Y, \o_Y(1))$. Note that
its divisor  $\bigl(s(a,c)=0\bigr)$ is independent of $a,c\in k^\times$.
However,  $s(a,c)=s(a',c')\cdot u$ for some $H^0(Y, \o_Y)^\times$
iff $a/c=a'/c'$. 
\end{say}

\begin{defn}\label{full.supp.sects.say}
  Let $Y$ be a  reduced scheme,
$B\subsetneq Y$ a closed subset and $L$ a line bundle on $Y$. 
For $m>0$ set
$$
\Gamma^B(Y,L,m):=\{s\in H^0(Y, L^m): \supp (s=0)=B\}.
\eqno{(\ref{full.supp.sects.say}.1)}
$$
We define analogously  $\Gamma^{\subset B}(Y,L,m)$,
$$
\Gamma^B(Y,L):=\oplus_m \Gamma^B(Y,L^m)\qtq{and} 
\Gamma^{\subset B}(Y,L):=\oplus_m\Gamma^{\subset B}(Y,L,m).
\eqno{(\ref{full.supp.sects.say}.2)}
$$
These are all  unions of $k[Y]^{\times}$-orbits (\ref{k[X].not}).

Note that, in view of (\ref{QD.defn}) and  (\ref{ample.Qiso.top.lem}),  $\Gamma^B(Y,L)$, is a very natural object to consider for us.
\end{defn}

\begin{lem}\label{full.supp.sects.say.3}
  Let $Y$ be a  reduced, projective  scheme,
$B\subset Y$ a closed subset that is  disjoint from $\Sigma(Y)$, and $L$ a line bundle on $Y$. 
\begin{enumerate}
\item If $B$ is irreducible   then  $\Gamma^B(Y,L,m)$ consists of at most one  $k[Y]^{\times}$-orbit. 
\item 
$\Gamma^{\subset B}(Y,L,m)/k[Y]^{\times}$ is  finite.
\item $\Gamma^{\subset B}(Y,L)/k[Y]^{\times}$ is a sub-semigroup of a  finitely generated group.
\end{enumerate}
\end{lem}

Proof.   The first claim follows from (\ref{section.zeros.say}.4). To see the other 2,  
let $B_i\subset B$ be the  irreducible, divisorial components. 
By (\ref{section.zeros.say}.4) 
$\Gamma^{\subset B}(Y,L)/k[Y]^{\times} \subset \oplus_i \z B_i$,
proving (\ref{full.supp.sects.say.3}.3). 
If $s\in \Gamma^{\subset B}(Y,L,m) $  and $(s=0)=\sum_i m_iB_i$ then,
computing the degrees (with respect to some ample divisor) gives that
$\sum_i m_i\deg B_i=\deg D$, hence  $m_i\leq \deg D$ for every $i$
and (\ref{full.supp.sects.say.3}.2) holds. \qed

\medskip

Next we look at the evaluation of a  section of a line bundle $L$ at a point or at a 0-dimensional subscheme $V$.
The twist is that we can not distinguish 2 sections if their zero sets have the same support, and we also can not distinguish various powers of $L$ from each other. Thus for us the outcome of evaluation is  not  a single element of
$H^0(V,\o_V)\otimes L$,  but  a subsemigroup of $H^0(V,\o_V)\otimes \oplus_m L^m$. Our aim is then to understand when this subsemigroup is small.

\begin{defn} \label{ZW.H.link.R.defn}
Let $X$ be a  normal, projective, irreducible  $k$-variety, $L$  a line bundle  and $D$ an effective divisor on $X$. Let  
$W\subset X$ be a closed, integral subvariety  and $V\subset W\setminus D$ a 0-dimensional subscheme.
Set
$$
\begin{array}{lcl}
\res^W_V(D,L,m)&:=&
\im\bigl[\Gamma^{W\cap D}(W, L|_W,m)\to  H^0(V, L^m|_V)\bigr]\qtq{and}\\
\res^W_V(D,L) &:=& \oplus_{m>0} \res^W_V(D,L,m).
\end{array}
$$
Note that $\res^W_V(D,L) $ is a semigroup that is closed under multiplication by $k[W]^{\times}=H^0(W, \o_W)^{\times}$ and, if $D\cap \Sigma(W)=\emptyset$, then
 $\res^W_V(D,L)/k[W]^\times $ is a subsemigroup of a  finitely generated group by (\ref{full.supp.sects.say.3}.3).
\end{defn}

The following elementary observations turn out to be crucial.

\begin{prop} \label{ZW.H.link.prop.2} Using the notation and assumptions of  (\ref{ZW.H.link.R.defn}), 
assume also that  $W$ is irreducible, $D:=(s=0)$ for some $s\in H^0(X, L^r)$, 
$D\cap W$ is irreducible and disjoint from  $\Sigma( W)$. Then 
\begin{enumerate}
\item $\res^W_V(D,L)\cong \langle s\rangle_{\q}\cdot  k[W]^{\times}$, 
where $\langle s\rangle_{\q} $ denotes the saturation of the subsemigroup generated by powers of $s$. 
\item $\res^W_V(D,L) $
 depends only on
$D,L$ and $V$  (but not  on $W$)
in the following cases
\begin{enumerate}
\item   $k[W]=k$, 
\item  $k[W]/k$ is Galois and $V$ is irreducible, or
\item $k[W]=k[\red V]$ and it is separable over $k$.
\end{enumerate}
\end{enumerate}
\end{prop}

Proof.  The first  assertion follows from (\ref{full.supp.sects.say.3}.1).
Indeed,  $s|_W$ is the unique section of $L^r|_W$ (up to $k[W]^{\times}$) that defines  $\supp(D\cap W)$. Therefore   
$$
\res^W_V(D,L,rm)= s^m|_W\cdot k[W]^{\times}|_V= s^m|_V\cdot k[W]^{\times}|_V.
\eqno{(\ref{ZW.H.link.prop.2}.3)}
$$
For other values, $\res^W_V(D,L,m')$ is either empty or 
consists of a single  $k[W]^{\times}$-orbit. We can write an element of it as
$(s|_W)^{m'/r}$, though it is well defined only up to an $r$th root of unity. 

If (\ref{ZW.H.link.prop.2}.2.a) holds then $k[W]^{\times}=k^{\times}$ and the
$k^{\times}$-action on $H^0\bigl(V, L^m|_V\bigr)$
is independent of $W$. 

In general the image of the restriction map 
$\sigma: k[W]\to k[V]$
depends on $W$, but not much. If $\red V=\{p_i:i\in I\}$
then we have finitely many choices for each embedding
$k[W]\into  k(p_i)$, giving finitely many possibilities for
$$
\red(\sigma):k[W]\to \oplus_i k(p_i)=k[\red V].
$$
(\ref{ZW.H.link.prop.2}.2.a--c) are the cases when there is a unique 
$k[W]\to \oplus_i k(p_i)$. 
If  $k[W]/k$ is a separable extension, then    
$\red(\sigma)$ has a unique lifting to
$\sigma:k[W]\to k[V]$ by (\ref{sep.residue.field.lem}.3). \qed
\medskip

Note that $\res^W_V(D,L) $ is nonempty only if 
$D$ is the support of some section of $L^m$ for some $m$.
Thus (\ref{ZW.H.link.prop.2}) is meaningful only if there are sections
with irreducible support. For most applications we need many such sections. This leads us to the following definition.

\medskip
{\bf Bertini-Hilbert dimension}
\medskip

\begin{say}\label{BH.defn.intro.say}
Let $X$ be a projective  variety and $L$ an ample line bundle on $X$.
We are looking for sections  $s\in H^0(X, L)$ that satisfy 3  properties.
\begin{enumerate}
\item The zero set $(s=0)$ is irreducible.
\item The values of $s$ at some points $x_i\in X$ are specified. More generally, given a 0-dimensional subscheme $Z\subset X$, we would like to specify $s|_Z$.
\item The zero set $(s=0)$  avoids a finite set of points $\Sigma\subset X$.
\end{enumerate}
To formalize these,  let $X$ be a scheme  over a field $k$,  $Z\subset X$ a  subscheme,  $L$  a line bundle on $X$ and  
 $s_Z\in H^0(Z, L|_Z)$. Set 
$$
H^0(X, L, s_Z):=\bigl\{s\in H^0(X, L):    s|_Z=c s_Z\mbox{ for some } c\in H^0(X, \o_X)\bigr\}.
\eqno{(\ref{BH.defn.intro.say}.4)}
$$
This is a vector subspace of $H^0(X, L)$. 
If $X$ is integral then 
for the corresponding linear systems we use the notation  $|L, s_Z|\subset |L|$. 
For a finite subset $\Sigma\subset X$, let  
$$
|L, s_Z, \Sigma^c|:=\{D\in  |L, s_Z|:  D\cap \Sigma=\emptyset\}
\eqno{(\ref{BH.defn.intro.say}.5)}
$$ 
denote the subset of those divisors that are disjoint from $\Sigma$.  Finally  we use 
$$
|L, s_Z, \Sigma^c|^{\rm irr}:=\{D\in  |L, s_Z, \Sigma^c|:  D \mbox{ is irreducible.} \}
\eqno{(\ref{BH.defn.intro.say}.6)}
$$ 
For our purposes we are free to replace $L$ by $L^m$. 
Thus  $H^0(X, L^m)\to H^0(Z, L^m|_Z)$ is surjective  
and $|L^m, s^m_Z|$ is  very ample on $X\setminus Z$ for $m\gg 1$. 
Hence conditions  (\ref{BH.defn.intro.say}.2--3) are easy to satisfy and
 the key question is the irreducibility condition (\ref{BH.defn.intro.say}.1).
Next we discuss 3 cases when we can guarantee irreducibility.
\end{say}

The optimal situation is when $\dim X\geq 2$.

\begin{lem}\label{bert.f.basic.prop.lem} Let $X$ be an irreducible,  projective variety of dimension $\geq 2$ over an infinite field $k$. Let $\Sigma\subset X$ be a finite subset and $Z\subset Z$ a finite subscheme. Let $L$ be an ample  line bundle on $X$ and  $s_Z\in H^0(Z, L|_Z)^\times$. Then, 
\begin{enumerate}
\item $|L^m, s^m_Z, \Sigma^c|^{\rm irr}$ contains an open and dense subset of $ |L^m, s^m_Z|$ for $m\gg 1$.
\end{enumerate}
\end{lem}

Proof.  $ |L^m, s^m_Z|$ is very ample on  $X\setminus Z$, hence
this follows from the usual Bertini theorems (\ref{gen.ci.say}).
\qed

\medskip

Next we consider Hilbertian  fields  (\ref{hilb.field.say}).
Here $|L^m, s^m_Z, \Sigma^c|^{\rm irr}$ need not be open, but it is still quite large. By (\ref{bert.f.basic.prop.lem}), we need to pay attention only to curves.

\begin{lem}\label{hilb.f.basic.prop.lem}  Let $C$ be an irreducible,    projective curve over a Hilbertian  field $k$. Let $\Sigma\subset C$ be a finite subset and $Z\subset C$ a finite subscheme. Let $L$ be an ample  line bundle on $C$ and  $s_Z\in H^0(Z, L|_Z)^\times$. Then, 
\begin{enumerate}
\item $|L^m, s^m_Z, \Sigma^c|^{\rm irr}$ contains the complement of a thin  subset  (\ref{thin.flt.defn}) of $ |L^m, s^m_Z|$ for $m\gg 1$.
\end{enumerate}
\end{lem}

Proof. As before, $ |L^m, s^m_Z,|$ is very ample on  $X\setminus Z$, hence
this follows from a basic property of Hilbertian fields 
(\ref{hilb.field.say}). \qed

\medskip

Analyzing the proofs is Sections~\ref{link.divs.res.sect}--\ref{min.rest.tra.sect}
shows that  a weaker version of (\ref{hilb.f.basic.prop.lem}.1) 
is sufficient. We only need $|L^m, s^m_Z, \Sigma^c|^{\rm irr}$ to be nonempty for some $m>0$. This led to the definition of
{\it weakly Hilbertian}  fields  (\ref{weak.hilb.say}). 
The following is essentially their definition; we state it as a lemma to emphasize the similarity to (\ref{hilb.f.basic.prop.lem}).

\begin{lem}\label{w.hilb.f.basic.prop.lem}  Let $C$ be an irreducible,    projective curve over a weakly Hilbertian  field $k$. Let $\Sigma\subset C$ be a finite subset and $Z\subset C$ a finite subscheme. Let $L$ be an ample  line bundle on $C$. Then
\begin{enumerate}
\item $|L^m, s^m_Z, \Sigma^c|^{\rm irr}$ is nonempty for some $m>0$. \qed
\end{enumerate}
\end{lem}

Note that although we ask for only 1 irreducible divisor, by enlarging
$\Sigma$ we see that we get infinitely many. In fact, the sets
$$
|L^m, s^m_Z, \Sigma^c|^{\rm irr}\subset |L^m, s^m_Z|
$$
seem to be quite large, though we do not have a precise way of stating what this means.
\medskip

For most of the proofs we need to know the smallest dimension where linear systems  are guarateed to have many irreducible  members. This leads to the following definition.

\begin{defn}\label{ber-hil.dim.defn}
Let  $k$ be a field that is not locally finite.
 We define the  {\it Bertini-Hilbert dimension} of  $k$---denoted by $\BH(k)$---by setting   
\begin{enumerate}
\item   $\BH(k)=1$ if $k$ is weakly Hilbertian (\ref{weak.hilb.say}), and
\item   $\BH(k)=2$ otherwise.
\end{enumerate}
In view of (\ref{bert.f.basic.prop.lem}), the distinction is only about curves.
If $k$ is Hilbertian  then $\BH(k)=1$ by (\ref{hilb.f.basic.prop.lem}).

We leave the  definition open for locally finite fields. If $k$ is locally finite and $L$ is an ample line bundle on an irreducible curve
$C$, then every smooth point $p\in C$ is the co-support of some section of
some $L^m$.  This would suggest that $\BH(k)$ should be $1$, but in some applications setting $\BH(\f_q)=2$  or even $\BH(\f_q)=\infty$  would seem the right choice.
\end{defn}

\section{Linear similarity}\label{lin.sim.sec}

\begin{defn}\label{QD.defn} 
Let $X$ be a normal variety.
Two divisors $D_1, D_2$ are {\it linearly similar} 
if there are nonzero integers $m_1, m_2$ such that $m_1D_1\sim m_2 D_2$. 
We denote it by  $D_1\sims D_2$.

If $\qrank\cl(X)=1$ then any 2 effective divisors are linearly similar.
Thus this notion is nontrivial only if $\qrank\cl(X)>1$.

 The set of all effective divisors linearly similar to a  fixed divisor $D$ is naturally an
infinite union of linear systems, we denote it by $|\q D|$.

Let $|\q D|^{\rm irr}\subset |\q D|$  be the  subset parametrizing irreducible (but not necessarily reduced) divisors. 

Some of the linear systems   $|D'|\subset |\q D|$ may be small and behave exceptionally. So let us call a subset $W\subset |\q D|$ {\it stably dense} if
$W\cap |mD|$ is dense in $|mD|$ for $m\gg 1$. 
Note that $|\q D|^{\rm irr}$ need not be dense in $ |\q D|$, but 
 if  $D$ is ample and $\dim X\geq 2$ (more generally, if $D$ is mobile and has Kodaira dimension $\geq 2$) then, by (\ref{gen.ci.say}.3),   $|\q D|^{\rm irr}$
is stably dense in   $ |\q D|$.  

If $\dim X=1$ then $|\q D|^{\rm irr}$ is frequently empty. 
This presents a serious technical difficulty in our treatment.
However, if $\deg D>0$
then $|\q D|^{\rm irr}$ is stably dense in $ |\q D|$
provided   $\BH(k)=1$.

Let $H^0\bigl (X, \o_X(\q D)\bigr)$ denote the  Cox ring of $\q [D]$, that is, the direct sum  of all
$H^0\bigl (X, \o_X(D')\bigr)$ where  we pick one $D'$ from each linear system in 
$|\q D|$. 

If $L=\o_X(D)$ is a line bundle, then  a power of $H^0\bigl (X, \o_X(\q D)\bigr)$
is contained in $H^0\bigl (X, \oplus_mL^m\bigr)$. It is sometimes  convenient 
 to  work with $H^0\bigl (X, \oplus_mL^m\bigr)$.

\end{defn}

\begin{say}[Restriction and linear similarity]\label{rest.linsim.say}
Let $X$ be a variety,  $Z\subset X$ a subvariety and 
$D_1, D_2$  effective divisors on $X$. If $D_1\sims D_2$ then
(aside from some problems that appear for non-Cartier divisors),
   $D_1|_Z\sims D_2|_Z$.

For us the main interest will be the converse: if $D_1|_Z\sims D_2|_Z$, when can we conclude that $D_1\sims D_2$?

Let $D$ be an irreducible divisor. We say that a  subvariety  $Z\subset X$
{\it detects}  linear similarity to $D$ if the following holds.
\begin{enumerate}
\item Let $D'$ be an effective divisor such that $\supp (D\cap Z)=\supp (D'\cap Z)$. Then $D'\sims D$.
\end{enumerate}
It is not always easy to see when this happens, but the following is quite useful.
\medskip

{\it Criterion \ref{rest.linsim.say}.2.}  Assume that $Z\cap \sing X$ has codimension $\geq 2$ in  $Z$, the kernel of $\cl(X)\to \pic(Z\setminus \sing X)$ is torsion, $D$ is disjoint from $\Sigma(Z)$ (\ref{section.zeros.say}.3) and $D\cap Z$ is irreducible. Then $Z$ detects  linear similarity to $D$.
\medskip

Proof.  If $Z\cap \sing X$ has codimension $\geq 2$ in  $Z$ then 
we have a restriction map from rank 1 reflexive sheaves on $X$ (that are 
locally free along  $\Sigma(Z)$) to rank 1 reflexive sheaves on $Z$ (that are 
locally free along  $\Sigma(Z)$) and  such a rank 1 reflexive sheaf on $Z$ is determined  by the divisors of its sections
   (\ref{section.zeros.say}.4). \qed
\medskip

 An easy argument (\ref{surf.detects.lem}) shows 
that this criterion applies  if  $Z$ is a general, ample,  complete  intersection subvariety of dimension $\geq 2$.

The  kernel of $\cl(X)\to \cl(C)$ is also torsion for many ample, complete  intersection curves $C\subset X$, but  $D\cap C$ is essentially never irreducible 
if the base field is algebraically closed.  There are 2 ways to go around this problem. 

If $k$ is  Hilbertian
(\ref{hilb.field.say}) then   $D\cap C$ is irreducible for many choices of $C$.

For arbitrary fields, we need the following.
\medskip

{\it Criterion \ref{rest.linsim.say}.3.}  Let $C\subset X^{\rm ns}$ be a smooth, projective curve.
Assume that the kernel of $\cl(X)\to \pic(C)$ is torsion and the following holds.
\begin{enumerate}
\item[(a)] Let $D\cap C=\{p_1,\dots, p_r\}$. Then the points  $p_1,\dots, p_{r-1}$ are linearly independent over
$ \im [\cl(X)\to \pic(C)]$. More precisely, 
$$
\langle p_1,\dots, p_r\rangle\cap  \im [\cl(X)\to \pic(C)] =\z[D|_C].
$$
 Then $C$ detects  linear similarity to $D$.
\end{enumerate}
\medskip

Proof. Let $D'$ be another effective divisor such that
$\supp (D'\cap C)=\supp(D\cap C)$. Note that both $D,D'$ are Cartier along $C$.
Thus $D|_C=\sum d_i[p_i]$ and $D'|_C=\sum d'_i[p_i]$. By (\ref{rest.linsim.say}.3.a)
$\sum d'_i[p_i]=m_1\sum d_i[p_i]$ for some $m_1$, hence  $D'-m_1D$ is in the kernel of
$\cl(X)\to \pic(C)$. Thus $m_2(D'-m_1D)\sim 0$ for some $m_2>0$. \qed 
\medskip

We see in (\ref{ner.for.ci.thm}) that the condition (\ref{rest.linsim.say}.3.a) holds for many ample-ci curves (\ref{ci.sci.say})  and use this in (\ref{ner.for.ci.thm}) to show that 
 if $k$ is not locally finite, then many
ample-ci curves detect linear similarity to a divisor. 

\end{say}

\begin{lem} \label{surf.detects.lem}  
 Let $X$ be a normal, projective variety  over a field $k$ and $D_1,\dots, D_r$ irreducible divisors on $X$.  Then  linear similarity to the $D_i$ is detected by 
\begin{enumerate}
\item   general, ample, complete intersections  of dimension $\geq 2$,   and
\item   a dense subset of complete  intersection curves if  $k$ is  Hilbertian.
\end{enumerate} 
\end{lem}

Proof. Let $Z\subset X$ be a  general, ample,  complete  intersection surface.
Then $\cl(X)\to \cl(Z)$ is an injection by (\ref{gnl-sri.thm}) and, if $\dim Z\geq 2$, then 
$Z\cap D_i$ is irreducible and reduced for every $i$ by  Bertini's theorem
(\ref{gen.ci.say}.3).

The injectivity of $\cl(X)\to \cl(C)$  also holds  for a dense subset of  complete  intersection curves $C\subset X$ (\ref{ner.for.ci.thm}), but for irreducibility we need
$k$ to be  weakly Hilbertian. So (\ref{rest.linsim.say}.3) applies in both cases. \qed

\medskip

We prove a much stronger form of (\ref{surf.detects.lem}.2)  in 
(\ref{ner.for.ci.thm}), but first we derive some consequences of  (\ref{surf.detects.lem}) and (\ref{ner.for.ci.thm}).


\begin{lem}\label{ample.is.top.lem}
 Let $X$ be a normal, projective variety over a field $k$.
Assume that  
\begin{enumerate}
\item either $k$ is not locally finite and $\dim X\geq 2$,
\item or $k$ is  locally finite and $\dim X\geq 3$.
\end{enumerate}
Then an irreducible divisor  $H$ is $\q$-Cartier and ample iff the following holds.
\begin{enumerate}\setcounter{enumi}{2}
\item For every divisor $D\subset X$ 
and closed points $p,q\in X\setminus D$, there is a divisor  $H(p,q)\subset X$ such that
\begin{enumerate}
\item  $H\cap D=H(p,q)\cap D$,
\item $p\notin H(p,q)$ and 
\item $q\in H(p,q)$.
\end{enumerate}
\end{enumerate}
\end{lem}

Proof.  If $H$ is $\q$-Cartier and ample then the restriction map
$$
H^0\bigl(X, \o_X(mH)\bigr)\onto H^0\bigl(D, \o_D(mH|_D)\bigr)+\o_X(mH)\otimes k(p)+\o_X(mH)\otimes k(q)
$$
is surjective for some $m>0$. We can thus find a section $s(p,q)\in H^0\bigl(X, \o_X(mH)\bigr)$ as needed.

Conversely, by (\ref{surf.detects.lem}) and (\ref{ner.for.ci.thm})  we can choose an ample divisor $D\subset X$ that detects  linear similarity to $H$.
Then assumption (\ref{ample.is.top.lem}.3.a) guarantees that $H(p,q)\sims H$. 
Assumption (\ref{ample.is.top.lem}.3.b) implies that $H$ is $\q$-Cartier at $p$.
Since $p,q$ are  arbitrary points (if we also vary $D$),   $H$ is $\q$-Cartier   and a multiple of it  separates points. \qed

\medskip

Together with (\ref{ample.is.top.lem}), the next result proves  (\ref{main.thm.pf.step.1}).

\begin{lem}\label{ample.Qiso.top.lem}
 Let $X$ be a normal, projective variety  over a field $k$ and $H_1, H_2$ irreducible, $\q$-Cartier,  ample  divisors.
Assume that  
\begin{enumerate}
\item either $k$ is not locally finite and $\dim X\geq 3$,
\item or $k$ is  locally finite and $\dim X\geq 5$.
\end{enumerate}
 Then the following are equivalent.
\begin{enumerate}\setcounter{enumi}{2}
\item   $H_1\sims H_2$.
\item $|\q H_1|^{\rm irr}=|\q H_2|^{\rm irr}$.
\item Let $Z_1, Z_2\subset X$ be 2 disjoint, irreducible subvarieties of dimension $\geq 2$ if $k$ is  locally finite and $\geq 1$ otherwise.
Then there is a $\q$-Cartier,  ample  divisor $H'$ such that
 $\supp(H'\cap Z_i)= \supp(H_i\cap Z_i)$ for $i=1,2$.
\end{enumerate}
\end{lem}

Proof. (\ref{ample.Qiso.top.lem}.3)  $\Leftrightarrow$ (\ref{ample.Qiso.top.lem}.4) is clear. 
If (\ref{ample.Qiso.top.lem}.3) holds then choose $m_1, m_2\gg 1$ such that $m_1H_1\sim m_2H_2$ and
$$
H^0\bigl(X, \o_X(m_1H_1)\bigr)\onto
H^0\bigl(Z_1, \o_X(m_1H_1)|_{Z_1}\bigr)+
H^0\bigl(Z_2, \o_X(m_2H_2)|_{Z_2}\bigr)
$$
is surjective.  We can then find $H'\in |m_1H_1|=|m_2H_2|$ whose restriction to  $Z_i$ is $m_i H_i|_{Z_i}$.

Finally assume (\ref{ample.Qiso.top.lem}.5). By (\ref{surf.detects.lem}) and (\ref{ner.for.ci.thm})  we can choose both $Z, W$ normal, disjoint  and such that they detect  linear similarity to the $H_i$.
Then we have the chain of linear similarities
$$
H_1\stackrel{(by\ Z_1)}{\sims} H'\stackrel{(by\ Z_2)}{\sims} H_2. \qed
$$

\begin{say}[Variants for reducible divisors]
\label{ample.Qiso.top.lem.var}
 With  $X$ as in (\ref{ample.Qiso.top.lem}), let  $B_1, B_2$ be effective   divisors. The above proof shows that
(\ref{ample.Qiso.top.lem}.5) is  equivalent to the following.
\begin{enumerate}
\item  There are $\q$-Cartier, ample, irreducible  divisors $H_1, H_2$ such that
$\supp H_i=\supp B_i$ and   $H_1\sims H_2$.
\end{enumerate}
We can thus recognize irreducible $\q$-Cartier divisors using the following criterion.
\medskip

{\it Claim \ref{ample.Qiso.top.lem.var}.2.} An irreducible  divisor
$D\subset X$ is $\q$-Cartier iff there are irreducible, $\q$-Cartier, ample divisors  $A_1, A_2$ such that   the above criterion holds for $B_1:=D+A_1$ and  $B_2:=A_2$. \qed
\end{say}

\begin{rem}\label{ample.Qiso.top.lem.rem} Using 
(\ref{ample.is.top.lem}) and (\ref{ample.Qiso.top.lem}) we get our first topological invariance claims. Namely,  
 let $X_K,Y_L$ be  normal, projective varieties such that $|X|\sim |Y|$.
Assume that  
 either $L$ is not locally finite and $\dim Y\geq 3$,
 or  $\dim Y\geq 5$.
Then 
\begin{enumerate}
\item  If $X$ is $\q$-factorial then so is $Y$.
\item If $\qrank\cl(X)=1$ then $\qrank\cl(Y)=1$.
\end{enumerate}

Note that by (\ref{quasifinite.field.thm.cor}),   $\p^2_{\f_p}$ is homeomorphic to  smooth surfaces
with arbitrary large Picard number, so some restriction on the dimension is necessary in (\ref{ample.Qiso.top.lem}). 
\end{rem}

\medskip{\bf N\'eron's theorem and consequences}\medskip

\begin{defn}\label{thin.flt.defn}
 Let $X$ be an irreducible variety. 
Following \cite{MR1002324},  a
 subset  $T\subset X(k)$ is 
called {\it thin}  if there is a morphism  $\pi:Y\to X$
such that $T\subset \pi(Y(k))$ and there is no rational section $\sigma:X\map Y$.  This notion is most interesting for
finitely generated, infinite fields. For such fields, 
$\a^1(k)\subset \a^1(k)$ is not thin; this is essentially due to Hilbert. 

A rather typical example to keep in mind is the following.
The map $\a^1\to \a^1$ given by $x\to x^2$ shows that the set of all squares is a thin subset of $\a^1(k)$.  

We also need a version of this for arbitrary fields $K$. 
Let us say that a subset  $T\subset X(K)$ is  {\it field-locally thin}
if for every finitely generated subfield  $k\subset K$, the intersection
$T\cap X(k)$ is thin.
\end{defn}

\begin{thm}\cite[Thm.6]{MR56951}\label{neron.thm}
 Let $k$ be a finitely generated, infinite field. 
Let $U\subset \p^1_k$ be an open subset and $\pi: T_U\to U$  a smooth, projective morphism of relative dimension 1.  
Then there is a dense set
$N(T_U)\subset U(k)$, such that  the restriction map
$$
\pic(T_U)\to \pic(T_u)\qtq{is injective for all $u\in N(T_U)$.}
$$ 
Moreover, $N(T_U)$ contains the complement of a  thin set.  \qed
\end{thm}

A  stronger version is proved in \cite[Thm.C]{MR703488}, though it applies only to number fields and finite extensions of $\f_p(t)$.

\begin{cor} \label{neron.thm.cor}
Let $K$ be a field that is not locally finite.  Let $S$ be a normal, projective surface over $K$ and 
 $|C|=\{C_u: u\in \p^1\}$   a mobile, linear pencil of curves with basepoints
$\{p_1,\dots, p_r\}$. Assume that  a general $C_u$ is smooth 
and $S$ is smooth along it. Let $\{B_j:j\in J\}$ be the irreducible components of the reducible members of $|C|$, plus one of the
irreducible members. 
Let $m_{ij}$ be the intersection multiplicity  of $B_j$ with a general
$C_u$  at $p_i$; this is independent of $u$.

Then there is a dense set  $N(S, |C|)\subset  \p^1(K)$ such that, for $u\in N(S, |C|)$,  
all the linear relations among 
 $$[ p_1(u)], \dots, [ p_r(u)]\in \frac{\pic(C_u)}{\im[\cl(S)\to \pic(C_u)]}
$$ 
are generated by
$\tsum_i m_{ij}[ p_i(u)]=0$ for all $j\in J$.

Moreover, $N(S, |C|)$ contains the complement of a field-locally thin set.

\end{cor} 

Proof. Note that the point $p_i$ is contained in every $C_u$; the notation
$[ p_i(u)]$ indicates that we take its class in $\pic(C_u) $, which depends on $u$.

The restriction of $B_j$ to $C_u$ is 
$\tsum_i m_{ij}[ p_i(u)]$, so we do need to have the equations $\tsum_i m_{ij}[ p_i(u)]=0$.
The interesting part is to show that there are no other relations.

Let $T$ be the normalization of the closure of the graph of
$|C|:S\map \p^1$. The projection $\pi_1:T\to S$ is birational, with exceptional curves $E_i\subset T$ sitting over $p_i$.   
Let $B_j^T\subset T$ denote the birational transform of $B_j$.  
Note that 
$$B_j^T\sim \pi_1^*B_j-\tsum_i m_{ij}E_i.
$$
The second projection $\pi_2:T\to \p^1$ is generically smooth and the irreducible components of its fibers are exactly the $B_j^T$. 

Let   $U\subset \p^1$ be the largest open set over which $\pi_2$ is smooth.
By restriction we get $T_U\to U$. 
The Picard group  of $T_U$ is then
$$
\pic(T_U)=\cl(T)\bigr\slash \langle B_j^T:j\in J\rangle.
$$ 
Choose now a finitely generated subfield $k\subset K$ such that  $S,|C|$, 
the $p_i$ and the $B_j$ are defined over $k$.

Note that $\cl(T_k)=\pi_1^*\cl(S_k)+\tsum_i[E_i]$, and killing 
$\pi_1^*\cl(S_k) $ gives  a natural surjection
$$
\cl(T_k)\bigr\slash \pi_1^*\cl(S_k)
\longrightarrow
\langle E_i:i\in I\rangle\bigr\slash\langle \tsum_i m_{ij}[p_i]:j\in J\rangle.
$$
Thus all the linear relations among 
 $[E_1], \dots, [E_r] \in \cl(T_k)/\pi_1^*\cl(S_k)$
are generated by
$\tsum_i m_{ij}[E_i]=0$ for all $j\in J$.

We now apply (\ref{neron.thm}) to get $N(T_U)\subset \p^1(k)$
such that, for $u\in N(T_U)$,  all the linear relations among 
 $$[ p_1(u)], \dots, [ p_r(u)]\in 
\pic(C_u)\bigr\slash\im[\cl(S_k)\to \pic(C_u)]$$ 
are generated by
$\tsum_i m_{ij}[ p_i(u)]=0$ for all $j\in J$.

This is not exactly what we want since $\cl(S_K)$ may be much bigger than
$\cl(S_k)$. However, if we have a linear relation
$$
\tsum_i \lambda_i[ p_i(u)]\sim [L_K]\qtq{where} L_K\in \im[\cls(S_k)\to \pics(C_u)](K),
$$
then in fact $L_K\in \im[\cls(S_k)\to \pics(C_u)](k)$, hence a power of $L_K$ is
in $\im[\cl(S_k)\to \pic(C_u)]$  (\ref{sch.onto.K.onto.lem}).  Letting $k$ vary now proves our claim. \qed

\begin{thm}\label{ner.for.ci.thm}
 Let $k$ be a field that is not locally  finite. Let
$X$ be a normal, projective variety of dimension $n\geq 2$ over $k$, $\{D_i:i=1,\dots, r\}$  irreducible Weil divisors and $H_1, \dots, H_{n-1}$  ample divisors. Then, for $m_1\gg, m_2,\dots, m_r\gg  1$, there is a dense set
$$
U\subset |m_1H_1|(k)\times\cdots\times |m_{n-1}H_{n-1}|(k)
$$ such that, for $u\in U$, the corresponding
complete intersection curve $C_u$ detects  linear similarity to each $D_i$.
\end{thm}

Proof.   By (\ref{surf.detects.lem}), there is a Zariski open 
$$
U_2\subset |m_2H_2|\times\cdots\times |m_{n-1}H_{n-1}|
$$
such that  $\cl(X)\to \cl(H_2\cap \cdots\cap H_{n-1})$ is an injection for
$(H_2, \dots, H_{n-1})\in U_2$ and the $D_i\cap H_2\cap \cdots\cap H_{n-1}$ are irreducible.   This reduces us  to the case $n=2$. 

Thus from now on we have  a normal, projective surface $X$ over $k$, $\{D_i:i=1,\dots, r\}$  irreducible Weil divisors on $X$ and  an ample divisor $H$ on $X$.

Now choose  a pencil $|C|\subset |mH|$ such that 
\begin{enumerate}
\item $D_0+\cdots+D_r\in |C|$ for some irreducible curve $D_0$,
\item all other members of $|C|$ are irreducible,
\item the general member of $|C|$ is smooth and $X$ is smooth along it.
\end{enumerate}
Applying (\ref{neron.thm.cor}) to it we get a dense set of $|C|(k)$ where
the requirements hold. \qed
\medskip

{\it Remark.} Most likely one can choose  
$U$ such that it  contains the complement of a field-locally thin set. 
\medskip

\section{Linkage of divisors and residue fields}\label{link.divs.res.sect}

\begin{defn}\label{t-link.gen.defn} Let $Y$ be a reduced, projective scheme (not necessarily pure dimensional) and $L$ an ample line bundle on $Y$.  
Two effective divisors $D_1, D_2\subset Y$ are called  {\it (topologically, directly)  $L$-linked} if there is a section $s\in H^0(Y, L^m)$ such that
$\supp (s=0)=\supp D_1\cup \supp D_2$. 
\end{defn}

Next we look at an unusual special case of linkage.

\begin{defn}  \label{t-link.ZW.defn}
Let $X$ be a  normal, projective, irreducible  $k$-variety and
$Z, W\subset X$ closed subsets. Let $L$ be an ample line bundle on $X$. Two effective divisors  $H_Z, H_W\in |\q L|$ are (topologically, directly)  {\it $L$-linked on $Z\cup W$} if  $H_Z\cap Z$ and $H_W\cap W$
are  $L$-linked on $Z\cup W$. That is, 
\begin{enumerate}
\item there is a section $s\in H^0(X, L^m)$ such that
$(s=0)\cap Z=H_Z\cap Z$ and $(s=0)\cap W=H_W\cap W$ (as sets).
\end{enumerate}
Note that (\ref{t-link.gen.defn}) would suggest working with a section
$s'\in H^0(Z\cup W, L^m|_{Z\cup W})$, but a power of $s'$ lifts to a section on $X$, so the 2 versions are equivalent.

As we see below, this notion is not interesting if $Z\cap W=\emptyset$
and it has various problems if $\dim (Z\cap W)\geq 1$. Thus we focus on the case when  $\dim (Z\cap W)=0$.

A key observation is that linkage  carries very significant information about \begin{itemize}
\item  residue fields of  $Z\cap W$ in every characteristic, and 
\item the structure of $Z\cap W$ in characteristic 0. 
\end{itemize}
\end{defn}

\begin{prop} \label{ZW.H.link.prop}
Let $X$ be a  normal, projective, irreducible  $k$-variety and
$Z, W\subset X$ closed subsets
such that $\dim (Z\cap W)=0$. Let $L$ be an ample line bundle  on $X$ and $H_Z, H_W\in |\q L|$.
Then $H_Z, H_W$ are  $L$-linked on $Z\cup W$ iff (using (\ref{ZW.H.link.R.defn}))
$$
\res ^Z_{Z\cap W}(H_Z, L)\cap \res ^W_{Z\cap W}(H_W, L)\neq\emptyset.
$$
\end{prop}

Proof. Assume that  $H=(s=0)$ gives the $L$-linkage for some $s\in H^0\bigl(X, L^m\bigr)$.
Then $s|_Z\in \Gamma^{Z\cap H_Z}\bigl(Z,  L|_Z,m\bigr)$ and
$s|_W\in\Gamma^{W\cap H_W}\bigl(W,  L|_W,m\bigr)$ have the same
restriction to $Z\cap W$. 

Conversely, if 
$s_Z\in H^0\bigl(Z,  L^m|_Z\bigr)$ and
$s_W\in H^0\bigl(W,  L^m|_W\bigr)$ have the same image in  
$H^0\bigl(Z\cap W,  L^m|_{Z\cap W}\bigr)$,
then they glue to a section
$s_{Z\cup W}\in H^0\bigl(Z\cup W,  L^m|_{Z\cup W}\bigr)$, 
 and then   $s_{Z\cup W}^{m'}$ lifts to a section of $H^0\bigl(X, L^{m'm}\bigr)$ for some $m'>0$. \qed
\medskip

\begin{say}\label{ZW.asymmetrty.say}
 The conditions in (\ref{ZW.H.link.prop})  give the strongest restriction if 
$$
\res ^Z_{Z\cap W}(H_Z, L)/k^{\times}\qtq{and}
\res ^W_{Z\cap W}(H_W, L)/k^{\times}
\eqno{(\ref{ZW.asymmetrty.say}.1)}
$$ both have $\q$-rank 1.
However, in general these objects are essentially extensions of a finite rank semigroup by $k[Z]^{\times}$  (resp. $k[W]^{\times}$).
Criteria for (\ref{ZW.asymmetrty.say}.1) are given in
(\ref{full.supp.sects.say}.2). This case can be used to identify the $k$-points of $X$  (\ref{ZW.H.link.prop.1p.cor.2}). 

We  get   further interesting consequences if we relax these restrictions. 
The general situation seems rather complicated. 
In our applications it is  advantageous to
work with a non-symmetric situation:
\begin{enumerate}\setcounter{enumi}{1}
\item $H^0(Z, \o_Z)=k$, and
\item $\res ^W_{Z\cap W}(H_W, L)/k[W]^{\times}$ has $\q$-rank 1.
\end{enumerate}
Note that  (\ref{ZW.asymmetrty.say}.2) holds if $Z$ is geometrically connected and reduced. In applications we achieve this by choosing   $Z$ to be  ample-ci  (\ref{ci.sci.say}). 

We saw in (\ref{ZW.H.link.prop.2}) that (\ref{ZW.asymmetrty.say}.3) holds 
if  $\dim W\geq \BH(k)$ (with  some additional mild genericity conditions).

\end{say}

Next we study the case when linking is always possible.

\begin{defn}[Free linking]\label{ZW.H.link.defn}
Let $X$ be a  normal, projective, irreducible  $k$-variety and
$L$  an ample line bundle  on $X$.
Let $Z, W\subset X$ be closed, integral  subvarieties
such that $\dim (Z\cap W)=0$.   

We say that {\it $L$-linking is free} on $Z\cup W$ iff
  two divisors $H_Z, H_W\in |\q L|$ are $L$-linked on $Z\cup W$ whenever they  \
are disjoint from $\Sigma(Z\cup W)$ (\ref{section.zeros.say}.3).
\end{defn}

In the rest of the section we discuss various cases when  the topological notion of free linking makes it possible to obtain information about the residue fields of closed points. 
The key is the following.

\begin{conj}\label{ZW.H.link.defn.conj}
Let $X$ be a  normal, projective, irreducible  $k$-variety and
$L$  an ample line bundle  on $X$.
Let $Z, W\subset X$ be closed, integral, positive dimensional  subvarieties
such that $\dim (Z\cap W)=0$.   The following are equivalent.
\begin{enumerate}
\item $L$-linking is  free  on $Z\cup W$. 
\item $k[Z\cap W]^\times\bigl\slash k[Z]^\times\cdot k[W]^\times$ is a torsion group.
\item One of the following holds.
\begin{enumerate}
\item $\chr k=0$,  $Z\cap W$ is  reduced, and \\
either $ k[Z\cap W]=k[Z]$ or $k[Z\cap W]=k[W]$.
\item $\chr k> 0$, and   either  $k\bigl[\red(Z\cap W)\bigr]/k[Z]$ 
or $k\bigl[\red(Z\cap W)\bigr]/k[W]$ is purely inseparable.
\item  $k$ is locally finite.
\end{enumerate}
\end{enumerate}
\end{conj}

Comments. Note  that (\ref{ZW.H.link.defn.conj}.2) $\Rightarrow$ (\ref{ZW.H.link.defn.conj}.1) holds by (\ref{ZW.H.link.prop})
and the equivalence of (\ref{ZW.H.link.defn.conj}.2) and (\ref{ZW.H.link.defn.conj}.3) is proved in  (\ref{2.fields.mult.gen.lem}).

Next we show that (\ref{ZW.H.link.defn.conj}.1) $\Rightarrow$ (\ref{ZW.H.link.defn.conj}.3) holds if $W$ is geometrically connected and  $\dim W\geq \BH(k)$. A careful study of the proof shows that the first assumption is not necessary, and the validity of (\ref{indep.zeros.conj})  would imply 
that (\ref{ZW.H.link.defn.conj}.1) $\Rightarrow$ (\ref{ZW.H.link.defn.conj}.2) holds unconditionally.

\begin{prop} \label{ZW.H.link.prop.cor.v1}
Let $X$ be a  normal, projective, geometrically irreducible  $k$-variety and
$L$  an ample line bundle  on $X$.
Let $Z, W\subset X$ be closed, integral, positive dimensional  subvarieties
such that $\dim (Z\cap W)=0$.   
Assume that   $W$ is geometrically connected and  $\dim W\geq \BH(k)$.

Then   (\ref{ZW.H.link.defn.conj}.1--4) are equivalent. \end{prop}

Proof.  As we noted, we need to show that (\ref{ZW.H.link.defn.conj}.1) $\Rightarrow$ (\ref{ZW.H.link.defn.conj}.3).
Thus assume (\ref{ZW.H.link.defn.conj}.1). 
First choose   $s_Z\in H^0(X, L^m)$ such that
$\supp(s_Z=0)=\supp H_Z$ is  disjoint from $\Sigma(Z\cup W)$.

By (\ref{bert.f.basic.prop.lem}) and (\ref{w.hilb.f.basic.prop.lem}), for every  $s_{Z\cap W}\in k[Z\cap W]^\times$,  
there is an $n>0$ such that 
 $ s_{Z\cap W}^n$ extends to 
 $s_W\in H^0(X, L^n)$ and 
 $W\cap H_W$ is irreducible, where  $H_W:=\supp(s_W=0)$. 

If $H_Z, H_W$ are $L$-linked, then there is an
$s\in H^0(X, L^m)$ as in (\ref{t-link.ZW.defn}). 
(We can use the same $m$, if we pass to  suitable powers of $s,s_Z, s_W$.)



By (\ref{full.supp.sects.say.3}), there are 
 $u_Z\in k[Z]^\times$,  $u_W\in k[W]^\times=k^\times$, a finitely generated subgroup  $\Gamma_Z\subset k(Z)^\times$,  $\gamma_Z\in \Gamma_Z$ and a  natural number $ r$ such that
$$
s_W^r=s^r|_W\cdot u_W\qtq{and} 
s_Z^r=s^r|_Z\cdot u_Z\cdot\gamma_Z.
\eqno{(\ref{ZW.H.link.prop.cor.v1}.1)}
$$ 
Therefore
$$
\begin{array}{rcl}
s_W^r|_{Z\cap W}\cdot s_Z^{-r}|_{Z\cap W}&=&u_W|_{Z\cap W}\cdot u_Z^{-1}|_{Z\cap W}\cdot \gamma_Z^{-1}|_{Z\cap W}\\[1ex]
&\in& k^\times\cdot k[Z]^\times\cdot  \Gamma_Z|_{Z\cap W}.
=k[Z]^\times\cdot \Gamma_Z|_{Z\cap W}.
\end{array}
\eqno{(\ref{ZW.H.link.prop.cor.v1}.2)}
$$
Next note that  $ s_W|_{Z\cap W}= s_{Z\cap W}^n$ where
$s_{Z\cap W}$ is arbitrary.  Thus $s_{Z\cap W}^ns_Z^{-1}|_{Z\cap W}$ is an arbitrary
element of $k[Z\cap W]^\times $ (up to $n$-torsion). Therefore 
 we get that 
$$
k[Z\cap W]^\times\big\slash \bigl(k[Z]^\times\cdot \Gamma_Z|_{Z\cap W}\bigr)
\qtq{is torsion.}
\eqno{(\ref{ZW.H.link.prop.cor.v1}.3)}
$$
Thus  we obtain that
$$
k[Z\cap W]^\times/k[Z]^\times
\qtq{has finite $\q$-rank.}
\eqno{(\ref{ZW.H.link.prop.cor.v1}.4)}
$$
By (\ref{ker.fin.Qrank.lem.2}) we are in one of 4 cases.
\begin{enumerate}\setcounter{enumi}{4}
\item   $k$ is locally finite; giving (\ref{ZW.H.link.defn.conj}.4.c).
\item $\chr k>0$  and $k[Z]\into k\bigl[\red(Z\cap W)\bigr]$ is
 a purely inseparable extension; giving (\ref{ZW.H.link.defn.conj}.4.b).
\item  $\chr k=0$  and $k[Z]\cong k\bigl[Z\cap W]$; giving (\ref{ZW.H.link.defn.conj}.4.a).
\item $\deg  (k/\q)<\infty$.
\end{enumerate}
 In the latter case
 $k$ is  Hilbertian.  Once $k$  is  Hilbertian,   at the beginning of the proof we can choose  $Z\cap H_Z$ to be irreducible; in which case $\Gamma_Z=\{1\}$  by  (\ref{full.supp.sects.say.3}.1).
   Thus (\ref{ZW.H.link.prop.cor.v1}.3) becomes
$$
k[Z\cap W]^\times\big\slash k[Z]^\times
\qtq{is torsion,}
\eqno{(\ref{ZW.H.link.prop.cor.v1}.9)}
$$
and 
(\ref{ker.fin.Qrank.lem.2}) implies that $Z\cap W$ is reduced.
\qed


\medskip

Using (\ref{ZW.H.link.prop.cor.v1})
we get a topological way of recognizing $k$-points.

\begin{cor} \label{ZW.H.link.prop.1p.cor.2}
Let $k$ be a perfect
field  that  is not locally finite, and 
 $X$  a  normal, projective, geometrically irreducible  $k$-variety of dimension $>1+\BH(k)$. 
Let $L$ be an ample line bundle and $p\in X$  a closed point. 
Assume that either $\chr k>0$ or $p$ is  a smooth point of $X$. 
The following are equivalent.
\begin{enumerate}
\item $p$ is a $k$-point.
\item  There are integral, ample-sci (\ref{ci.sci.say}) 
subvarieties
$Z, W$   such that 
\begin{enumerate}
\item $\dim Z= 1$,  $\dim W  = \BH(k)$,
\item  $\supp(Z\cap W)=\{p\}$ and 
\item $L$-linking is free on $Z\cup W$. 
\end{enumerate}
\end{enumerate}
\end{cor}

Proof.  (\ref{ZW.H.link.prop.1p.cor.2}.2) $\Rightarrow$ (\ref{ZW.H.link.prop.1p.cor.2}.1) follows from (\ref{ZW.H.link.prop.cor.v1}).
Conversely, we can take $Z, W$ to be general complete intersections of  ample divisors  containing $p$. \qed

\medskip

{\it Remark \ref{ZW.H.link.prop.1p.cor.2}.3.} If $\chr k=0$ and (\ref{ZW.H.link.prop.1p.cor.2}.2) holds then
$Z\cap W$ is a $k$-point, even if $X$ is singular there. However, 
for a singular $k$-point it may not be possible to find $Z, W$  such that $Z\cap W=\{p\}$ (as schemes). 
Thus the method does not yet provide a topological way of identifying singular $k$-points if $\chr k=0$.

\medskip

\begin{cor} \label{ZW.H.link.prop.1p.cor.3}
Let $k$ be a perfect 
field  that  is not locally finite, and 
 $X$  a  normal, projective, geometrically  irreducible  $k$-variety of dimension $>1+\BH(k)$. 
Let $L$ be an ample line bundle and $p,q\in X$   closed points. 
Assume that either $\chr k>0$ or $p$ is  a smooth point. 
The following are equivalent.
\begin{enumerate}
\item There is a $k$-embedding $k(p)\into k(q)$.
\item  There are irreducible  subvarieties
$Z, W$  such that 
\begin{enumerate}
\item  $\dim Z=1, \dim W= \BH(X)$, 
\item $\supp(Z\cap W)=\{p\}$,
\item   $q\in Z$, 
\item $W$ is geometrically connected,  and
\item $L$-linking is free on $Z\cup W$.
\end{enumerate}
\end{enumerate}
\end{cor}

Proof. If (\ref{ZW.H.link.prop.1p.cor.3}.2) holds then 
$k(p)\cong k[Z]$ by (\ref{ZW.H.link.prop.cor.v1}) and (\ref{ZW.H.link.prop.1p.cor.3}.2.c) gives an embedding $k[Z]\into k(q)$.

Conversely, given $k(p)\into k(q)$, the required $Z$ is constructed in (\ref{ZW.H.link.prop.1p.cor.3.lem}) and then choose $W$ to be  a general complete intersection  containing $p$. \qed

\medskip

Reversing the role of $p,q$ we the obtain a criterion to decide whether
$k(p)\cong k(q)$. Note, however, that we get no information about $\deg \bigl(k(p)/k\bigr)$. Nonetheless, we have enough information  to obtain the following.

\begin{say}[Isomorphism of 0-cycles from  $|X|$ and $\simsa$] \label{iso.0-cyc.sinsa.say}
Let $k$ be a perfect field  that  is not locally finite, and 
 $X$  a  normal, projective, irreducible  $k$-variety of dimension $>1+\BH(k)$.   Let $Z_1, Z_2\subset X$ be reduced 0-dimensional subschemes.
Assume that either $\chr k>0$ or $Z_1, Z_2\subset X^{\rm ns}$.

We can then decide, using only $\bigl(|X|,\simsa\bigr)$, whether
 $Z_1, Z_2$ are isomorphic as $k$-schemes.  \qed
\end{say}

\begin{say}[Imperfect fields] \label{ZW.H.link.prop.cor.v1.imp}
If $k$ is an imperfect field, we can apply the above results to $k^{\rm ins}$. This results in the following obvious changes in the statements.

In (\ref{ZW.H.link.prop.1p.cor.2}.1) we characterize  $k^{\rm ins}$-points.

In (\ref{ZW.H.link.prop.1p.cor.3}.1) we characterize  $k$-embeddings $k^{\rm ins}(p)\into k^{\rm ins}(q)$.

In (\ref{iso.0-cyc.sinsa.say}) we characterize isomorphisms
$Z_1\times_kk^{\rm ins}\cong Z_2\times_kk^{\rm ins}$.
\end{say}

\begin{rem} \label{char.p.or.not.rem}
Let  $X$  a  normal, projective, geometrically irreducible  $k$-variety of dimension $>1+\BH(k)$.  Then $\chr k>0$ iff the following holds.
\begin{enumerate}
\item 
There is an integral  curve   $C\subset X$ and a point $p\in C$,  
 such that   $L$-linking is free on $C\cup W$ 
for every ample-sci subvariety $W$ of dimension $\BH(k)$, 
for which  $\supp(C\cap W)=\{p\}$. 
\end{enumerate}
Indeed, if $\chr k=0$ then for any $p\in C$ we can choose $W$ such that
$C\cap W$ is non-reduced, and then  $L$-linking is not free on $C\cup W$
by (\ref{ZW.H.link.prop.cor.v1}).

Conversely, we use (\ref{ZW.H.link.prop.1p.cor.3.lem}) to get
$p\in C$ such that $k(p)^{\rm ins}=k[C]^{\rm ins}$, and then (\ref{ZW.H.link.prop.cor.v1}) and (\ref{ZW.H.link.prop.cor.v1.imp})
show that $L$-linking is  free on $C\cup W$
if $\chr k>0$. 
\end{rem}

\section{Minimally restrictive linking and transversality}\label{min.rest.tra.sect}

\begin{defn} \label{red.int.top.prop.defn}
Let $X$ be a  normal, projective, irreducible  $k$-variety and
$Z, W_1, W_2\subset X$ closed, irreducible, geometrically connected  subvarieties
such that  $\dim (Z\cap W_i)=0$. 
We say that  {\it $L$-linking on  $W_2$  determines $L$-linking on  $W_1$} if 
the following holds.
\begin{enumerate}
\item Let  $H_Z, H_W\in |\q L|^{\rm irr}$ be   effective  divisors disjoint from the $\Sigma(Z\cup W_i)$ such that $W_2\cap H_W$ is irreducible. Then 
$$
\left(
\begin{array}{c}
H_Z, H_W \mbox{ are}\\
\mbox{linked on } Z\cup W_2
\end{array}
\right)
\Rightarrow
\left(
\begin{array}{c}
H_Z, H_W \mbox{ are}\\
\mbox{linked on } Z\cup W_1
\end{array}
\right).
$$
\end{enumerate}
In applying this notion we always assume that $\dim W_i\geq \BH(k)$, hence the above conditions are not empty. 

We say that   {\it $L$-linking is minimally restrictive on  $W_1$,}
if  $L$-linking on  $W_2$  determines $L$-linking on  $W_1$, whenever
$\supp (Z\cap W_1)=\supp (Z\cap W_2)$. 
\end{defn}

The key result---and rationale for the definition---is the following.

\begin{prop}\label{red.int.top.prop}
 Let $k$ be a field of characteristic 0,
 $X$  a  normal, projective, irreducible  $k$-variety  and $L$  an ample line bundle  on $X$. Let $Z, W_1, W_2\subset X$ be closed, integral, geometrically connected  subvarieties
such that  $\dim Z\geq 1$,  $\dim W_i\geq \BH(k)$  and $\dim (Z\cap W_i)=0$. 
The following are equivalent.
\begin{enumerate}
\item $Z\cap W_1\subset Z\cap W_2$  as schemes.
\item  $L$-linking on  $W_2$  determines $L$-linking on  $W_1$.
\end{enumerate}
\end{prop}

Proof.  Pick $H_Z=(s_Z=0)$ and  $H_W=(s_W=0)$.
By assumption  $W_2\cap H_W$ is irreducible and disjoint from
$\Sigma(Z\cup W_2)$.  Thus, by (\ref{ZW.H.link.prop.2}.1), 
$$
\res^{W_2}_{Z\cap W_2}(H_Z, L)= \bigl\langle s_W|_{Z\cap W_2}\bigr\rangle_{\q} \cdot k[W_2]^\times=\bigl\langle s_W|_{Z\cap W_2}\bigr\rangle \cdot k^\times.
$$
So, by (\ref{ZW.H.link.prop}),  $H_Z, H_W$ are linked on $Z\cup W_2$ iff, for some $r>0$, 
$$
s_W^r|_{Z\cap W_2}\in \res^{Z}_{Z\cap W_2}(H_Z, L)
\subset H^0\bigl(Z\cap W_2, \oplus_m L^m|_{Z\cap W_2}\bigr).
\eqno{(\ref{red.int.top.prop}.3)}
$$
Note that (\ref{red.int.top.prop}.1) holds iff there is a  natural surjection
$$
 H^0\bigl(Z\cap W_2, \oplus_m L^m|_{Z\cap W_2}\bigr)\onto 
 H^0\bigl(Z\cap W_1, \oplus_m L^m|_{Z\cap W_1}\bigr).
$$ 
Thus  (\ref{red.int.top.prop}.3) implies that
$$
s_W^r|_{Z\cap W_1}\in \res^{Z}_{Z\cap W_1}(H_Z, L)
\subset H^0\bigl(Z\cap W_1, \oplus_m L^m|_{Z\cap W_1})\bigr),
\eqno{(\ref{red.int.top.prop}.4)}
$$
proving (\ref{red.int.top.prop}.2).


To see the converse, let $N$ be the kernel of
$$
 \rho: H^0\bigl(Z\cap W_1, \o_{Z\cap W_1}\bigr)^\times\to
 H^0\bigl(Z\cap W_1\cap W_2, \o_{Z\cap W_1\cap W_2}\bigr)^\times.
$$
It is a direct sum of a commutative, unipotent group over $k$ and of the
$k(p_i)^\times$ for every $p_i\in W_1\setminus W_2$.  $N$ is positive dimensional
iff   $Z\cap W_1\not\subset Z\cap W_2$.
We distinguish 2 cases, depending on whether 
 $\qrank N(k)=\infty$ or not. 

Next choose any  $\sigma_2\in H^0\bigl(Z\cap W_2, L|_{Z\cap W_2}\bigr)^\times$.
If  $\qrank N(k)=\infty$  then the restriction of $\sigma_2$ to
$Z\cap W_1\cap W_2 $ can be lifted in 2 different ways to 
$$
\sigma_Z, \sigma_W \in H^0\bigl(Z\cap W_1, L|_{Z\cap W_1}\bigr)^\times
$$
such that $\sigma_Z\sigma_W^{-1}$ is non-torsion in 
$$
H^0\bigl(Z\cap W_1, \oplus_mL^m|_{Z\cap W_1}\bigr)^\times\bigr\slash
\res^{Z}_{Z\cap W_1}(H_Z, L).
$$
We can now glue $\sigma_2, \sigma_W$ to a section of
$H^0\bigl(Z\cap (W_1\cup W_2), L|_{Z\cap (W_1\cup W_2)}\bigr)^\times$ and then lift (some power of) it to $s_W\in H^0(X, L)$ such that both $W_i\cap (s_W=0)$ are irreducible and disjoint from $\Sigma(Z\cup W_1\cup W_2)$. Similarly,
we can  glue $\sigma_2, \sigma_Z$ to a section of
$H^0\bigl(Z\cap (W_1\cup W_2), L|_{Z\cap (W_1\cup W_2)}\bigr)^\times$ and then lift (some power of) it to $s_Z\in H^0(X, L)$ such that  $Z\cap (s_Z=0)$ is irreducible and disjoint from $\Sigma(Z\cup W_1\cup W_2)$.

By construction, $s_Z|_{Z\cap W_2}=s_W|_{Z\cap W_2}$, hence  
$(s_Z=0)$ and $(s_W=0)$ are $L$-linked on  $Z\cup W_2$,
but $s_Z|_{Z\cap W_1}$ and $s_W|_{Z\cap W_1}$ are multiplicatively independent, hence  
$(s_Z=0)$ and $(s_W=0)$ are not $L$-linked on  $Z\cup W_2$. 

We are left with the case when $\qrank N(k)<\infty$. In this case
$\deg (k/\q)<\infty$ by (\ref{Q.rank.say}), hence $k$ is Hilbertian (\ref{hilb.field.say}). We can thus choose  $s_Z$ such that  
 $Z\cap H_Z$ is irreducible and disjoint from
$\Sigma(Z\cup W_2)$.  So, by (\ref{ZW.H.link.prop.2}.1), 
$$
\res^{Z}_{Z\cap W_2}(H_Z, L)= \bigl\langle s_Z|_{Z\cap W_2}\bigr\rangle_{\q} \cdot k[Z]^\times=\bigl\langle s_Z|_{Z\cap W_2}\bigr\rangle_{\q} \cdot k^\times.
$$
This implies that $\res^{Z}_{Z\cap W_2}(H_Z, L)$ has trivial intersection with $N$.
We can thus again choose $\sigma_Z, \sigma_W \in H^0\bigl(Z\cap W_1, L|_{Z\cap W_1}\bigr)^\times$
such that $\sigma_Z\sigma_W^{-1}$ is non-torsion, and complete the proof as before. 
\qed

\medskip

\begin{cor}\label{red.int.top.prop.cor}
 Let $k$ be a field of characteristic 0,
 $X$  a  normal, projective, irreducible  $k$-variety and $L$ an ample line bundle. Let
$Z, W\subset X$ be closed, integral, geometrically connected  subvarieties
such that  $\dim (Z\cap W)=0$.
Assume that $\dim X>\dim Z+ \BH(k)$ 

Then  $L$-linking is minimally restrictive on  $Z\cup W$ iff $Z\cap W$ is reduced. 
\end{cor}

Proof. Set $W_1:=W$ in (\ref{red.int.top.prop}) and let $W_2$ run through all 
ample-sci subvarieties of dimension $\BH(k)$  (\ref{ci.sci.say}) that intersect $Z$ exactly along $Z\cap W_1$. 
Then apply   (\ref{red.int.top.prop.lem.1}). \qed

\begin{lem}\label{red.int.top.prop.lem.1}
 Let $X$ be a projective $k$-variety, $Z\subset X$  a subscheme of codimension $>r$ and $P\subset Z$ a reduced, finite subscheme. 
Let  ${\mathcal W}(Z,P)$ be the set of all irreducible, $r$-dimensional, ample-sci (\ref{ci.sci.say}) subvarieties  $W\subset X$ for which
$\supp(Z\cap W)=P$. 
Then 
$$
\bigcap_{W\in {\mathcal W}(Z,P)}  W=P  \qtq{(scheme theoretically).} \qed
$$
\end{lem}

The following consequence of (\ref{red.int.top.prop.cor})  allows us to understand  intersection multiplicities topologically.

\begin{say}[Determining transversality from $|X|$ and $\simsa$]
\label{num.eq.top.say.det}
 Let $k$ be a field of characteristic 0,
 $X$  a  normal, projective, geometrically irreducible  $k$-variety of dimension $>1+\BH(k)$.
Let $H\subset X$ be an irreducible, ample divisor and
$C\subset X$ an irreducible, geometrically connected   curve.
Assume that $C\cap H\subset X^{\rm ns}$.  The following are equivalent.
\begin{enumerate}
\item   All intersections of $C\cap H$ are transversal.
\item  $C\cap H$ is reduced.
\item  $L$-linking is minimally restrictive on  $C\cup H$ for some ample line bundle $L$.
\end{enumerate}
Indeed, (\ref{num.eq.top.say.det}.2) $\Leftrightarrow$  (\ref{num.eq.top.say.det}.3) follows from (\ref{red.int.top.prop.cor}) and
(\ref{num.eq.top.say.det}.1) $\Leftrightarrow$  (\ref{num.eq.top.say.det}.2)
is a bsic property of intersection mutiplicities; see for example \cite[8.2]{MR1644323}. \qed
\medskip

Note that there are several weaknesses of the current form of the above equivalences.

First, we do not yet know how to decide which are the  smooth  points of $X$.  We usually go around this by saying that some assertion holds outside some codimension $\geq 2$ subset.

Second, we also do not yet know how to decide whether a curve $C$ is
 geometrically  connected  or not. However, if $C$ is  ample-sci (\ref{ci.sci.say}),
then $C$ is geometrically  connected (\ref{ci.sci.say}.1).

\end{say}

The above arguments also show the following.

\begin{cor}\label{red.int.top.prop.cor.4}  
 Let $k$ be a field of characteristic 0,
 $X$  a $k$-variety, $Z\subset X$  an irreducible, geometrically connected  subvariety of codimension $r> \BH(k)$ and $p\in Z$ a closed point  such that $X$ is smooth at $p$. 
Then $Z$ is smooth at $p$ iff there is an irreducible, ample-sci subvariety $W\subset X$  of dimension $r$
such that $p\in \supp(Z\cap W)$ and 
$L$-linking is minimally restrictive on  $W$. \qed
\end{cor}



Interchanging the roles of $Z,W$ gives the following dual version.

\begin{cor}\label{red.int.top.prop.cor.6}   
 Let $k$ be a field of characteristic 0,
 $X$  a $k$-variety, $W\subset X$  an irreducible, geometrically connected  subvariety of dimension $r>\BH(k)$ and $p\in W$ a closed point  such that $X$ is smooth at $p$. 
 Then $W$ is smooth at $p$ iff there is an irreducible, ample, complete intersection subvariety $Z\subset X$ 
of codimension $r$
such that $p\in \supp(Z\cap W)$ and 
$H$-linking is minimally restrictive on  $W$. \qed
\end{cor}



\begin{say}\label{red.int.top.prop.p} The    argument in  (\ref{red.int.top.prop}) also applies in positive characteristic, except that then the kernel of
$$
H^0\bigl(Z\cap W, \o_{Z\cap W}\bigr)^{\times}\to H^0\bigl(\red(Z\cap W), \o_{\red(Z\cap W)}\bigr)^{\times}
$$
is $p$-power torsion. Thus multiplicative  independence is not changed
as we pass from $ Z\cap W$ to $\red(Z\cap W) $.  We get that, if $k$ is not locally finite, then the following are equivalent.
\begin{enumerate}
\item $\supp(Z\cap W_1)\subset \supp(Z\cap W_2)$.
\item  $L$-linking on  $W_1$  determines $L$-linking on  $W_2$. 
\end{enumerate}
Thus, while  $L$-linking carries scheme-theoretic information in characteristic 0, it does not in positive characteristic. 
\end{say}

\section{Degrees of curves and divisors}

Next we try to determine intersection numbers $(C\cdot H)$ of curves and divisors. First we use  
(\ref{num.eq.top.say.det}) to check that all 
intersections  are transversal. 
If the field is algebraically closed,  then $(C\cdot H)=\# (C\cap H)$, and we are done. 

If the field is not algebraically closed, we would need to compute
$\deg k(p)/k$ for the intersection points. This we can not do, but, by
(\ref{iso.0-cyc.sinsa.say}),   we can determine when $\deg k(p_1)/k=\deg k(p_2)/k$ for 2 points. 
The following result says that, although  we are not  able to compute
$(C\cdot H)$, we can decide whether $(C_1\cdot H)=(C_2\cdot H)$ for 2  curves. This is enough for many applications.

\begin{thm} \label{secs.isom.zeros.prop.c1}
Let $X$ be a projective variety over an infinite field, $H$ an irreducible, ample, Cartier divisor,  $\{C_i:i\in I\}$ finitely many irreducible, geometrically reduced  curves. The following are equivalent.
\begin{enumerate}
\item  $(C_i\cdot H)$ is independent of $i$. 
\item There is an     irreducible divisor
$G\simsa H$ such that the (scheme theoretic) intersections
$\{C_i\cap G :i\in I\}$ are reduced and  isomorphic to each other. 
\end{enumerate}
Moreover, we can choose $G$ to be   disjoint from any finite subset  $\Sigma\subset X$.
\end{thm}

Proof. Assume that (\ref{secs.isom.zeros.prop.c1}.2)  holds and $G\sim mH$.
Then  $m(C_i\cdot H)=\deg_k (C_i\cap G )$, proving (\ref{secs.isom.zeros.prop.c1}.1).  To see the converse, let $C$ denote the union of the $C_i$. 
Set $L:=\o_X(H)|_C$.
Choose $m\gg 1$ such that  $H^0\bigl(X, \o_X(mH)\bigr)\to  H^0\bigl(C, \o_C(mH|_C)\bigr)$ is surjective and 
 there is a section $s\in H^0\bigl(C, \o_C(mH|_C)\bigr)$ as in (\ref{secs.isom.zeros.prop.c1.cor}). Then  $G:=(s=0)$ works.  \qed

\begin{cor} \label{secs.isom.zeros.prop.c2}
Let $X$ be a normal, projective variety  over an infinite field, $H$ an irreducible, ample,  Cartier divisor,
 $\{D_i:i\in I\}$ irreducible, geometrically reduced  divisors
and $B\subset X$  a closed  subset of dimension $\leq n-2$ containing $\sing X$. 
The following are equivalent.
\begin{enumerate}
\item $(D_i\cdot H^{n-1})$ is independent of $i$. 
\item There are irreducible,    $H$-sci curves $A$ (\ref{ci.sci.say}) that are disjoint from $B$ and such that the (scheme theoretic) intersections
$A\cap D_i$ are reduced and  isomorphic to each other. 
\end{enumerate}
\end{cor}

Proof. As before, (\ref{secs.isom.zeros.prop.c2}.2)  $\Rightarrow$ (\ref{secs.isom.zeros.prop.c2}.1) is clear.
For the converse, we look for  $A$ contained in a 
general complete intersection surface $S\subset X$. 
This reduces us to the special case when $\dim X=2$. Then the $D_i$ are curves, so  (\ref{secs.isom.zeros.prop.c2}) 
follows from (\ref{secs.isom.zeros.prop.c1}). \qed

\medskip

The following observation allows us to describe $\q$-linear equivalence topologically.

\begin{say}[Determining numerical equivalence from $|X|$ and $\simsa$]
\label{num.eq.top.say}
We state 4 pairs of results. In each of them the first is the information we seek, the second shows how it can be obtained using 
(\ref{secs.isom.zeros.prop.c1}--\ref{secs.isom.zeros.prop.c2}) and the previous ones.  We assume that  $X$ is a projective variety over a field
of characteristic 0. The latter is necessary in order to use
(\ref{red.int.top.prop.cor}). 
We use several  claims about intersections of curves and divisors; these are recalled in (\ref{lin.eq.amp.lems}). 

\medskip
{\it Claim \ref{num.eq.top.say}.1.}
Let  $C_1, C_2$ be irreducible curves, not contained in $\sing X$.  Then  
\begin{itemize}
\item  $(C_1\cdot H)=(C_2\cdot H)$  iff
\item  for every finite subset $\Sigma\subset C_1\cup C_2$ 
there is  an  irreducible divisor
$G\sims H$ disjoint from $\Sigma$, such that 
\begin{enumerate}
\item[(a)]   $H$-linking is minimally restrictive on $C_i\cup G$  for $i=1,2$, and
\item[(b)] $C_1\cap G$ and $C_2\cap G$ are isomorphic as $k$-schemes.\qed
\end{enumerate}
\end{itemize}

\medskip
{\it Claim \ref{num.eq.top.say}.2.}
 Let  $C_1, C_2$ be irreducible curves, not contained in $\sing X$. Then  
\begin{itemize}
\item  $C_1\equiv C_2$ iff
\item   $(C_1\cdot H)=(C_2\cdot H)$  for every ample  divisor $H$. \qed
\end{itemize}

\medskip
{\it Claim \ref{num.eq.top.say}.3.}
Let  $D_1, D_2$ be irreducible, geometrically connected  divisors on $X$ and $H$ an ample divisor. 
Then  
\begin{itemize}
\item  $(D_1\cdot H^{n-1}) = (D_2\cdot H^{n-1})$ iff 
\item  For every  codimension $\geq 2$ subset $B\subset X$ 
there are irreducible,   $H$-sci  curves $A$  that are disjoint from $B$ and such that
\begin{enumerate}
\item[(a)]   $H$-linking is minimally restrictive on $A\cup D_i$   for $i=1,2$, and
\item[(b)] $A\cap D_1$ and $A\cap D_2$ are isomorphic as $k$-schemes.\qed
\end{enumerate}
\end{itemize}

\medskip
{\it Claim \ref{num.eq.top.say}.4.}
Let  $D_1, D_2$ be irreducible, geometrically connected  divisors on $X$. 
Then  
\begin{itemize}
\item  $D_1\equiv D_2$ iff
\item  $(D_1\cdot H^{n-1}) = (D_2\cdot H^{n-1})$ for every ample  divisor $H$. \qed
\end{itemize}

\medskip
{\it Claim \ref{num.eq.top.say}.5.}
Let  $H_1, H_2$ be irreducible,  ample  divisors.
Then  
\begin{itemize}
\item  $H_1\simq H_2$ iff
\item  $H_1\simsa H_2$ and  $(H_1\cdot H_2^{n-1})=(H_2\cdot H_2^{n-1})$. 
\qed
\end{itemize}

\end{say}


\begin{lem}  \label{lin.eq.amp.lems}
Let $X$ be a normal, projective variety.
\begin{enumerate}
\item  Let $C_1, C_2$ be 1-cycles. Then  $C_1\equiv C_2$ iff   $(C_1\cdot H)=(C_2\cdot H)$ for every ample Cartier divisor $H$. 
\item Let $D_1, D_2$ be $\q$-Cartier divisors. Then  $D_1\equiv D_2$ iff   $(D_1\cdot A)=(D_2\cdot A)$ for every irreducible, ample-sci curve  (\ref{ci.sci.say}).  
\end{enumerate}
\end{lem}

Proof. We use  the theory of cones and divisors; see for example \cite[1.4.C]{laz-book}. 
Let $N^1(X)=\ns(X)\otimes\q$ be the N\'eron-Severi space (the space of $\q$-Cartier $\q$-divisors  modulo algebraic equivalence) and $N_1(X)$
the space of curves with $\q$-coefficients, modulo  algebraic equivalence.
By Kleiman's theorem \cite[1.4.23]{laz-book} they are dual to each other, hence they have the same dimension.  Ample divisors span $N^1(X)$, thus (\ref{lin.eq.amp.lems}.1) holds. For (\ref{lin.eq.amp.lems}.2) it remains to show that irreducible, ample-sci curves span $N_1(X)$.

First we claim that if $H$ is ample then
$H^{n-2}:N^1(X)\to N_1(X)$ is an isomorphism. This is a special case of the Hard Lefschetz theorem  \cite[3.1.39]{laz-book}, which is usually stated for
smooth varieties. Since $N^1(X)$ and $ N_1(X)$ have the same dimension, it is enough to show that $H^{n-2}$ gives an injection.  This follows from the
Grothendieck-Lefschetz hyperplane theorem: 
if $S\subset X$ is a general $H$-complete intersection surface then
$\pic(X)\to \pic(S)$ is an injection.

Finally, let $D$ be any Cartier divisor. Then
$H+\epsilon D$  is ample for $0<\epsilon\ll 1$ and
$$
H^{n-2}D\in \bigl\langle (H^iD^{n-1-i}) : i=0,\dots, n-1\bigr\rangle =
\bigl\langle (H+\epsilon D)^{n-1} : 0<\epsilon\ll 1\bigr\rangle. \qed
$$

\medskip
{\bf Zero cycles on curves}
\medskip

Let $L$ be a very ample line bundle on a reduced, projective curve over an algebraically closed field. The zero set of a general section
of $L$ consist of $\deg_CL$ distinct points.
However, if we work over a non-closed field $k$, then 
the zero set of a general section $s\in H^0(C, L)$ is a union of points of the form
$\spec_k k_i$ for some field extensions $k_i/k$, that depend on the choice of the section in a rather unpredictable way. 
We may thus aim to find sections $s\in H^0(C, L)$ whose zero set is arithmetically simple.  If $k$ is Hilbertian, we can chose the zero set to be irreducible. 
Another direction would be to find zero sets that consist of low degree points.
This is, however, impossible
already for genus 1 curves over $\q$.

The next result says that, for any finite set of curves $C_i$ and line bundles
$L_i$, a uniformly nice choice of sections is possible.

\begin{thm} \label{secs.isom.zeros.prop.c1.cor}
Let $C$ be a  geometrically reduced, projective   curve over   a field   $k$ 
with  irreducible components $\{C_i: i\in I\}$.
Let $L$ be an ample line bundle on $C$ and $\Sigma\subset C$ a finite set.
 Then there is an $m>0$, a separable field extension $K/k$ 
and a section $s\in H^0(C, L^m)$ such that 
\begin{enumerate}
\item $(s=0)$ is disjoint from $\Sigma\cup \sing C$, and
\item  $C_i\cap (s=0)$ is  isomorphic to the disjoint  union of $\bigl(m/\deg (K/k)\bigr)\cdot \deg_{C_i}L $ copies of  $\spec_kK$ for every $i$. 
\end{enumerate}
\end{thm}

Proof.   For $m_1$ large enough  there is a 
 separable morphism  $\pi:C\to \p^1$ such that $L^{m_1}\cong \pi^*\o_{\p^1}(1)$.
By (\ref{secs.isom.zeros.prop.v2})  there is  a separable point $p\in \p^1$ 
that is disjoint from $\pi(\Sigma\cup \sing C)$ and 
such that 
 $\pi^{-1}(p)$ is a reduced, disjoint union of  copies of $p$.
Let $s'\in H^0(\p^1, \o_{\p^1}(m_2))$ be a defining equation of $p$.
Then $s:=\pi^*s'\in H^0(C, L^{m_1m_2})$ has the required properties. \qed

\begin{cor} \label{secs.isom.zeros.prop.c1.cor.cor}
Let $X$ be a projective variety over a field $k$, 
$L$ an ample line bundle, $\{C_i: i\in I\}$ a finite set of  geometrically reduced  curves
and $\Sigma\subset X$  a finite subset. 
 Then there is an $m>0$, a section  $s\in H^0(X, L^m)$ and a separable field extension $K/k$  such that 
\begin{enumerate}
\item $(s=0)$ is disjoint from $\Sigma\cup \sing (\cup_i C_i)$, and
\item  $C_i\cap (s=0)$ is  isomorphic to the disjoint  union of $\bigl(m/\deg (K/k)\bigr)\cdot\deg_{C_i}L $  copies of  $\spec_kK$ for every $i$.  \qed
\end{enumerate}
\end{cor}

\begin{say} Another consequence of  (\ref{secs.isom.zeros.prop.c1.cor}) is the following.  
Let $X$ be a projective variety over a field $k$. Then every element of the Chow group of 0-cycles of degree 0 has a representative of the form
$$
\tsum_i  \bigl([p_i]-[q_i]\bigr)\qtq{where}  k(p_i)\cong k(q_i).
$$
\end{say}

The following can be viewed as a generalization of a special case of Chebotarev's density theorem: there are infinitely many completely split primes in any separable field extension.

\begin{prop}\cite{MR1851662}  \label{secs.isom.zeros.prop.v2}
Let  $C$  be a geometrically reduced $k$-curve
and $\pi:C\to \p^1$ a quasi-finite, separable  morphism. Then there are  infinitely many  separable points $p_j\in \p^1$ such that 
 $\pi^{-1}(p_j)$ is a reduced, disjoint union of  copies of $p_j$ for every $j$.
\end{prop}

Proof.   Let $C_i$ be the  irreducible components of $C$. 
Let $D$ be the normalization of $\p^1$ in the Galois closure of a composite of the $k(C_i)/k(\p^1)$.  If $p_j$ works for $\sigma: D\to \p^1$ and 
$\pi$ is \'etale over $p_j$  then  $p_j$ works for $C\to \p^1$.

If $k$ is infinite, then a general pencil of very ample divisors gives
  a  separable morphism $\rho:D\to {\mathbf P}^1:=\p^1$ such that
$$
(\rho, \sigma):D\to D'\subset {\mathbf P}^1\times \p^1
\eqno{(\ref{secs.isom.zeros.prop.v2}.1)}
$$
is birational onto its image. (We use the notation ${\mathbf P}^1$ to distinguish the 2 factors.)  Let $S\subset D$ be the union of the preimage of $\sing D'$,  the ramification locus of $\sigma$ and the ramification locus of $\rho$.

Pick any   $c\in {\mathbf P}^1(k)\setminus \rho(S)$. Let $p_D\in \rho^{-1}(c)$ be any closed point, $p'_D$ its image in $D'$ and $p:=\sigma(p_D)\in \p^1$. 
Then 
$$
k(p_D)=k(p'_D)=k(c)\otimes_k k(p)=k(p).
$$
Since $D/\p^1$ is Galois, the same holds for all points in
$\sigma^{-1}(p)$. 

If $k$ is finite, choose $q=p^r$ such that $D$ decomposes into $m$ 
irreducible components that are geometrically  irreducible.
Then $D$ has  about $mq$ points in $\f_q$.  All these map to
$\f_q$ points in $\p^1$. We show that for most of them, their image is not  defined over a subfield of $\f_q$. All subfields of $\f_q$ have at most $\sqrt{q}$ elements and the number of maximal ones equals the  number of prime divisors of  $r$, so there are at most $\log_2 r$ of them.
Thus at most $\log_2 r \cdot \sqrt{q} $ points of $\p^1(\f_q)$ are in a smaller subfield and these have at most $\deg \pi\cdot \log_2 r \cdot \sqrt{q} $
preimages in $D$. So for $q\gg 1$, almost all $\f_q$ points of $D$ map to points of $\p^1$ whose residue field is $\f_q$. \qed

\begin{say}[Variants of  (\ref{secs.isom.zeros.prop.v2})] \label{secs.isom.zeros.prop.v2.var}  The following versions are also useful.
\medskip

{\it Claim \ref{secs.isom.zeros.prop.v2.var}.1.} Let  $X, Y$  be  geometrically reduced $k$-schemes
and $\pi:X\to Y$ a quasi-finite, separable  morphism. Then there is a Zariski dense set of  separable points $p_j\in Y$ such that 
 $\pi^{-1}(p_j)$ is a reduced, disjoint union of  copies of $p_j$.
\medskip

Proof. We can  replace $Y$ by a general curve $B\subset Y$.  The curve case is reduced to  (\ref{secs.isom.zeros.prop.v2})  by composing with a  quasi-finite, separable  morphism  $B\to \p^1$. \qed

\medskip

{\it Claim \ref{secs.isom.zeros.prop.v2.var}.2.} 
Let  $X, Y$  be  geometrically reduced $k$-schemes
and $\pi:X\to Y$ a quasi-finite, separable  morphism. 
Then there are infinitely  many irreducible  divisors  $D_j\subset Y$ such that 
 $\pi^{-1}(D_j)$ is a reduced union of  divisors $D_j^i$, such that, each
$D_j^i\to D_j$ is birational. 
\medskip

Proof.  As in (\ref{secs.isom.zeros.prop.v2}) and in (\ref{secs.isom.zeros.prop.v2.var}.1), we may assume that
$Y=\p^n$ and $X\to \p^n$ is Galois.
 If $k$ is infinite,
the proof works as in (\ref{secs.isom.zeros.prop.v2}),  but we replace (\ref{secs.isom.zeros.prop.v2}.1) by
$$
(\rho, \sigma): X\to  X'\subset {\mathbf P}^1\times \p^n. 
$$
We can also reduce (\ref{secs.isom.zeros.prop.v2.var}.2) to
(\ref{secs.isom.zeros.prop.v2.var}.1) by choosing a coordinate projection
$p:\p^n\map \p^{n-1}$ and work with the generic fibers of $Y\map \p^{n-1}$ and
$X\map \p^{n-1}$. \qed
\end{say}

\section{Topological pencils}\label{top.pencils.sec}

\begin{defn}[Pencils] \label{pencil.defn}
Let $k$ be a field and $X$  an integral $k$-variety.
Let  $C$  be an integral  curve and $p:X\map C$
a dominant map with  indeterminacy locus  $B\subset X$.
We care about the map $p$ up to birational equivalence, so we may as well assume that $C$ is nonsingular  and projective.

For a  point $c\in C(\bar k)$ let $D_c\subset X$ denote the closure of $p^{-1}(c)$.  In traditional terminology  (see, for example, \cite{zar-41}) 
$$
|D|_{\rm alg}:=\{D_c : c\in C(\bar k)\}
$$ 
is an {\it (algebraic) pencil of divisors} parametrized by the curve $C$.
We use $|D|_{\rm alg}$ to emphasize that  $C$ can be a non-rational curve.
We refer to the map $p:X\map C$ itself as an {\it (algebraic) pencil.}
The pencil is called {\it linear}
if $C\cong\p^1_k$. These are  usually denoted by $|D|$.

If $p$ factors as  $X\map C'\map C$ then $X\map C'$ determines another pencil $|D'|_{\rm alg}$ and members of  $|D|_{\rm alg}$ are certain unions of members of $|D'|_{\rm alg}$. If $C'\map C$ is not birational, we say that   $|D|_{\rm alg}$  is  {\it  composite with}  
$|D'|_{\rm alg}$.

A theorem of Bertini---which seems to have been first fully proved by \cite{vdW-36}---states that the following are equivalent.
\begin{enumerate}
\item Almost all members of $|D|_{\rm alg}$ are  irreducible and generically reduced.
\item  $|D|_{\rm alg}$  is  not  composite with any other pencil.
\item  $k(C)$ is algebraically closed in $k(X)$.
\end{enumerate}
 As in \cite{zar-41}, we call such a pencil {\it non-composite.}  

Every pencil is composite with a unique non-composite pencil;
nowadays this is  pretty clear 
 using Stein factorization.

Note. The  definition of linear system frequently allows fixed components; those without fixed components are called {\it mobile.} In this terminology,
 our pencils are the {\it mobile pencils.}
\end{defn}

\begin{defn}[Topological pencil]\label{t.pencil.defn}
Let  $X$ be  a  projective, geometrically normal  $k$-variety.
 A {\it t-pencil} is a collection of effective divisors
$\{D_{\lambda}:\lambda\in\Lambda\}$ such that
\begin{enumerate}
\item Every closed point of $X$ is contained in some $D_{\lambda}$.
\item Almost all of the $D_{\lambda}$ are irreducible.
\item  There is a closed subset $B\subset X$ of codimension $\geq 2$ such that
$D_{\lambda}\cap D_{\mu}\subset B$ for every $\lambda\neq  \mu \in \Lambda\}$.
The smallest such $B$ is called the {\it base locus.} 
\item Each $D_{\lambda}\setminus B$ is connected.
\end{enumerate}
We call a t-pencil  {\it ample} if almost all members are ample $\q$-Cartier divisors.
\end{defn}

\begin{exmp} \label{t.penc.true.meb.defn}
Let $X$ be  a  projective, geometrically normal $k$-variety
and  $|D|_{\rm alg}$ a  non-composite algebraic pencil with parameter curve $C$
and base locus $B$. 

For $c\in C(\bar k)$  let $\{D_{c,i}:i=1,\dots, r_c\}$ be the connected components of $\red D_c\setminus B$.   As $\sigma$ runs through all embeddings
$k(c)\into \bar k$, we get $k^{\rm ins}$-divisors
$\cup_{\sigma}  D_{c,i}^\sigma$.
The set of all $D_{c,i}$ forms a t-pencil with base locus $B$.
We denote it by  $|D|^t$ and write its members as 
$\{D_{\lambda}:\lambda\in\Lambda\}$.
These are the   {\it algebraic  t-pencils.} 

Note that $|D|^t$ is ample iff the members of  $|D|_{\rm alg}$ are ample $\q$-Cartier divisors.
\end{exmp}

The following is proved in  \cite[5.2.8]{k-lo-2}.

\begin{prop} \label{t.penc.alg.deg.char}
Let $X$ be a normal, projective variety of dimension $n$ over an infinite field and $H$ an ample divisor. A  t-pencil $\{D_{\lambda}:\lambda\in\Lambda\}$  is algebraic iff 
there is an infinite subset $\Lambda^*\subset \Lambda$ such that  the intersection numbers $\bigl(H^{n-1}\cdot D_{\lambda}\bigr)$ are independent of  $\lambda\in\Lambda^*$.
\end{prop}

Proof. Let $p:X\map C$ be an algebraic pencil. 
$C$ has infinitely many points of degree $d$ for some $d>0$,
these give infinitely many divisors of the same degree in the corresponding t-pencil. 

The converse is proved in \cite[5.2.8]{k-lo-2}.\qed
\medskip

Combining (\ref{t.penc.alg.deg.char}) with (\ref{secs.isom.zeros.prop.c2}) we get the following.

\begin{thm} \label{t.penc.alg.top.prop.thm}
Let  $X$ be a normal, projective, geometrically integral  variety of dimension $>1+\BH(k)$ over a  field $k$ of characteristic 0.  Then algebraicity of ample t-pencils is a property of $\bigl(|X|, \simsa\bigr)$. \qed
\end{thm}

\noindent{\it Remark \ref{t.penc.alg.top.prop.thm}.1.}  
We use ampleness only to guarantee that almost all members are geometrically connected;
 a little more work should establish the theorem  for non-ample t-pencils as well.
If the pencil has a basepoint, then most likely almost all members are
geometrically connected, hence we can use (\ref{num.eq.top.say}.3)  do decide when they have the same degree.
Basepoint-free t-pencils are always algebraic by (\ref{PTBPS.thm}).
\medskip

The following is a linearity criterion for algebraic t-pencils.

\begin{lem}[Linearity test]\label{lin.test.lem}
Let  $X$ be a normal, projective variety over a perfect  field $k$. 
 Let $|D|_{\rm alg}$ be a  non-composite pencil with parameter curve $C$ and corresponding topological pencil $|D|^t$.  
\begin{enumerate}
\item  If  $|D|^t$ has a smooth basepoint then $g(C)=0$.
\item If $g(C)=0$ and there is a geometrically irreducible subvariety  $W\subset \bs |D|^t$ such that some $D_\lambda\in |D|^t$ is generically smooth along $W$, then
$C\cong \p^1$.
\end{enumerate}
\end{lem}

Proof. The first claim is classical,  see for example \cite[VI.1.9]{rc-book} for a more general assertion. 

  If $g(C)=0$ and $C$ has a $k$-point then
$C\cong \p^1$.  Every member of $|D|^t$ is obtained as
$D_\lambda=\cup_{\sigma}  D_{c,i}^\sigma$  (\ref{t.penc.true.meb.defn}). 
If $D_{c,i}$ contains $W$ then so do its conjugates, so
$D_\lambda$ is singular along $W$ iff $k(c)\neq k$.
 \qed

\medskip

The following is a very simple topological formula computing the intersection number of a curve with a  pencil of divisors.

\begin{lem}\label{int.number.curve.pencil}
  Let  $X$ be a normal, projective variety over a  field $k$ of characteristic 0,  $|D|_{\rm alg}$  a pencil with parameter curve $C$ and corresponding topological pencil $|D|^t$. Let $A\subset X$ be an irreducible curve  disjoint from  $\bs|D|^t$. Then 
 $$
\bigl(A\cdot |D|_{\rm alg}\bigr)=\max\{\# |A\cap D_{\lambda}|:D_\lambda\in  |D|^t\}.
$$
\end{lem}

Proof.  The inequality $\geq $ is clear. 
Conversely, the pencil is given by a map  $X\map C$ that induces a dominant morphism $\pi: A\to C$.  We need to find a $c\in C(\bar k)$ such that
$\pi$ is \'etale over $c$ and $k(a_i)=k(c)$ for every
$a_i\in \pi^{-1}(c)$.  This is possible by (\ref{secs.isom.zeros.prop.v2}). \qed

\medskip{\bf True members}\medskip

\begin{defn}\label{true.memb.defn}
Let $X$ be  a normal, projective, geometrically integral $k$-variety
and  $|D|_{\rm alg}$ a  non-composite algebraic pencil with
corresponding t-pencil  $|D|^t=\{D_{\lambda}:\lambda\in\Lambda\}$.
$D_{\lambda}$ is called a {\it true member} of  $|D|^t$
if there is a $k$-member  $D_c\in |D|_{\rm alg}$ such that $D_{\lambda}=\red D_c=D_c$. 

If $|D|_{\rm alg}$  is a linear pencil with parameter curve
 $C\cong \p^1$,  then $D_c$ is a true member for  all but finitely many  $c\in \p^1(k)$, but in general there may not be any true members.
\end{defn}

\begin{defn} Let $X$ be  a normal, projective, geometrically integral $k$-variety over  an infinite field and $L$ an ample line bundle on $X$. Let  $|D|^t$  be a linear, algebraic  t-pencil.  Define its {\it $L$-degree}  as
 the smallest number in the set $\{\deg_L D_{\lambda}: \lambda\in\Lambda\}$ that is  taken infinitely many times.
We denote it by 
$\deg_L |D|^t$.

Note that $\deg_L |D|^t=\deg_L |D|$.  (For non-linear pencils the definition makes sense but frequently gives a multiple of $\deg_L |D|$.)
\end{defn}

\begin{lem}[True membership test] \label{true.memb.char.lem}
Let $X$ be a normal, projective variety over a perfect field $k$ and 
$|D|^t=\{D_{\lambda}: \lambda\in \Lambda\}$ a linear,  algebraic t-pencil. 
\begin{enumerate}
\item Almost all of the  $\{ D_\mu\colon \deg_H D_{\mu}=\deg_L |D|^t\}$ are true $k$-members of $|D|^t$. 
\item  If such a $D_\mu$ is generically smooth along a geometrically irreducible subvariety  $W\subset \bs |D|^t$,   then
it is a true member of $|D|^t$.
\end{enumerate}
\end{lem}

Proof. The first claim is clear. As in  (\ref{lin.test.lem}),
$D_\lambda=\cup_{\sigma}  D_{c,i}^\sigma$  (\ref{t.penc.true.meb.defn})
is  smooth along $W$ iff $c\in C(k)$. If this holds then
$\deg_LD_\lambda\leq \deg_LD_c=\deg_L|D|$ and equality holds iff
$D_\lambda=D_c$, that is, iff $D_\lambda$ is a true member of $|D|^t$.
 \qed

\medskip

The following example shows that in general there can be false members
that are linearly equivalent to true members.

\begin{exmp} \label{non.true.emp}
Start with the pencil  on $\a^2$
$$
\bigl| (u^2+v^2+u)(v^2+u), (u^2+v^2+v)(u^2+v)\bigr|.
$$
Its general member is a quartic with a node at the origin, but it has 2 members that split into conics.

Next make a change of variables  $u=x+iy, v=x-iy$. 
The resulting pencil  $|D|$ is still defined over $\q$
but now it has  a conjugate pair of reducible members. 
Thus we obtain that
$$
\begin{array}{l}
\bigl((x+iy)^2+(x-iy)^2+(x-iy)\bigr)
\bigl((x+iy)^2+(x-iy)^2+(x+iy)\bigr)=\\
(2x^2+2y^2+x-iy)(2x^2+2y^2+x+iy)=\\
(2x^2+2y^2+x)^2+y^2\qtq{and}\\[1ex]
\bigl((x-iy)^2+(x+iy)\bigr)\bigl((x+iy)^2+(x-iy)\bigr)=\\
\bigl(x^2-y^2+x+i(y-2xy)\bigr)\bigl(x^2-y^2+x+i(y-2xy)\bigr)=\\
(x^2-y^2+x)^2+(y-2xy)^2
\end{array}
$$
both give degree 4 false members of  $|D|^t$.

\end{exmp}

\begin{thm} \label{lin.eq.thm}
Let $X$ be a normal, projective variety over an infinite, perfect   field $k$. 
Two reduced divisors $A_1, A_2$ are linearly equivalent iff the following holds.

There is a codimension 2 subset $\Sigma\subset X$ such that for every
irreducible curve $C\not\subset \Sigma\cup A_1\cup A_2$  there are 
irreducible  divisors $A', H, H'$ and  t-pencils  $|D_1|^t$ and $|D_2|^t$, such that, for $i=1,2$,
\begin{enumerate}
\item $H, H'$ are $\q$-Cartier and ample,
\item $ |D_i|^t$ are algebraic, 
\item $A_i+A'+H$ and $H'$ are reduced members of $ |D_i|^t$,  
\item $(A_i\cup A'\cup H)\setminus (\Sigma\cup H')$ are connected, 
\item  $C\subset \bs |D_i|^t$, and
\item $A_i+A'+H$ and $H'$ are   generically smooth along $C$.
\end{enumerate}
\end{thm}

Proof. If $A_1, A_2$ are linearly equivalent then first choose
$\Sigma=\sing X$ and $A'$ such that $A_i+A'$ are reduced and  $\q$-Cartier.  Then choose 
$H$ sufficiently ample, generically smooth along  $C$ and
such that $|A_1+A'+H|=|A_2+A'+H|$ is very  ample and has an irreducible 
 member $H'$ that is  generically smooth  along  $C$ and such that the $(A_i\cup A'\cup H)\setminus (\Sigma\cup H')$ are connected. Set $|D_i|:=|A_i+A'+H, H'|$. 
Then $A_i+A'+H$ and $H'$ are reduced members of $ |D_i|^t$  and the rest is clear.

Conversely, choose $C$ to be geometrically irreducible such that $X$ is generically smooth along $C$. Then the $|D_i|^t$ are linear pencils by
(\ref{lin.test.lem}). Furthermore,   $A_i+A'+H$ and $H'$ are true members of $ |D_i|^t$ by (\ref{true.memb.char.lem}). Thus  
$A_1+A'+H\sim H'\sim A_2+A'+H$ by (\ref{lin.eq.thm}.3),  hence
$A_1\sim A_2$. \qed


\begin{thm} \label{lineq.top.prop.thm}
Let  $X$ be a normal, projective variety of dimension $>1+\BH(k)$ over a  field $k$ of characteristic 0.  Then linear equivalence of divisors is determined by  $\bigl(|X|, \simsa\bigr)$. 
\end{thm}

Proof. Let $A=\sum a_iA_i$ be an effectve divisor and fix an ample divisor $H$.
For every $i$ there is an $m_i$ such that $|A_i+m_iH|$ is mobile.
Let $A_i^*$ be the sum of $a_i$ distinct, irreducible members of $|A_i+m_iH|$.
Then
$A+(\tsum a_im_i)H\sim \tsum A^*_i$ and $\tsum A^*_i$ is reduced. Thus,
as in \cite[5.3.8]{k-lo-2}, we see that linear equivalence of all  divisors is determined by  linear equivalence of reduced divisors.  

Thus it remains to show that linear equivalence of reduced divisors is determined by  $\bigl(|X|, \simsa\bigr)$. 
For this we need to check that the conditions
(\ref{lin.eq.thm}.1--6) are determined by $\bigl(|X|, \simsa\bigr)$. 

For (\ref{lin.eq.thm}.1) we use (\ref{ample.is.top.lem}) and  for 
(\ref{lin.eq.thm}.2) we need (\ref{t.penc.alg.top.prop.thm}). 
Conditions (\ref{lin.eq.thm}.3--5) are topological.
 Finally for (\ref{lin.eq.thm}.6) we can use 
(\ref{red.int.top.prop.cor.6}). \qed

\section{Complete intersections}\label{compl.int.sec}

In this section we collect various results on complete intersections that were used earlier.

\begin{say}[Complete intersections]\label{ci.sci.say}
 Let $X$ be an irreducible variety.
A subscheme  $Z\subset X$ of codimension $r$ is a {\it complete intersection}
(resp.\ {\it set-theoretic complete intersection}) if there are effective Cartier divisors $D_1,\dots, D_r$ such that
$Z=D_1\cap\cdots\cap D_r$  scheme theoretically  (resp.\ 
$\supp Z=\supp(D_1\cap\cdots\cap D_r)$).

If the $D_i$ are ample, we call $Z$ an  ample (set-theoretic)  complete intersection, usually abbreviated as   {\it ample-ci} resp.\ {\it ample-sci.}

If $H$ is a Cartier divisor and  $D_i\in |m_iH|$   for every $i$, then we say that  $Z$ is a {\it complete $H$-intersection.}
These we usually abbreviate as  {\it $H$-ci, } the set-theoretic versions as {\it $H$-sci.}

Ample complete intersections inherit many properties of a variety, but
the strongest results are for general complete intersections; that is, when the
$D_i\in |m_iH|$ are sufficiently general.
\end{say}

\begin{say}[General  complete intersections] \label{gen.ci.say}
We say that a 
{\it general  complete intersection has property} ${\mathcal P}$ if the following holds.
\begin{enumerate}
\item Let $X$ be a  projective  variety and $H_i$ ample divisors on $X$.
Then for $m_i\gg 1$ there is an open, dense subset
$U\subset \prod_i |m_iH_i|$ such that if $D_i\in |m_iH_i|(\bar k)$ and
$(D_1,\dots, D_r)\in U$  then
$Z:= D_1\cap\cdots\cap D_r$ satisfies ${\mathcal P}$.
\end{enumerate}
There is a long list of such properties ${\mathcal P}$; below is a  list of those that we use. They usually go by the name  {\it Bertini theorems.} Most of these are proved in 
\cite[II.8.18, III.10.9]{hartsh}, but for (\ref{gen.ci.say}.5) one needs to use that a general hypersurface section of an $S_2$ scheme is also $S_2$; see
\cite[IV.12.1.6]{ega}.

\begin{enumerate}\setcounter{enumi}{1}
\item  Given finitely many irreducible subvarieties  $W_j\subset X$, we have $\dim (Z\cap W_j) = \dim W_j-r$  (or the intersection is empty).
\item  Given finitely many irreducible subvarieties  $W_j\subset X$ of dimension $\geq r+1$, the $Z\cap W_j$ are irreducible.
\item Given finitely many locally closed, smooth subvarieties $W_j\subset X$,
the $Z\cap W_j$ are smooth.
\item Given finitely many locally closed, normal subvarieties $W_j\subset X$,
the $Z\cap W_j$ are normal.
\end{enumerate}

In many applications, we need complete intersections that are special in some respects but general in some others.  Here we deal with local conditions.
\begin{enumerate}\setcounter{enumi}{5}
\item A set ${\rm LC}$  of {\it local conditions} consist of  finitely many smooth points  $p_j\in X$ with maximal ideals
$m_j\subset \o_{p_j,X}$, natural numbers $n_j$  and  regular sequences
$g_{j1},\dots,g_{jr} \in \o_{p_j,X}$.  
\end{enumerate}
For example, we can specify whether $p_i\in Z$ or that $Z$ be smooth at $p_i$ with given tangent space.

A codimension $r$ complete intersection  $Z=D_1\cap\cdots\cap D_r$  {\it satisfies} ${\rm LC}$ if  $g_{ji}$ is a  local equation for $D_i$ at $p_j$, modulo  $m_j^{n_j}$.

For every $i$ the local conditions define a linear subspace
$\bigl(|m_iH_i|,{\rm LC}\bigr) \subset |m_iH_i|$  (which may be empty).
The complete intersections that satisfy ${\rm LC}$ form an open, dense subset of
$$
\tprod_i \bigl(|m_iH_i|,{\rm LC}\bigr) \subset \tprod_i |m_iH_i|.
$$
The following combines the local conditions with (\ref{gen.ci.say}.2--5).
\begin{enumerate}\setcounter{enumi}{6}
\item Let $X$ be a  projective  variety,  $H_i$ ample divisors on $X$,
 ${\rm LC}$ a set of local conditions at the points $\{p_j\}$ 
and   ${\mathcal P}$  any of the properties (\ref{gen.ci.say}.2--5). 
Then for $m_i\gg 1$ there is an open, dense subset
$U\subset \prod_i \bigl(|m_iH_i|,{\rm LC}_i\bigr)$ such that if 
$(D_1,\dots, D_r)\in U(\bar k)$  then
\begin{enumerate}
\item[(a)]  $Z:= D_1\cap\cdots\cap D_r$ satisfies ${\mathcal P}$ on $X\setminus\{p_j\}$, and 
\item[(b)]  $\bigl(I_{D_i}, m_j^{n_j}\bigr)= (g_{ji}, m_j^{n_j})$ for every $j$. \qed
\end{enumerate}
\end{enumerate}
\medskip

{\it Finite fields \ref{gen.ci.say}.8.}  If $k$ is finite, then a dense open subset of $\a^n$ may be disjoint from $\a^n(k)$. Nonetheless, the results of 
\cite{MR2385639, cha-poo}  
say that the open sets $U$ above have $k$-points. 

Our proofs have other problems with finite fields, so we will be able to make only very limited use of these cases.
\end{say}

\begin{say}[Connectedness]\label{sci.say}
 Let  $Z$ be a scheme. Connectedness and irreducibility of $Z$ depends only on the topological space $|Z|$, but
geometric connectedness and geometric irreducibility can not be determined using $|Z|$ only. 

We frequently need to guarantee that certain schemes   are 
geometrically connected.  
The next criterion can be proved by repeatedly using \cite[II.7.8]{hartsh}; see also \cite{MR0142547}. 
\medskip

{\it Claim \ref{sci.say}.1.} Let $X$ be a normal, projective, geometrically irreducible variety and  $Z\subset X$  a positive dimensional,  
 ample-sci.  Then $Z$ is
geometrically  connected.\qed
\medskip

 Note that  a proper $k$-scheme $Y$ is geometrically connected iff 
$H^0(Y, \o_Y)$ is a  local, Artin $k$-algebra
such that 
$H^0(Y, \o_Y)\big\slash\sqrt{0}$ is a purely inseparable field extension of $k$.  We can thus restate (\ref{sci.say}.1) as follows.
\medskip

{\it Claim \ref{sci.say}.2.} Let $X$ be a normal, projective, geometrically irreducible variety and  $Z\subset X$  a positive dimensional, reduced, 
 ample-sci.  Then $k[Z]/k[X]$ is
purely inseparable.\qed
\end{say}

We also use the following  variant of the Lefschetz hyperplane theorem.
For the Picard variety, this is proved in \cite{sga2}. 
For normal varieties in characteristic 0, the class group version  is proved in \cite{MR2219849, MR2567426}. For positive characteristic see \cite{ji-tocome}. 

\begin{thm}\label{gnl-sri.thm} 
Let $X$ be a geometrically normal, projective variety
and $|H|$ an ample  linear system  on $X$. Let $Z\subset X$ be a normal, $\geq 2$ dimensional,  complete $H$-intersection. Then 
the restriction map 
$$
\cl(X)\into \cl(Z)\qtq{is injective.}
$$ 
\end{thm}

\begin{say}[Disjointness of conjugates]\label{disj.conj.say}

Let $K/k$ be a finite, separable  field extension and 
$A\in \operatorname{Mat}_{r\times n}(K)$ an $r\times n$ matrix.
Then  $A{\mathbf x}=0$ defines a  linear subspace
$L_A\subset K^n$, usually of codimension $r$.
Let  $\sigma:K\into \bar k$ be a $k$-embedding, it defines
$A^\sigma$ and  $L_A^\sigma\subset K^n$.  We aim to show that usually
$L_A\cap L_A^\sigma$ has codimension $2r$.

  After choosing a $k$-basis $e_i$ for $K$, we can write
$A=\sum_i e_i A_i$ where  $A_i\in \operatorname{Mat}_{r\times n}(k)$.

\medskip
{\it Claim \ref{disj.conj.say}.1.} There is a Zariski open subset $U\subset \oplus_i \operatorname{Mat}_{r\times n}(k)$
such that if $A\in U$ then $\codim (L_A\cap L_A^\sigma)=\min\{2r,n\}$ for all non-identity $k$-embeddings $\sigma:K\into \bar k$.
\medskip

Proof. Let us work out the critical case when $n=2r$. 
We use  $x_1,\dots, x_r, y_1,\dots, y_r$ as variables.
Write  the equations as  $A{\mathbf x}+B{\mathbf y}=0$ where
$A, B\in \operatorname{Mat}_{r\times r}(K)$. 
The conjugate equations  are  $A^\sigma{\mathbf x}+B^\sigma{\mathbf y}=0$. We need to show that  
$$
\rank\left(
\begin{array}{cc}
A & B \\
A^\sigma & B^\sigma
\end{array}
\right)=2r
$$
is a non-empty, Zariski  open condition. It  is clearly an open condition, so need to find one example with maximal rank.
Take  $A, C\in \operatorname{Mat}_{r\times r}(k)$ and $B=\alpha C$. By row reduction
$$
\left(
\begin{array}{cc}
A & B \\
A^\sigma & B^\sigma
\end{array}
\right)
=
\left(
\begin{array}{cc}
A &  \alpha C\\
A & \alpha^\sigma C
\end{array}
\right)
\rightsquigarrow
\left(
\begin{array}{cc}
A &  \alpha C\\
0 & (\alpha^\sigma-\alpha)C
\end{array}
\right)
$$
 whose determinant is  $(\alpha^\sigma-\alpha)^r\det A\det C$. \qed

\medskip

The intrinsic way of passing from  $\operatorname{Mat}_{r\times n}(K)$ to $\oplus_i \operatorname{Mat}_{r\times n}(k)$ is the {\it Weil restriction,} denoted by $\res^K_k(\ )$; see, for instance  \cite[Sec.7.6]{blr}. 
We can thus globalize
(\ref{disj.conj.say}.1) first to projective spaces and then to 
their subvarieties as follows.

\medskip
{\it Claim \ref{disj.conj.say}.2.} Let $X$ be a $k$-variety of pure dimension $n$
and $|M_1|, \dots, |M_r|$ basepoint-free linear systems.  Let $K/k$ be  a finite, separable field extension.  Then  there is a dense, Zariski open subset 
$$
U\subset \res^K_k |M_1|\times\cdots\times   \res^K_k |M_r|,
$$
such that, if $(D_1, \dots, D_r)\in U$, then 
$$
\codim_X \bigl(D_1\cap \cdots\cap D_r\cap D_1^\sigma\cap \cdots\cap D_r^\sigma\bigr)=\min\{2r,n+1\},
$$
for all non-identity $k$-embeddings $\sigma:K\into \bar k$. \qed
\medskip
\end{say}

\begin{lem}\label{ZW.H.link.prop.1p.cor.3.lem}
Let $k$ be an infinite  field  and 
 $X$  a  normal, projective, irreducible  $k$-variety of dimension $n>2r$. 
Let  $p,q\in X$  be  closed points
such that there are embeddings  $k\subset k(p)\subset k(q)^{\rm ins}\subset \bar k$. 

Then there is an irreducible, $r$-dimensional  $k$-variety  $W\subset X$ such that 
\begin{enumerate}
\item $p, q\in W$ and 
\item  $k(p)/k[W]$ is purely inseparable.
\end{enumerate}
Furthermore, if $p$ is a smooth, separable  point of $X$ then we can also assume that 
\begin{enumerate}\setcounter{enumi}{2}
\item $p$ is a smooth point of $W$.
\end{enumerate}
\end{lem}

Proof. Let $k\subset K_p\subset k(p)$ and
$k\subset K_q\subset k(q)$ be maximal separable subextensions.
After base change to $K_p$, we have a degree 1 point $\bar p$ lying over $p$ and a degree $=\deg (K_q/K_p)$  point $\bar q$ lying over $q$. 
Let $\{\sigma\}$ be the set of all $k$-embeddings
$\sigma:K_p\into \bar k$. 
Thus $\bar p$ and $\bar q$ each have  $\deg (K_p:k)$ conjugates over $k$ and these are disjoint from each other.

Next  take a general ample-ci variety   (\ref{ci.sci.say})  $W_1\subset X_{K_p}$ that contains
$\bar p$ and $\bar q$.  By (\ref{disj.conj.say}.2) 
the $W_1^\sigma$ are disjoint from each other. Thus their union
$W_{K_p}=\cup_{\sigma}W_1^\sigma$  descends to a $k$-subvariety $W\subset X$
with the required properties.\qed

\section{Picard group, class group and Albanese variety}\label{pic.cl.alb.sec}

For the  Picard group and Picard scheme, \cite[Lects.V-VI]{FGA}, \cite[Sec.19]{mumf66}, 
\cite{mumf-abvar}
 or \cite{blr} contain proofs and details; for these we just fix our notation. 
Modern references for the class group and Albanese variety are harder to find; about these we give more details.

\begin{say}[Picard group  of a normal variety]\label{pic.cl.norm.say}
The group of line bundles on a scheme  $X$ is the {\it Picard group} of $X$, denoted by $\pic(X)$.  If $X$ is proper then $\pic^\circ(X)\subset \pic(X)$ denotes  the subgroup of divisors that are algebraically equivalent to 0. 
The quotient   $\ns(X):=\pic(X)/\pic^\circ(X)$ is   the {\it N\'eron-Severi} group of $X$.
It is a finitely generated abelian group.

If $X$ is proper over a field $k$
with algebraic closure $\bar k$
then $\pic(X_{\bar k})$ 
 has a natural $k$-scheme structure, denoted by $\pics(X)$. The identity component  is denoted by $\pico(X)$, it is a commutative algebraic group. 
If $X$ is geometrically normal and $\chr k=0$ then $\pico(X)$ is an Abelian variety.
If $X$ is geometrically normal and $k$ is perfect then $\red\pico(X)$ is an Abelian variety.
The  non-reduced structure of  $\pico(X)$ will pay no role in our questions.

There is a natural inclusion $\pic(X)\into \pics(X)(k)$ which is an
isomorphism if $X$ has a $k$-point. In general the quotient
$\pics(X)(k)/\pic(X)$ is a torsion group.

\medskip

{\it Claim \ref{pic.cl.norm.say}.1.} For Abelian varieties, $A\mapsto \pico(A)$ is a duality. 
\medskip

Proof. Over $\c$, this is a special case of the  Appell-Humbert theorem
\cite[pp.21-22]{mumf-abvar}. In general see \cite[Sec.13]{mumf-abvar}. \qed
\end{say}

\begin{say}[Class group of a normal variety]\label{cl.norm.say}
For  a  normal and proper $k$-variety $X$,  let $\cl(X)$ denote the group of Weil divisors modulo linear equivalence. It is also  isomorphic to  the group of  rank 1 reflexive sheaves, where the product is the double dual of the tensor product.

Let  $\clo(X)\subset \cl(X)$  be the subgroup  of divisors that are algebraically equivalent to 0.
I call  the quotient  $\ns^{\rm cl}(X):=\cl(X)/\clo(X)$ the {\it N\'eron-Severi class group}\footnote{The literature seems inconsistent. Frequently this is called the N\'eron-Severi  group.} of $X$, and  its  $\q$-rank  the
{\it class rank} of $X$,  denoted  by  $\rho^{\rm cl}(X)$.

Note the we have natural inclusions  
$$
\pic^{\circ}(X)\subset \clo(X)\qtq{and} \ns(X)\subset \ns^{\rm cl}(X),
\eqno{(\ref{cl.norm.say}.1)}
$$
that are 
isomorphisms iff every Weil divisor is Cartier, for example
when $X$ is smooth.

Basic results about these groups are the following.
\medskip

{\it Claim \ref{cl.norm.say}.2.} Let $p:Y\to X$ be a birational morphism of  normal, proper varieties over a perfect field $k$. Then
\begin{enumerate}
\item[(a)] $p_* :\clo(Y)\to \clo(X)$ is an isomorphism and
\item[(b)] $p_*: \ns^{\rm cl}(Y)\onto  \ns^{\rm cl}(X)$ is onto.
\end{enumerate}
\medskip

\medskip

{\it Claim \ref{cl.norm.say}.3.} Let $X$ be a normal, proper variety over a perfect field $k$. Then there is normal, proper variety $Y$ and a birational morphism $p:Y\to X$ such that  $\clo(Y_K)=\pic^\circ(Y_K)$ for every $K\supset k$. 
\medskip

It is quickest to prove these by  using
the Albanese variety, see  (\ref{pic.alb.say}.5--6).  
As a consequence, we can define the scheme structure of $\clo$ by
$$
\clos(X)\cong \clos(Y)\cong \red\pico(Y).
\eqno{(\ref{cl.norm.say}.4)}
$$
 In the complex case  these results  go back to  Picard \cite{picard-95} and Severi \cite{MR1511373}, but 
the most complete references may be the  papers of
Matsusaka   \cite{mats-pic} and  of N\'eron \cite{neron-pic};
see also \cite[Sec.3]{k-mumf} for some discussions.

 More recent results on various aspects of the class group of singular varieties are discussed in \cite{MR1203911, MR2457299, MR2567426}.

\end{say}

\begin{defn} \label{pic.Z.X.def}
Let $X$ be a normal, proper $k$-variety and $\Sigma\subset X$ a  subset. Let
$\wdiv(X, \Sigma)\subset \wdiv(X)$  and
 $\cl(X, \Sigma)\subset \cl(X)$ denote the subgroup of those Weil divisors that are Cartier at every point $x\in \Sigma$. (The concept and notation are not standard.)

Note that  $\cl(X, \Sigma)$ is isomorphic to  the group of those rank 1 reflexive sheaves that are locally free  at every point $x\in \Sigma$.

We see in (\ref{pic.Z.X.lem}) that 
$\cl(X_{\bar k}, \Sigma_{\bar k})$ is naturally identified with  a closed $k$-subgroup $\cls(X, \Sigma)\subset \cls(X)$.
We denote its identity  component by  $\clos(X, \Sigma)$.  
Note that in general  $\cls(X, \Sigma)\cap \clos(X)$ may be disconnected.

The quotient
$\ns^{\rm cl}(X, \Sigma):=\cls(X, \Sigma)/\clos(X, \Sigma)$  is finitely generated.
\end{defn}

\begin{lem} \label{pic.Z.X.lem}
Let $X$ be a  normal, proper variety over a perfect field $k$ and $\Sigma\subset X$ an arbitrary  subset. Then there is a closed, algebraic  $k$-subgroup
$\cls(X, \Sigma)\subset \cls(X)$ such that $\cl(X_{\bar k}, \Sigma_{\bar k})=\cls(X, \Sigma)(\bar k)$.
\end{lem}

Proof.  Assume first that  $\Sigma=\{x\}$   is a closed point  and
there is a universal family  ${\mathbf L}$ on $X\times \clos(X)$ that is flat  over  $\clos(X)$. The set of points  $V\subset X\times \clos(X)$ where ${\mathbf L}$ is not locally free is closed. Since 
$\cls(X, \{x\})$ is the complement of the image of
$V\cap (\{x\}\times \clos(X))$, it is constructible.  It is also a subgroup and a constructible subgroup is closed. 

In general such an ${\mathbf L}$ does not exist, but we check in (\ref{pic.Z.X.lem.0}) that a flat universal family exists after a finite field extension and a constructible subdivision  $\tau:\amalg_j W_j\to \clos(X)$.
The argument above then shows that $\cls(X, \{x\})$ is constructible, hence closed as before.

If $\Sigma$   is any set of closed  points then $\cls(X, \Sigma)=\cap_{x\in \Sigma} \cls(X, \{x\})$.

If $\eta$   is a non-closed  point,  then
$\cls(X, \{\eta\})$ is the union of all
$\cls(X, \Sigma_U)$, where $U$ runs through all open subsets of $\bar\eta$ and
$\Sigma_U$ denotes the set of closed points of $U$. By the Noetherian property,
$\cls(X, \{\eta\})=\cls(X, \Sigma_U)$ for some $U$. 

Finally,
$\cls(X, \Sigma)=\cap_{x\in \Sigma} \cls(X, \{x\})$ holds for any set of points $\Sigma$. \qed


\begin{lem} \label{pic.Z.X.lem.0}
Let $X$ be a normal, proper variety over an algebraically closed field $K$.
There is a locally closed decomposition  $\tau:\amalg_j W_j\to \clos(X)$ such for every $j$  there is a universal family  ${\mathbf L}_j$ on $X\times W_j$ that is flat  over  $W_j$ and whose fiber over $w\in W_j$ is the reflexive sheaf corresponding to
$\tau(w)\in \clos(X)$.  
\end{lem}

Proof. By (\ref{cl.norm.say}.3), there is a proper, birational morphism from a normal variety
$p:Y\to X$ such that $\clos(X)=\pico(Y)$. Let
${\mathbf L}$ be the universal line bundle over $Y\times \pico(Y)$.
Pushing it forward we get a rank 1  sheaf
$$
{\mathbf L}_X:=\bigl(\pi_*{\mathbf L}\bigr)^{[**]}\qtq{over} X\times \clos(X).
$$
In general $ {\mathbf L}_X$ is not flat over $\clos(X)$. However, by generic flatness, $ {\mathbf L}_X$ is  flat with reflexive fibers over a dense, open subset   $W_1\subset \clos(X)$. Repeating this with $\clos(X)\setminus W_1$ we get the required
locally closed decomposition. \qed

\medskip{\bf Albanese variety}

\begin{say}[Albanese variety]\label{pic.alb.say}
Let $X$ be a proper, normal variety  over a perfect field $k$.
There are 2 different notions of the {\it Albanese variety} of $X$ in the literature.
In \cite[VI.3.3]{FGA} it is the target of the universal {\em morphism} from $X$ to an Abelian torsor; that is, a principal homogeneous space under an Abelian variety. We denote this by
$$
\albm^{\rm gr}_X: X\to \alb^{\rm gr}(X).
\eqno{(\ref{pic.alb.say}.1)}
$$
Pull-back by $\albm^{\rm gr}_X$ gives an isomorphism
$$
\pico\bigl(\alb^{\rm gr}(X)\bigr)\cong \red\pico(X).
\eqno{(\ref{pic.alb.say}.2)}
$$
If $X$ has a  $k$-point then  $\alb^{\rm gr}(X)$  an Abelian variety.
$\alb^{\rm gr}(X)$ is a birational invariant for smooth, proper varieties, but not a birational invariant for  normal varieties.

In the older literature, for example  \cite{mats-pic, serre-albanese}, 
the Albanese map  is the universal {\em rational map} from $X$ to an Abelian torsor, called the {\it classical Albanese variety}
$$
\albm_X: X\map \alb(X).
\eqno{(\ref{pic.alb.say}.3)}
$$
If $X$ has a smooth $k$-point then so does $\alb(X)$ and then
it is an Abelian variety.

$\alb(X)$ is a birational invariant of $X$ (for  normal, proper varieties) and 
the two versions coincide if $X$ is smooth. Therefore, 
if $X'\to X$ is a resolution then  $\alb(X)=\alb(X')=\alb^{\rm gr}(X')$. In any case, we get a morphism over the smooth locus
$$
\albm_X: X^{\rm sm}\to \alb(X).
\eqno{(\ref{pic.alb.say}.4)}
$$
Let $X'$ be the normalization of the closure of the graph of $\albm_X$.
Then we have a commutative diagram
$$
\begin{array}{lcl}
& X' &\\
\  p\swarrow && \searrow \albm_{X'}   \\
X & \stackrel{\albm_X}{\dasharrow} & \alb(X)=\alb(X'),
\end{array}
\eqno{(\ref{pic.alb.say}.5)}
$$
where $\albm_{X'}$ is a morphism. In particular (\ref{pic.alb.say}.2) gives that
$$
\clos(X')=\red\pico(X')\cong \pico\bigl(\alb(X')\bigr).
\eqno{(\ref{pic.alb.say}.6)}
$$
Therefore
$$
\clos(X)\cong \pico\bigl(\alb(X)\bigr).
\eqno{(\ref{pic.alb.say}.7)}
$$
Let $p:Y\map X$ be a map of normal varieties. As long as
$p(Y)$ is not contained in the singular locus of $X$, the composite
$\alb_X\circ\ p:Y\map \alb(X)$ is defined, hence we get a morphism
$$
\albm_p:\alb(Y)\to \alb(X).
\eqno{(\ref{pic.alb.say}.8)}
$$
\end{say}

\begin{say} \label{pic.alb.rel.say} Let $p:Y\to X$ be a morphism of  normal varieties.
Let 
$$
\albm^{\rm gr}_{Y\backslash X}: X\to \alb^{\rm gr}(Y\backslash X)
\eqno{(\ref{pic.alb.rel.say}.1)}
$$
denote the universal  morphism from $X$ to an Abelian torsor that maps every irreducible component of $Y$ to a point. Thus we get
$$
 \alb^{\rm gr}(Y)\to \alb^{\rm gr}(X)\to \alb^{\rm gr}(Y\backslash X).
\eqno{(\ref{pic.alb.rel.say}.2)}
$$
We claim that the dual sequence
$$
0\to \pico\bigl(\alb^{\rm gr}(Y\backslash X)\bigr)\to \pico(X)\to \pico(Y)
\eqno{(\ref{pic.alb.rel.say}.3)}
$$
is exact.
To see this, let $K$ denote the kernel of $ \pico(X)\to \pico(Y)$.
It is clear that $\pico\bigl(\alb^{\rm gr}(Y\backslash X)\bigr)\subset K$.
To see the converse, we may assume that $X$ has a $k$-point. 
By (\ref{pic.cl.norm.say}.1) we get an exact sequence
$$
\alb^{\rm gr}(Y)\to \alb^{\rm gr}(X)\to \pico(K)\to 0.
$$
The resulting $Y\to \pico(K)$ maps every irreducible component of $Y$ to a point, so it factors through
$\alb^{\rm gr}(Y\backslash X)$. By duality  we get 
$K\to \pico\bigl(\alb^{\rm gr}(Y\backslash X)\bigr)$. \qed
\end{say}

We would like to know when $\albm_p$ is dominant. Lefschetz theory suggests that this should hold if $p(Y)$ is   ample-ci (\ref{ci.sci.say}).
We will need two generalizations of this.

\begin{lem} \label{subvar.albm.injective.lem}
Let $X,Y$ be normal,   projective varieties  and
$p:Y\to X$ a morphism. Assume  that $p(Y)$ has nonempty intersection with every nonzero divisor
in $X$ and $\albm_X$ is a morphism along $p(Y)$. Then
$\albm_p:\alb(Y)\to \alb(X)$ is surjective.
\end{lem}

Proof. 
If $\albm_p$ is not surjective then the quotient 
$\alb(X)/\albm_p(\alb(Y))$ is positive dimensional. Hence 
there is a nonzero, effective  divisor 
 $D\subset \alb(X)/\albm_p(\alb(Y))$  whose pull-back
to $\alb(X)$ is disjoint from $\albm_p(\alb(Y))$. Then its pull-back to $X$ is a
divisor which is disjoint from $p(Y)$. \qed

\begin{cor}\label{ci.curve.onto.alb.cor}
 Let $X$ be a normal, projective variety  and
$C\subset X^{\rm ns}$ an irreducible  ample-sci curve (\ref{ci.sci.say}). 
Then $\clos(X)\to \jacs(\bar C)$ is injective modulo torsion. \qed
\end{cor}

\begin{lem} \label{subvar.albm.injective.cor}
Let $X$ be a normal,   proper variety. Then there is a finite subset $\Sigma\subset X$ such that the following holds.

Let $Y\subset X$ be an irreducible divisor that is disjoint from $\Sigma$. Assume  that $Y$ has nonempty intersection with every nonzero divisor
in $X$. Then
$\alb(\bar Y)\to \alb(X)$ is surjective.
\end{lem}

Proof.  Consider the normalization of the closure of the graph  of $\albm_X$
$$
X\stackrel{\pi}{\longleftarrow}  X' \stackrel{\albm_X}{\longrightarrow} \alb(X).
$$
Let $E'_i\subset X'$ be the $\pi$-exceptional divisors. Choose $\Sigma$ 
to contain the generic point of  each $\pi(E'_i)$  and 
every  non-Cartier center (\ref{non.car.cent.defn}). 

Then $Y$ is a Cartier divisor and $\pi^{-1}(Y)=\pi^{-1}_*(Y)$. 
Therefore, if $D'\subset X'$ is a divisor then
$$
\pi\bigl(\pi^{-1}_*(Y)\cap D'\bigr)=\pi\bigl(\pi^{-1}(Y)\cap D'\bigr)=Y\cap \pi(D')\neq\emptyset.
$$
Let $\bar Y'$ denote the normalization of $ \pi^{-1}_*(Y)$.
Then $\alb(\bar Y')\to \alb(X)$ is surjective by (\ref{subvar.albm.injective.lem}) and
$\alb(\bar Y')\cong \alb(\bar Y)$. \qed

\begin{defn}[Partial Albanese variety]\label{pic.alb.say.Z}
Let $X$ be a proper, normal variety  over a perfect field $k$
and $\Sigma\subset X$ a  subset.

Define  the {\it Albanese map} of $(X,\Sigma)$   as the universal  rational map  from $X$ to an Abelian torsor,  that is a morphism along $\Sigma$
$$
\albm_{X,\Sigma}: X\map \alb(X,\Sigma).
\eqno{(\ref{pic.alb.say.Z}.1)}
$$
If $\Sigma\subset X^{\rm sm}$ then $\alb(X,\Sigma)=\alb(X)$. In general
$\alb(X,\Sigma)$ is a quotient of $\alb(X)$.
\end{defn}

\begin{thm}\label{pic.alb.Z.lem} Let $X$ be a normal, proper variety over a perfect field $K$ and $\Sigma\subset X$ a  subset. Then 
pull-back by $\albm_{X,\Sigma}$ gives an isomorphism
$$
\albm_{X,\Sigma}^*: \pico\bigl(\alb(X,\Sigma)\bigr) \cong \clos(X,\Sigma).
\eqno{(\ref{pic.alb.Z.lem}.1)}
$$
\end{thm}

Proof. Consider $\albm_{X,\Sigma}: X\map \alb(X,\Sigma)$. By assumption it is a morphism along $\Sigma$, thus the pull-back of a line bundle on $\alb(X,\Sigma)$ is locally free along $\Sigma$. That is, 
$$ \albm_{X,\Sigma}^*\pico\bigl(\alb(X,\Sigma)\bigr)\subset \clos(X,\Sigma).
\eqno{(\ref{pic.alb.Z.lem}.2)}
$$
For the converse, assume first that $\Sigma=\{x\}$ is a closed point. 
As in (\ref{pic.alb.say}.5) we have
$$
\begin{array}{lcl}
& X' &\\
\  p\swarrow && \searrow \albm_{X'}   \\
X & \stackrel{\albm_X}{\dasharrow} & \alb(X)=\alb(X'),
\end{array}
$$
where $\albm_{X'}$ is a morphism. Let $Y'$ be the normalization of
$p^{-1}(x)$. 
Since $X'\to  \alb^{\rm gr}(Y'\backslash X')$ contracts  every irreducible component of $Y'$ to a point, the composite  $X\map X'\to \alb^{\rm gr}(Y'\backslash X')$ is a morphism at $x$ by Zariski's main theorem. This gives $\alb(X,\{x\})\to  \alb^{\rm gr}(Y'\backslash X') $.  Thus we get a commutative diagram
$$
\begin{array}{ccc}
X' & \to & \alb^{\rm gr}(Y'\backslash X')\\
\downarrow && \uparrow \\
X & \to & \alb(X,\{x\}).
\end{array}
\eqno{(\ref{pic.alb.Z.lem}.3)}
$$
If $L\in \clos(X,\{x\})(\bar k)$ then  its pull-back to
$X'$ is trivial along $Y'$, hence it is obtained as the pull-back of a line bundle on $\alb^{\rm gr}(Y'\backslash X')$. Factoring through $\alb(X,\{x\})$
shows that 
$ \clos(X,\{x\})\subset \albm_{X,\{x\}}^*\pico\bigl(\alb(X,\{x\})\bigr)$.

The same argument works for any finite number of closed points.
If $\Sigma$ is an infinite set of closed point then, by the Noetherian property,  $\clos(X,\Sigma)=\clos(X,\Sigma')$ for every large enough finite subset $\Sigma'\subset \Sigma$.

Finally assume that $y$ is a non-closed point. 
Then $\clos(X,\{y\})$ is the union of all $\clos(X,\Sigma_U)$ where
$\Sigma_U$ us the set of all closed points in some open subset $U\subset \overline{\{y\}}$. By the Noetherian property, we have equality $\clos(X,\{y\})=\clos(X,\Sigma_U)$ for some fixed $U$. 
\qed

\begin{cor} \label{pico.bir.mod.cor}
Let $X$ be a normal, projective variety over a perfect field and $Z\subset X$ a closed, reduced subscheme with generic points $g_Z$. 
Then there is a  normal, projective variety $X'$, a birational morphism
$p:X'\to X$ and a  closed, reduced subscheme $Z'\subset X'$ 
 with generic points $g_{Z'}$
such that
\begin{enumerate}
\item $p$ is a local isomorphism at all generic points of $Z'$,
\item  $Z=p(Z')$,
\item  $\albm_{X',g_{Z'}}$ is a morphism along $Z'$ and
\item   $\clos(X,g_Z)=\clos(X',g_{Z'})=\clos(X',Z')$.
\end{enumerate}
If either $\dim Z=1$ or the characteristic is 0, we can also achieve that
\begin{enumerate}\setcounter{enumi}{4}
\item $Z'$ is smooth, 
\end{enumerate}
\end{cor}

Proof.  We can take $X'$ to be the
normalization of the closure of the graph of $\albm_{X,g_Z}$. 
Then we can  resolve the singularities of $Z'$  if desired. 
\qed

\medskip{\bf Non-Cartier centers}\medskip

\begin{defn}[Non-Cartier centers]\label{non.car.cent.defn} 
Let $X$ be a reduced scheme and $D$ an effective  Weil divisor.  There is a unique largest open subscheme $X^{\rm car}_D\subset X$, called the {\it Cartier locus} of $D$   such that the restriction of $D$ to  
$X^0_D$ is Cartier. The complement $X\setminus X^{\rm car}_D$ is 
the {\it non-Cartier locus} of $D$. A point
$x\in X$ is a {\it non-Cartier center} of $X$ if 
there is a  Weil divisor $D$ such that $x$ is the generic point of an irreducible component of the 
non-Cartier locus of $D$.

For example, let $X=(xy=0)\subset  \a^3_{xyz}$ and set
$D_c:=(x=z-c=0)$. Its non-Cartier locus is the point  $(x=y=z-c=0)$. 
Thus every closed point of the $z$-axis is a non-Cartier center of $X$.
The generic point of the $z$-axis is also a  non-Cartier center of $X$
for the divisor $(x=y=0)$. 

In direct analogy   one can define the notions of   {\it $\q$-Cartier locus} and  {\it non--$\q$-Cartier center.}
\end{defn}

The next result of \cite[6.7]{2011arXiv1104.1861B} shows that the situation is quite different for normal varieties.  (Note that \cite{2011arXiv1104.1861B}
works over an algebraically closed field, but this is not necessary.)

\begin{thm} \label{non.car.cent.thm}
A geometrically normal variety  has only finitely many non-Cartier or
non--$\q$-Cartier centers.
\end{thm}

Proof. We may assume that $X$ is proper and irreducible. 
Let $U\subset X$ be an open subset such that 
$X$ has only finitely many non--($\q$-)Cartier or  centers in  $U$. We show that there is a strictly larger  open subset 
$U\subsetneq U'\subset X$ such that 
$X$ has only finitely many non--($\q$-)Cartier centers in  $U'$.
We can start with the smooth locus $U=X^{\rm sm}$,  since it is disjoint from every non--($\q$-)Cartier center. Noetherian induction then gives that
$X$ has only finitely many non--($\q$-)Cartier centers.

Let $Z\subset X\setminus U$ be an irreducible component. 
By (\ref{non.car.cent.lem}) there is a  dense, open subset $Z^0\subset Z$ such that if a  Weil divisor $D$ is ($\q$-)Cartier at the generic point $g_Z\in Z$ then it is ($\q$-)Cartier along $Z^0$.   We may assume that $Z^0$ is disjoint from
every other irreducible component of $X\setminus U$. Then $U':=U\cup Z^0$ is open in $X$ and  $g_Z$ is the only possible new  non--($\q$-)Cartier center in $U'$.  \qed

\begin{lem} \label{non.car.cent.lem}
Let $X$ be a normal, proper variety over an algebraically closed field and $Z\subset X$ an irreducible subvariety. Then there is a dense, open subset $Z^0\subset Z$ such that
the following holds.

Let $D$ be a Weil divisor that is ($\q$-)Cartier at the generic point  $g_Z\in Z$. Then it is ($\q$-)Cartier everywhere along $Z^0$.
\end{lem}

Proof. As in (\ref{pic.Z.X.def}), let $\cls(X,g_Z)\subset\cls(X)$ be the subgroup of those divisors that are Cartier at the generic point of $Z$ and 
$\clos(X,g_Z)\subset\clos(X)$  the identity  component. 

As we noted in (\ref{pic.Z.X.def}), 
the quotient $\cls(X,g_Z)/\clos(X,g_Z)$  is finitely generated; say by the divisors $D_i$. There is a dense, open subset $Z^0_1\subset Z$ such that every $D_i$ is Cartier along $Z^0_1$, hence the same holds for every linear combination of the $D_i$. 

Next we show that there is a dense, open subset $Z^0_2\subset Z$ such that every
divisor in $\clos(X,g_Z)$ is Cartier along $Z^0_2$.
Consider the Albanese map
$\albm_{X,Z}:X\map  \alb(X,Z)$. By (\ref{pic.alb.say.Z}) it is defined at $g_Z$, hence
on a dense, open subset $Z^0_2\subset Z$. By (\ref{pic.alb.Z.lem}), 
$\clos(X,g_Z)$ is the pull-back of $\pico\bigl(\alb(X,g_Z)\bigr)$,
hence every member of $\clos(X,g_Z)$ is locally free along $Z^0_2$.
Finally  $Z^0=Z^0_1\cap Z^0_2$ is the dense, open subset that we need.

If $D$ is $\q$-Cartier at $g_Z$ then $mD$ is Cartier at $g_Z$ for some
$m>0$, hence $D$ is $\q$-Cartier along $Z^0$ by the previous results. 
\qed

\begin{defn}[Relative factorial locus] \label{non.car.cent.thm.cor} Let $X$ be a normal variety  over a perfect field 
and $Z\subset X$ a closed subset. Let $\{w_i\in X: i\in I\}$ be those non-Cartier centers of $X$ whose closure  does not contain any irreducible component of $Z$.
Then 
$$
\operatorname{Fact}(Z\subset X):=Z\setminus \cup_i \bar w_i
\eqno{(\ref{non.car.cent.thm.cor}.1)}
$$
is the unique largest dense open subset  $Z^\circ\subset Z$ such that if a Weil divisor $D$ is Cartier at some point of each irreducible component of $Z$ then it is
Cartier  everywhere along $Z^\circ$. We call $\operatorname{Fact}(Z\subset X)$ the
{\it relative factorial locus} of $X$ along $Z$.
The factorial locus of $X$ along $X$ is the largest open subscheme $X^\circ\subset X$ whose local rings are factorial. 
\end{defn}



\section{Commutative algebraic groups}\label{comm.ag.sec}

\begin{say}[Structure of commutative algebraic groups]\label{com.al.gr.say}
Let  $A$ denote  a  commutative algebraic group over a perfect field $k$,
  $A^\circ\subset A$ the identity component   and
$A^{\rm lin}\subset A$ the largest connected linear algebraic subgroup. Then 
$A/A^\circ$ is a finite abelian group and  $A^\circ/A^{\rm lin}$ is an Abelian variety.

Let  $A^{\rm unip}\subset A^{\rm lin}$ be the largest unipotent subgroup
 and    $A^{\rm tor}\subset A^{\rm lin}$ the largest semisimple subgroup.
Then 
$A^{\rm tor}$ is a {\it torus} (that is, isomorphic to   $\gm^r$ over $ k^{\rm sep}$ for some $r$)   and
$A^{\rm lin}=A^{\rm unip}+A^{\rm tor}$.

$A$ is called  {\it semi-Abelian} if $A^{\rm unip}=0$.

Let $A^{\rm prop}\subset A$ denote the largest proper, connected subgroup.
Then $A^{\rm prop}\cap A^{\rm lin}$ is finite but usually
$A^{\rm prop}+ A^{\rm lin}$ does not equal $A$. 

See    \cite{borel}, \cite[Chap.8]{MR3729270} or \cite{MR3645068}
for details and proofs.
\end{say}

\begin{say}[$\q$-rank]\label{Q.rank.say} For Abelian varieties, the $\q$-rank of $A(k)$ is a subtle invariant of $A$ and $k$; see for example
(\ref{m.w.f.defn}) and (\ref{anti-mordell.defn}). By contrast the $\q$-rank of a linear algebraic group is easy to compute.
\begin{enumerate}
\item   $G(k)$ is torsion for every algebraic group $G$ over a locally finite field $k$.
\item  $U(k)$ is $p^{\infty}$-torsion for every  unipotent algebraic group $U$ over a  field $k$ of characteristic $p>0$.
\item  $\rank_{\q}U(k)=\dim U\cdot \deg(k/\q)$ for every unipotent algebraic group $U$ over a  field $k$ of characteristic $0$.
\item  (Larsen) $\rank_{\q}T(k)=\infty$ for every  positive dimensional  torus over a field that is not locally finite.
\end{enumerate}
Of these only the last claim is nontrivial. More precise versions follow from   weak approximation, but here is a 
shorter argument. 

$T$ is defined over a finitely generated subfield, so we may as well assume that $k$ is either a number field or the field of functions  of a geometrically integral 
curve $C$ over a subfield $k_0\subset k$. 
  Over a dense, open, regular  subset
$U\subset C$ we have a torus $T_U\to U$. 

Assume to the contrary that $t_1,\dots, t_s\in T(k)$ generate a maximal rank subgroup.  We can view the $t_i$ as rational sections of $T_U\to U$. After further shrinking $U$ we may assume that they are all regular sections. 
If $t\in T(k)$ then $t^n\in \langle t_1,\dots, t_s\rangle$ for some $n>0$. Thus $t$ is a rational section that is also finite over $U$, hence a regular section. 
Thus every rational section of $T_U\to U$ is regular. 

Next we show that this is not the case. 
$T$ is isomorphic to   $\gm^r$ over $ k^{\rm sep}$, hence there is a 
finite, separable field extension  $K/k$   such that
$T_K\cong  \gm^r$. (Such a $K$ is called a {\it  splitting field} of $T$.)
After further shrinking $U$ we may assume that we have a finite morphism
$\pi:V\to U$ such that $T_V\cong V\times \gm^r$. 
Let now $p\in U$ be a point with preimages $p_1,\dots, p_r\in V$
such that $k(p_i)\cong k(p)$ for every $i$ and $\pi$ is \'etale over $p$.
In the geometric cases this is possible  by (\ref{secs.isom.zeros.prop.v2.var}.1),
 for number fields this follows from the  Chebotarev density theorem.
 
Then $T_V$ has a rational section $s_V$ that has a pole along $p_1$ but regular at $p_2,\dots, p_r$.  Then
$\norm_{K/k}(s_V)$ is a rational section of $T_U$ with a pole at $p$, a contradiction.  \qed

\end{say}



 

\begin{say}[Jacobians of curves] \label{jac.say}
Let $C$ be a  proper scheme of dimension 1 over a field $k$. Then $\pico(C)$ is a called the {\it  Jacobian}  or  {\it generalized Jacobian}  of $C$  and denoted by
 $\jacs(C)$. 

Let $C^{\rm wn}\to C$ denote the {\it weak normalization} and
$\bar C\to C^{\rm wn}\to C$ the {\it normalization.}
The pull-back maps give
$$
\jacs(C)\to \jacs(C^{\rm wn})\to \jacs(\bar C).
$$
The kernel of $\jacs(C)\to \jacs(C^{\rm wn})$ is $\jacs(C)^{\rm unip}$ and
the  kernel of $\jacs(C)\to \jacs(\bar C)$ is $\jacs(C)^{\rm lin}$.
Thus  (\ref{Q.rank.say}) gives the following.
\begin{enumerate}
\item If $k$ is locally finite then $\jacs(C)(k)$ is torsion.
\item If $\chr k>0$ and $C$ is geometrically integral, then $\jacs(C)(\bar k)$ is torsion iff  $C$ is rational and $C^{\rm wn}=\bar C$. 
\item If $\chr k=0$ then $\jacs(C)(\bar k)$ is torsion iff 
$h^1(C, \o_C)=0$. This implies that 
every irreducible component of $C_{\bar k}$ is smooth and rational.
\end{enumerate}
\end{say}

\begin{defn}  \label{tr.A.defgn}
Let $P$ be a 0-cycle on $A$ and  $P_{\bar k}=\cup_i m_i[p_i]$. 
Set
$$
\tr_AP:=\tsum_i m_i[p_i] \qtq{(summation in $A$.)}
\eqno{(\ref{tr.A.defgn}.1)}
$$
Note that $\tr_AP\in A(k^{\rm ins})$. 
If $A$ is the additive group $\ga$  then this is the usual trace, but for the multiplicative group $\gm$  this is the norm. Since we usually use additive notation, trace seems a better choice. 
For $Z\subset A$ set  
$$
\tr_AZ:=\{\tr_AP: P \mbox{ is a 0-cycle on } Z\}.
\eqno{(\ref{tr.A.defgn}.2)}
$$
Let $C$ be a smooth, projective curve.
There is a natural embedding  $j:C\into \jacs_1(C)\subset \pics(C)$. 
If $P$ is a 0-cycle on $C$ then 
$$
\tr_{\pics(C)} \bigl(j(P)\bigr)=[\o_C(P)]\in \jac_{\deg P}(C).
\eqno{(\ref{tr.A.defgn}.3)}
$$
\end{defn}

The following is a restatement of  \cite[4.9]{MR3650225}.

\begin{lem} \label{sch.onto.K.onto.lem}   $A\mapsto A(k)\otimes\q$ is an exact functor on the category  of commutative algebraic groups.
\end{lem}

Proof. The only nontrivial claim is that if 
 $g:A\to B$ is  a dominant morphism then  $g(k):A(k)\to B(k)$ is surjective modulo torsion. To see this pick  $b\in B(k)$  and let
$P$ be a 0-cycle on the fiber $A_b$.  Then
$\tr_AP\in A(k)$ and  $g(\tr_AP)=\deg P\cdot b$. \qed

\medskip{\bf The multiplicative group of Artin algebras}\medskip

\begin{defn} \label{gm.not.defn}  Let $A$ be an  Artin  algebra.
The  multiplicative group of its invertible elements is denoted by $A^\times$.

If $A$ is a $k$-algebra, it is frequently convenient to view $A^\times$ as an algebraic group over $k$; tis is called {\it Weil restriction.} We denote it by $\res^A_k\gm$. That is, 
$$
\bigl(\res^A_k\gm\bigr)(B)=(A\otimes_kB)^\times
\qtq{for a $k$-algebra $B$.}
$$
Note that $\dim \res^A_k\gm=\dim_k A$. 

For example, if $K/k$ is a field extension of degree $n$, choose a basis $e_i\in K$. As a variety, $\res^K_k\gm$ is 
$\a^n\setminus\bigl(\norm_{K/k}(\tsum x_ie_i)=0\bigr)$.
\end{defn}

\begin{say}\label{sep.residue.field.lem}
Let $(A, m)$ be a semi-local, Artin $k$ algebra and $K=A/m$. 
There is an exact sequence
$$
1\to U\to A^\times \to K^\times\to 1,
\eqno{(\ref{sep.residue.field.lem}.1)}
$$
where the map  $a\mapsto 1+a$ identifies $m$ with $U$. Note that 
$a\mapsto 1+a$ is a group isomorphism if $m^2=0$ but not otherwise.
In characteristic 0 one can correct this by taking
$a\mapsto \exp(a)$. We will think of (\ref{sep.residue.field.lem}.1) as the
$k$-points of an exact sequence of algebraic $k$-groups
 $$
1\to U\to \res^A_k\gm \to \res^K_k\gm\to 1,
\eqno{(\ref{sep.residue.field.lem}.2)}
$$
where $U$ is a unipotent group.
In positive characteristic the algebraic groups
$(m, +)$  and $ (U, \cdot)$ need not be isomorphic.

We also use the following variant of Hensel's lemma.

\medskip
{\it Claim \ref{sep.residue.field.lem}.3.}
Let  $k\subset k'\subset K$ be 
a subfield that is separable over $k$. Then there is a unique lifting
$k'\to A$. \qed
\end{say} 

Combining the above with the previous discussion on algebraic groups yield the following lemmas. 

\begin{lem} \label{ker.fin.Qrank.lem.2}
Let $k$ be a field and $A\to B$ a homomorphism of Artin $k$-algebras.  Then  $\coker[A^{\times}\to B^{\times}]$ is torsion iff one of the following holds. 
\begin{enumerate}
\item $A\to B$ is surjective.
\item $\chr k>0$ and  $B/\sqrt{0}$ is purely inseparable over $A/\sqrt{0}$.
\item $k$ is locally finite. 
\end{enumerate}
Moreover, $\coker[A^{\times}\to B^{\times}]$ has finite $\q$-rank in one additional case:
\begin{enumerate}\setcounter{enumi}{3}
 \item $\deg (k/\q)<\infty$ and  $A/\sqrt{0}\to B/\sqrt{0}$ is surjective. \qed
\end{enumerate}
\end{lem}

\begin{lem} \label{ker.fin.Qrank.lem}
Let $k$ be a field and $A\to B$ a homomorphism of Artin $k$-algebras.  Then  $\ker[A^{\times}\to B^{\times}]$ is torsion iff one of the following holds. 
\begin{enumerate}
\item $A\to B$ is  injective.
\item $\chr k>0$ and $A/\sqrt{0}\to B/\sqrt{0}$ is injective.
\item $k$ is locally finite. 
\end{enumerate}
Moreover, $\ker[A^{\times}\to B^{\times}]$ has finite $\q$-rank in one additional case:
\begin{enumerate}\setcounter{enumi}{3}
 \item $\deg (k/\q)<\infty$ and  $A/\sqrt{0}\to B/\sqrt{0}$ is injective. \qed
\end{enumerate}
\end{lem}

\begin{lem} \label{2.fields.mult.gen.lem}
Let $k$ be a field that is not locally finite. 
Let  $A\supset k$ be a finite, reduced $k$-algebra. 
Let  $k\subset L_1, L_2\subset A$ be subfields.  The following are equivalent. 
\begin{enumerate}
\item $A^\times/(L_1^\times\cdot L_2^\times)$ is torsion.
\item $A^\times/(L_1^\times\cdot L_2^\times)$ has finite $\q$-rank.
\item $A$ is a field and $A/L_i$  is purely inseparable for some $i=1,2$. 
\end{enumerate}
\end{lem}

Proof.  If  $A/L_i$  is purely inseparable then $A^q\subset L_i$ for some
power $q$ of $\chr k$, hence $A^\times/L_i^\times$ is torsion. This proves
(\ref{2.fields.mult.gen.lem}.3) $\Rightarrow$ (\ref{2.fields.mult.gen.lem}.1) and  (\ref{2.fields.mult.gen.lem}.1) $\Rightarrow$ (\ref{2.fields.mult.gen.lem}.2) is clear. 

Assume (\ref{2.fields.mult.gen.lem}.2). We may replace $A$ by its maximal separable subalgebra.
Thus assume that $A/k$ is separable.  If $A=L_i$ for some $i$ then we are done.
Otherwise
$\dim_k A=\dim_{L_i}A\cdot \dim_kL_i\geq 2\dim_kL_i$. 

If $B$ is a reduced, separable $k$-algebra then  $B^\times$ is identified with the
$k$-points of the $k$-torus  $\res^B_k\gm$. Thus
$L_1^\times\cdot L_2^\times\to A^\times$ can be viewed as the $k$-points of a morphism of $k$-tori
$$
\mu:  \res^{L_1}_k\gm\times \res^{L_2}_k\gm
{\longrightarrow} \res^{A}_k\gm.
$$
Both of the $L_i$ contain $k$, thus
$$
\dim_k \im(\mu)\leq \dim_k L_1+\dim_k L_2-1<\dim A.
$$
Thus $\coker(\mu)$ is a positive dimensional $k$-torus, hence
$ \qrank \bigl(\coker(\mu)(k)\bigr)=\infty$ by (\ref{Q.rank.say}.4). 
Finally (\ref{sch.onto.K.onto.lem}) shows that
$$
\qrank \bigl(A^\times/(L_1^\times\cdot L_2^\times)\bigr)=\qrank \bigl(\coker(\mu)(k)\bigr)=\infty. \qed
$$

\section{Special fields}\label{spec.fields.sec}

We discuss various classes of fields that were used earlier.

\begin{say}[Locally finite fields]  \label{loc.fin.defn} 
A field $k$ is called {\it locally finite}\footnote{The terminology is not standard in English; it is an  analog of the notion of locally finite group.} if the following equivalent conditions hold.
\begin{enumerate}
\item Every finitely generated subfield of $k$ is finite.
\item $k$ is an algebraic extension of a finite field.
\item $k$ is isomorphic to a subfield of $\bar \f_p$ for some $p>0$.
\item $A(k)$ is a torsion group for every Abelian variety $A$ over $k$.
\item There is a  $C>0$ such that $\qrank A(k)\leq C$ for
 every Abelian variety $A$ over $k$.  
\end{enumerate}
The only non-obvious claims are (\ref{loc.fin.defn}.4--5). 
If $k$ is not locally finite then it
contains either  $\q$ or $\f_p(t)$. In both cases, there is an 
Abelian variety $A$ over $k$ with arbitrarily large 
$\rank_{\q} A(k)$. For example, if $E$ is an elliptic curve of  rank $\geq 1$  then $E^m$ has rank $\geq m$.
\end{say}

\begin{say}[$\q$--Mordell-Weil fields]  \label{m.w.f.defn}
A field $k$ is  {\it Mordell-Weil}  (resp.\ {\it $\q$--Mordell-Weil}) if for every Abelian variety $A$ over $k$, the group of its $k$-points 
$A(k)$ is finitely generated  (resp.\ has finite $\q$-rank). 

By \cite{MR0102520}, every  finitely generated  field is 
Mordell-Weil.

Weil restriction (cf.\ \cite[Sec.7.6]{blr})  shows that these properties are invariant under finite field extensions.
Since every Abelian variety is a quotient of a Jacobian, 
it is equivalent to ask that $\jac(C)$ have these properties  for every smooth projective curve $C$ over $k$ (\ref{sch.onto.K.onto.lem}).

It is nor clear how much the class of $\q$--Mordell-Weil fields differs from the class of Mordell-Weil fields.
For example, $\bar \f_p$ is $\q$--Mordell-Weil but not Mordell-Weil.
\end{say}

\begin{say}[Anti--Mordell-Weil fields] \label{anti-mordell.defn} Following \cite{2019arXiv190204011I} a   field $k$ is called {\it anti--Mordell-Weil} if 
\begin{enumerate}
\item the $\q$-rank of
$A(k)$ is infinite  for every
positive dimensional abelian variety. 
\end{enumerate}
In particular, $k$ is not locally  finite. 
If the latter holds then  the $\q$-rank of
$T(k)$ is infinite  for every  $k$-torus  $T$ (\ref{Q.rank.say}.4), hence (\ref{anti-mordell.defn}.1) can be restated as: 
\begin{enumerate}\setcounter{enumi}{1}
\item The $\q$-rank of
$A(k)$ is infinite  for every
positive dimensional semi-abelian variety $A$. 
\end{enumerate}
Finally, if $\chr k=0$, then  $k$ is not a finite extension of $\q$ by the Mordell-Weil theorem, hence the  $\q$-rank of
$U(k)$ is infinite  for every unipotent group (\ref{Q.rank.say}.3). Thus  (\ref{anti-mordell.defn}.1--2) are further equivalent to:
\begin{enumerate}\setcounter{enumi}{2}
\item The $\q$-rank of
$A(k)$ is infinite  for every
positive dimensional commutative algebraic group  $A$. 
\end{enumerate}
{\it Warning.} Note that if $\chr k=p>0$ then $U(k)$ is $p$-power torsion  for every unipotent group; this creates a crucial difference between 0 and positive characteristics for us.

Examples of anti--Mordell-Weil fields are the following.
\begin{enumerate}\setcounter{enumi}{3}
\item algebraically closed fields, save for $\bar \f_p$,
\item $\r$ and all real closed fields,
\item $\q_p$, more generally quotient fields of  Henselian, local domains,
\item large fields    \cite{MR2284973, MR2652895, 2019arXiv191210710F} that are not locally finite, where a field $k$ is {\it large} (also called ample, fertile  or anti-mordell) if  $C(k)$ is either empty or infinite for every smooth curve $C$.
\end{enumerate}
The last case implies  the earlier ones.

 \end{say}

\begin{say}[Hilbertian fields]\label{hilb.field.say}
 A field $k$ is {\it Hilbertian} if
for every irreducible polynomial  $f(x,y)\in k[x,y]$ such that
$\partial f/\partial y\neq 0$, there are  infinitely many $c\in k$ such that $f(x,c)\in k[x]$ is irreducible. (We follow \cite[Chap.12]{fri-jar} with adding the separability condition.)

Equivalently, for every smooth, irreducible curve $C$ and every 
basepoint-free linear system $|M|$ that defines a separable map $C\to \p^1$, there are 
 infinitely many  irreducible members $M_c\in |M|$. 
This also implies that,  for every irreducible variety $X$ over $k$, and every 
mobile linear system $|M|$  that defines a separable map $X\map \p^N$, there is a dense set
$\Lambda\subset |M|(k)$ such that $M_\lambda\in |M|$
is  irreducible for $\lambda\in \Lambda$.

Hilbert proved that number fields are  Hilbertian. More generally, 
 every finitely generated, infinite field is Hilbertian.
A finite, separable extension of a Hilbertian field is Hilbertian, and so is any purely inseparable extension. 
 See \cite[Chap.VIII]{MR0142550} or \cite[Chaps.12--13]{fri-jar} for these and many other facts.  
\end{say}

In our proofs the Hilbertian condition is mostly used through the following consequence.

\begin{lem}\label{hilb.f.basic.prop.lem.2}  Let $C$ be an irreducible,  geometrically reduced,   projective curve  over a Hilbertian  field $k$. Let $\Sigma\subset C$ be a finite subset and  $Z\subset C$ a finite subscheme. Let $L$ be a  line bundle on $C$  such that $\deg L\geq \deg Z+\deg \omega_C+3$.   
Then every $s_Z\in H^0(Z, L|_Z)^\times$ can be lifted to   $s_C\in H^0(C, L)$
such that $(s_C=0)$ is irreducible, reduced  and disjoint from $\Sigma$. 
\end{lem}

Proof. The condition $s|_Z=c\cdot s_Z$ for some $c\in k$ determines a linear subsystem  $|L, s_Z|\subset |L|$. The degree condition guarantees that it is basepoint-free and separable. Hence it has infinitely many irreducible members. Almost all of them are disjoint from $\Sigma$. \qed

\medskip
It turns out that  versions of  (\ref{hilb.f.basic.prop.lem.2}) hold for some  non-Hilbertian fields and 
in our  proofs  a  weakening of 
it  is sufficient.   In order to state it, we need a definition.

\begin{defn}\label{hilb.f.basic.prop.defn} Let $C$ be an irreducible,  geometrically reduced,  projective curve over a field $k$. Let $\Sigma\subset C$ be a finite subset, $Z\subset C$ a finite subscheme and
 $L$ an ample line bundle on $C$.  Let
 $$
\Gamma^{\rm irr}(C,L, Z, \Sigma^c)\subset  
 H^0\bigl(Z, \oplus_{m> 0}L^m|_Z\bigr)^\times
\eqno{(\ref{hilb.f.basic.prop.defn}.1)}
$$
be the subset consisting of  the restrictions $s|_Z$ of those sections 
$s\in H^0(C, L^m)$, for which   $(s=0)$ is irreducible 
 and disjoint from $\Sigma\cup Z$.

If $k$ is Hilbertian then,  by
(\ref{hilb.f.basic.prop.lem.2}), 
$$
\Gamma^{\rm irr}(C,L, Z, \Sigma^c)\supset  
 H^0\bigl(Z, \oplus_{m\geq m_0}L^m|_Z\bigr)^\times
\qtq{for some $m_0$.}
\eqno{(\ref{hilb.f.basic.prop.defn}.2)}
$$



\end{defn}


\begin{defn}[Weakly Hilbertian fields]\label{weak.hilb.say}
Let $C$ be a projective curve over  a field $k$ and $L$ an ample line bundle on $C$. We say that $(C, L)$ is {\it weakly Hilbertian}
if the following holds.
\begin{enumerate}
\item  For every    $s_Z\in H^0\bigl(Z, L^m|_Z\bigr)^\times$,
some power of $s_Z$ is in $\Gamma^{\rm irr}(C,L, Z, \Sigma^c)$. 
\end{enumerate}
  We call a field $k$ 
{\it weakly Hilbertian} if $(C, L)$ is  weakly Hilbertian
for every  irreducible,  geometrically reduced,  projective curve  over  $k$ and every ample line bundle $L$ on $C$.
We see in (\ref{wH.goes.down.lem.0}) that it is enough to check this for 
smooth curves.

By (\ref{hilb.f.basic.prop.lem.2})
every  Hilbertian field  is also  weakly Hilbertian. 
Below we show that all   Galois extensions  $\q\subset K\subsetneq \bar\q$ are  weakly Hilbertian. More generally,  Galois extensions of Hilbertian fields, save those that are separably closed, are  weakly Hilbertian; see (\ref{wH.goes.down.cor.1}). The fields $\q_p$ and $\f_p((t))$  are not weakly Hilbertian, but all smooth curves over them satisfy property (\ref{weak.hilb.say}.2) below; see (\ref{hens.wH.res.f.lf.thm}). 
See (\ref{poonen.letter.say}) for finite fields.

{\it Comment.} We do not claim that the above  is the `right' notion of weakly Hilbertian fields. It works well for our purposes. However, from the  field-theory point of view, the following even weaker version may seem more natural.
\begin{enumerate}\setcounter{enumi}{1}
\item Let $C$ be an irreducible,  geometrically reduced,  projective curve over $k$ and  $L$ an ample line bundle on $C$. Then there is an $m>0$ and infinitely many $s_i\in H^0(C, L^{m_i})$ whose zero sets  $(s_i=0)$ are irreducible and disjoint. 
\end{enumerate}
It is clear that  (\ref{weak.hilb.say}.1) $\Rightarrow$ (\ref{weak.hilb.say}.2); we do not know whether the converse holds.
\end{defn}

\medskip
{\bf Subfields of weakly Hilbertian fields}
\medskip

Let $K/k$ be a  finite, separable field extension. If $k$ is Hilbertian, then so is $K$, but the converse fails. We prove in
(\ref{wH.goes.down.lem}) that if $K$  is weakly Hilbertian, then so is $k$.
We do now know whether the converse holds.

\begin{lem}\label{wH.goes.down.lem.0}
Let  $C$ be an irreducible,  geometrically reduced,  projective curve  and $L$ an  ample line bundle  on $C$. Let $\pi:\bar C\to C$ denote the normalization.
If $(\bar C, \bar L)$ is  weakly Hilbertian, then so is $(C, L)$.
\end{lem}

Proof. 
We are free to enlarge $Z$, hence we may assume that 
 it contains the conductor subscheme of $\pi$ and 
its ideal sheaf is also an ideal sheaf on $\bar C$, defining $\bar Z$.  We have a natural inclusion  $\tau: \o_{Z}\into \pi_*\o_{\bar Z}$. Thus if $\bar s$ is a section of $\bar L^m$
whose restriction to $\bar Z$ agrees with $\tau(s_Z^m)$, then
it descends to a section $s$ of $L^m$. \qed


\begin{lem}\label{wH.goes.down.lem}
 Let $K/k$ be a  separable, algebraic field extension. If $K$ 
is weakly Hilbertian  then so is $k$.
\end{lem}

Proof. Start with  $C_k, L_k, Z_k, \Sigma_k$ and  $s_{Z,k}\in H^0(Z_k, L_k|_{Z_k})^\times$. 
By base change we get $C_K, L_K, Z_K$ and   $s_{Z,K}\in H^0(Z_K, L_K|_{Z_K})^\times$.  Let $C'_K\subset C_K$ be one of its irreducible components.

Assume that we have $s'_K\in H^0(C'_K, L'_K)$ such that   $s'_K|_{Z'_K}=s_{Z,K}^m$.
This $s'_K$ is defined over a finite degree subextension
$k\subset K_1\subset K$.

Next assume first that  $C_k$ is smooth. Then $C'_K$ is smooth and
$C'_K\to C_k$ is flat. Thus $\norm_{K_1/k}$ sends sections of
$(L'_K)^m$ to sections of
$L_k^{md}$  where $d=\deg (K_1/k)$. 
Set   $s_k:=\norm_{K_1/k}(s'_K)$.  Then $s_k\in H^0(C_k, L_k^{md})$ and
$s_k|_{Z_k}=s_{Z,k}^{md}$.

The singular case follows from (\ref{wH.goes.down.lem.0}).
\qed

\begin{cor}\label{wH.goes.down.cor.1}  Let $k$ be a  Hilbertian field and 
 $K/k$  a Galois extension that is not    separably closed. Then $K$ is 
 weakly Hilbertian.
\end{cor}

Proof. By \cite{MR667458},  every nontrivial finite extension of such a field $K$  is  Hilbertian.  Thus, if  $K$ is not separably closed, then   
$K$  is  weakly Hilbertian  by (\ref{wH.goes.down.lem}). 
\qed

\medskip
As an  example,  let $\q^{\rm solv} $ denote the composite of all  Galois extensions
of $\q$ with solvable Galois group. It is   weakly Hilbertian by
(\ref{wH.goes.down.cor.1}), but  not Hilbertian, as shown by the polynomial $x^2-y$.

\begin{cor}\label{wH.goes.down.cor}  Let $K\subset \bar \q$ be a subfield whose Galois closure is different from $\bar\q$. Then $K$ is   weakly Hilbertian.  \qed
\end{cor}

\medskip
{\bf Not weakly Hilbertian fields}
\medskip

Next we  show  that $\bar\q$ is not weakly Hilbertian, 
it does not even satisfy   (\ref{weak.hilb.say}.2) for smooth curves.
 Note that  a much stronger variant of (\ref{barQ.not.wh.cor}) could be true; 
see (\ref{few.zeros.conj.2}) and (\ref{barQ.not.wh.cor.cor}).
The same methods work for  $\r\cap \bar\q$. We do not know any other subfield of $\bar\q$ that does not  satisfy   (\ref{weak.hilb.say}.2) for smooth curves, though presumably there are many.

\begin{prop} \label{barQ.not.wh.cor}  There is a smooth projective curve $C$ and an ample line bundle $L$, defined over $\bar\q$, such that  every section of $L^m$ has at least 2 distinct zeros  for every $m>0$.
\end{prop}

Proof.  Let $\pi:C\to B$  be a nonconstant morphism between smooth, projective curves
such that $g(C)\geq g(B)+2$ and $g(B)\geq 1$.   Let $\Gamma\subset \pic(B)$ be as in (\ref{2.zeroes.lem.1}). 

The $\q$-rank of $\pic(B)$ is infinite by (\ref{anti-mordell.defn}.4),  so 
 there is an ample 
$L\in \pic(B)$ no power of which is  in $\Gamma$.  
Then  $\pi^*L$ has the required property. \qed
\medskip

\begin{lem} \label{2.zeroes.lem.1}
Let $\pi:C\to B$ be a nonconstant morphism between smooth, projective curves
such that $g(C)\geq g(B)+2$.  Then
$$
\Gamma:=\bigl\langle L\in \pic(B): \pi^*L^m\sim \o_C\bigl(n[c]\bigr)
\mbox{ for some } c\in C(k), n, m>0\bigr\rangle\subset \pic(B)
$$
has finite $\q$-rank.
\end{lem}

Proof. This is clear if $C(k)$ is empty. Otherwise fix a point $c_0\in C(k)$
and embed $C\into \jacs(C)$ sending $c_0$ to the origin.
Set $A:=\jacs(C)/\pi^*\jacs(B)$ with quotient map
$\sigma:C\to A$. If $\pi^*L^m\sim \o_C\bigl(n[c]\bigr)$ then
$\sigma(c)\in A$ is a torsion point. Thus there are only finitely many such $c\in C(k)$ by \cite{MR1609518}, and the $\q$-rank of $\Gamma$ is at most the number of such torsion points.  \qed

\medskip

For nodal rational curves, there is an elementary proof.

\begin{prop} \label{barQ.not.wh.cor.prop} Nodal rational curves over $\bar\q$ are not weakly Hilbertian.
\end{prop}

Proof. A rational curve with $r$ nodes is obtained form $\p^1$ by identifying
$r$ point pairs. Thus we start with  $2r$ distinct points
$a_1,\dots, a_{2r}\in \a^1$  and identify $a_i$ with $a_{r+i}$ to get a nodal rational curve $C$. 

A line bundle on $C$ is obtained by starting with some $L=\o_{\p^1}(m)$
and specifying isomorphisms  $L|_{a_i}\cong L|_{a_{r+i}}$. 
Thus sections of the resulting line bundle are given by
polynomials $p(x)$ of degree $\leq m$ such that 
$p(a_i)=u_ip(a_{r+1})$ for every $i$ where  $u_i\in k^\times$ specify the line bundle.

 A polynomial with  zeros $\{z_j:j\in J\}$ is 
$s(x)=\gamma \tprod_j (x-z_j)^{m_j}$. Thus
for nonzero  $u_1,\dots,  u_r$ we aim to solve the  $r$ equations
$$
\prod_{j\in J} \Bigl(\frac{a_i-z_j}{a_{r+i}-z_j}\Bigr)^{m_{j}}
=u_i^n,
\eqno{(\ref{barQ.not.wh.cor.prop}.1)}
$$
where $m_j, n\in \z$ and $z_j\in \bar\q$ are unknowns with  $n\neq 0$. 

For every $p$ choose an extension $v_p$ of the  $p$-adic valuation to $\bar\q$. 
The $v_p$-valuation of  any $d\in \bar\q^\times$ is 0 for all but finitely many $p$. We can thus choose $p$ such that 
$$
v_p(a_i)=v_p(a_i-a_j)=0 \qtq{for every} i\neq j.
\eqno{(\ref{barQ.not.wh.cor.prop}.2)}
$$
 Thus taking the valuation of (\ref{barQ.not.wh.cor.prop}.1) we get the equations
$$
\sum_j m_j v_p \Bigl(\frac{a_i-z_j}{a_{r+i}-z_j}\Bigr)=n v_p(u_i).
\eqno{(\ref{barQ.not.wh.cor.prop}.3)}
$$ 
Choose the $u_i$ such that $v_p(u_i)\neq 0$ for every $i$.
By (\ref{barQ.not.wh.cor.prop}.5), for every $i$ we get a
$\sigma(i)\in J$ such that
$$
v(a_i-z_{\sigma(i)})>0 \qtq{or} v(a_{r+i}-z_{\sigma(i)})>0.
\eqno{(\ref{barQ.not.wh.cor.prop}.4)}
$$
If $r>|J|$  then the same $z_j$ appears twice. 
Thus we have  $v(a_{i_1}-z_j)>0$ and $v(a_{i_2}-z_j)>0$ for some $i_1\neq i_2$ and $j$.  Then 
$$
v(a_{i_1}-a_{i_2})\geq \min\{v(a_{i_1}-z_j), v(a_{i_2}-z_j)\}>0
$$
gives a contradiction.\qed

\medskip
{\it Claim \ref{barQ.not.wh.cor.prop}.5.}  Let 
$(R, v)$ be a valuation ring and 
$a,c\in R^\times$ such that  $v(a)=v(c)=v(a-c)=0$. Then
$$
v\Bigl(\frac{a-z}{c-z}\Bigr)>0 \ \Leftrightarrow\  v(a-z)>0\qtq{and}
v\Bigl(\frac{a-z}{c-z}\Bigr)<0 \ \Leftrightarrow\  v(c-z)>0. \qed
$$

The above proof shows the following stronger claim.

\begin{cor} \label{barQ.not.wh.cor.cor}  Let $C$ be a rational curve with $g$ nodes  over $\bar\q$. Then for `most' ample line bundles $L$ over $C$, 
every section of $L^m$ has at least $g$ distinct zeros. \qed
\end{cor}



\medskip
{\bf Quotient fields of discrete valuation rings}
\medskip

$\q_p$, $\f_p((t))$, and, more generally, quotient fields of excellent DVRs, have very interesting behavior.  Smooth curves have weak Hilbertian properties but singular curves do not. 

Note that a DVR   of  characteristic 0 is excellent   \cite[07QW]{stacks-project}. In positive characteristic,  local rings of smooth curves are excellent and so are power series rings $K[[t]]$. However, 
there are many non-excellent DVRs; see
\cite{dat-smi} for especially simple examples.
(The proof below uses excellence, but the result might hold for any DVR.)

\begin{prop} \label{hens.wH.res.f.lf.thm}
Let  $(R,m)$ be an excellent  DVR with quotient field $K$ and  locally finite residue field $k$.   Let $C_K$ be a smooth projective, irreducible curve over $K$ and $L_K$ an ample line bundle on $C_K$. Then $|L_K^n|$ has infinitely  many irreducible members for  $n$ sufficiently divisible.
\end{prop}

Proof. We extend $C_K$  to a flat morphism $C_R\to \spec R$. We may assume that $C_R$ is regular by \cite{shaf-MR0217068}. Then $L_K$ extends to a line bundle $L_R$ on $C_R$. 

 Let $E_1,\dots, E_r$ be the irreducible components of the central fiber  $C_k=\tsum m_i E_i$.  The intersection matrix
$(E_i\cdot E_j)$ is negative semidefinite with 
$[C_k]$ as the only null-vector. Thus the intersection matrix of
$E_2,\dots, E_r$ is negative definite. We can thus find a divisor
$F$ supported on $E_2,\dots, E_r$ such that  $L^{n_1}_R(F)$ has degree 0
on $E_2,\dots, E_r$. 

By \cite{Artin62}  the curves $E_2,\dots, E_r$ can be contracted 
$C_R\to C_R^*$ and
a suitable power of $L^{n_1}_R(F)$ descends to a line bundle $L^*$ on $C_R^*$.

Now we have a normal scheme with a flat morphism
$\pi:C_R^*\to R$ whose generic fiber is $C_K$ and whose central fiber  $C^*_k$ is an irreducible curve.    Furthermore there is a line bundle $L^*$
whose restriction to $C_K$ is a power of $L_K$.

Set  $E^*=\red C^*_k$ and pick any point $p\in E^*$  that is regular both on $E^*$ and on $C^*_R$.  Since $\qrank\pic(E^*)=1$, after passing to 
a power of $L^*$, we may assume that 
\begin{enumerate}
\item $L^*|_{E^*}$ has a section $\bar s_0$ that  vanishes only at $p$,
\item $L^*(-E^*)/ L^*(-2E^*)$ has a section  $\bar s_1$ that does not vanish at $p$, and 
\item $H^1\bigl(C^*_R, L^*(-E^*)\bigr)=H^1\bigl(C^*_R, L^*(-2E^*)\bigr)=0$.
\end{enumerate}
By (\ref{hens.wH.res.f.lf.thm}.3) we can lift $\bar s_0$ and $\bar s_1$  to
$s_0\in H^0\bigl(C^*_R, L^*\bigr)$ and $s_1\in H^0\bigl(C^*_R, L^*(-E^*)\bigr)$.
For all but 1 residue value of $\lambda\in R$, 
$D_R(\lambda):=(s_0+\lambda s_1=0)$ is regular at $p$. Since $p$ is its sole point over $k$, $D_R(\lambda)$  is irreducible and reduced. \qed

\begin{prop} \label{hens.dvr.nodal.notwh} Let  $(R,m)$ be a Henselian  DVR with quotient field $K$ and   residue field $k$.  
For $2g\leq |k|$ 
there are  rational curves $C$  with $g$ nodes  over $K$ such that, for `most' ample line bundles $L$ over $C$, 
every section of $L^n$ has at least $g$ distinct zeros. \qed
\end{prop}

Proof. We choose  $a_i\in R\subset \a^1(K)\subset \p^1_K$ such that  $\bar a_i\in k$ (their reduction mod $m$)  are all distinct. 
As in (\ref{barQ.not.wh.cor.prop}), 
identifying  $a_i$ with $a_{g+i}$ to get a nodal rational curve $C$ and  
a line bundle $L$ on $C$ is obtained by starting with some $\o_{\p^1}(r)$
and specifying isomorphisms  $\o_{\p^1}(r)|_{a_i}\cong \o_{\p^1}(r)|_{a_{g+i}}$. 
Thus sections of $L^n$  are given by
polynomials $f(x)$ of degree $\leq nr$ such that 
$$
f(a_i)=u_i^nf(a_{g+1})\qtq{for} i=1,\dots, g,
\eqno{(\ref{barQ.not.wh.cor.prop}.1)}
$$
where the $u_i\in K^\times$  determine $L$.
We may assume that $f(x)\in R[x]\setminus m[x]$ and denote by  $\bar f$ its image in $k[x]$.
If the $m$-adic valuation of $u_i$ is not 0 for every $i$,  this implies that  
$$
 \bar f(\bar a_i)=0 \qtq{or} \bar f(\bar a_{g+i})=0
\qtq{for} i=1,\dots, g.
 \eqno{(\ref{barQ.not.wh.cor.prop}.1)}
$$
Thus  $\bar f$ has at least $g$ distinct zeros, so, by Hensel's lemma,
$f$ has at least  $g$ distinct prime factors.  \qed
\medskip

Note that, in the above case, the precise meaning of `most' line bundles is the following.
Choosing the point at infinity fixes an isomorphism
$\pic(C)\cong \pic^\circ(C)\times \z$.

For a  node $p_i$ the map $\rho_i:\pic^\circ(C)\to \mg$ depends on the ordering of the preimages of the node. So we get 2 maps, that differ only by composing with the inverse map of $\mg$.  Composing $\rho_i$ with the absolute value of the $m$-adic valuation gives a well-defined map
$v_i:\pic^\circ(C)\to \n$. The $g$ nodes together give
$v:\pic(C)\to \n^g$.  If all coordinates of $v(L)$ are nonzero, then
every section of $L^n$ has at least $g$ distinct zeros for every $n$.

\begin{say}\label{poonen.letter.say}
B.~Poonen explained to us that, using geometric class field
theory
and the function field Chebotarev density theorem, one can prove that if $C$ is a smooth projective curve over a finite field $k$ and $L$ is an ample line bundle, then $L^m$ has a section with irreducible and reduced zero set for all $m\gg 1$.  However, the probability that a random section 
has this property tends to 0 as $m \to\infty$.
\end{say}

\def\cprime{$'$} \def\cprime{$'$} \def\cprime{$'$} \def\cprime{$'$}
  \def\cprime{$'$} \def\dbar{\leavevmode\hbox to 0pt{\hskip.2ex
  \accent"16\hss}d} \def\cprime{$'$} \def\cprime{$'$}
  \def\polhk#1{\setbox0=\hbox{#1}{\ooalign{\hidewidth
  \lower1.5ex\hbox{`}\hidewidth\crcr\unhbox0}}} \def\cprime{$'$}
  \def\cprime{$'$} \def\cprime{$'$} \def\cprime{$'$}
  \def\polhk#1{\setbox0=\hbox{#1}{\ooalign{\hidewidth
  \lower1.5ex\hbox{`}\hidewidth\crcr\unhbox0}}} \def\cdprime{$''$}
  \def\cprime{$'$} \def\cprime{$'$} \def\cprime{$'$} \def\cprime{$'$}
\providecommand{\bysame}{\leavevmode\hbox to3em{\hrulefill}\thinspace}
\providecommand{\MR}{\relax\ifhmode\unskip\space\fi MR }
\providecommand{\MRhref}[2]{%
  \href{http://www.ams.org/mathscinet-getitem?mr=#1}{#2}
}
\providecommand{\href}[2]{#2}

\bigskip

  Princeton University, Princeton NJ 08544-1000

\email{kollar@math.princeton.edu}

\end{document}